\documentclass[12pt,letterpaper]{article}

\usepackage{natbib}
\setcitestyle{authoryear,open={(},close={)},citesep={,}}
\usepackage{hyperref}
\usepackage{xcolor}
\hypersetup{colorlinks=true,citecolor=blue,linkcolor=blue,urlcolor=blue}
\usepackage{amsmath,amsthm,amssymb,mathtools,mathrsfs}
\usepackage{bm,graphicx,booktabs,multirow,dsfont}
\usepackage{cleveref}
\usepackage{enumitem}
\usepackage{placeins}
\usepackage[a4paper,left=1.75cm,right=1.75cm,top=2.5cm,bottom=2.5cm]{geometry}

\usepackage{caption}
\newtheorem{theorem}{Theorem}[section]
\newtheorem{lemma}{Lemma}[section]

\newtheorem{proposition}{Proposition}[section]

\newtheorem{example}{Example}[section]
\newtheorem{remark}{Remark}[section]
\newtheorem{condition}{Condition}[section]
\newenvironment{customcon}[1]{\renewcommand\thecondition{#1}\condition}{\endcondition}

\DeclareMathOperator*{\argminop}{arg\,min}

\renewcommand{\P}{\mathbb P}
\newcommand{\E}{\mathbb E}
\newcommand{\R}{\mathbb R}
\newcommand{\M}{\mathbb M}
\newcommand{\Sph}{\mathbb S}

\newcommand{\cH}{\mathcal H}
\newcommand{\cL}{\mathcal L}

\newcommand{\cT}{\mathcal T}
\newcommand{\cN}{\mathcal N}
\newcommand{\tran}{^{\mathstrut\scriptscriptstyle\top}}
\newcommand{\mx}{\mathbf x}
\newcommand{\mz}{\mathbf z}
\newcommand{\mb}{\mathbf b}
\newcommand{\mX}{\mathbf X}
\newcommand{\mh}{\mathbf h}
\newcommand{\md}{\mathbf d}

\newcommand{\zero}{\mathbf 0}

\newcommand{\mI}{\mathbf I}
\newcommand{\mB}{\mathbf B}
\newcommand{\mH}{\mathbf H}

\newcommand{\mQ}{\mathbf Q}
\newcommand{\mP}{\mathbf P}

\newcommand{\cK}{\mathcal K}

\newcommand{\mZ}{\mathbf Z}
\newcommand{\mGamma}{\bm\Gamma}
\newcommand{\dd}{\,\mathrm d}

\newcommand{\cF}{\mathcal F}

\newcommand{\1}{\mathds 1}

\renewcommand{\theequation}{\thesection.\arabic{equation}}
\renewcommand{\theHequation}{\thesection.\arabic{equation}}

\begin{document}

\title{Local Fr\'echet regression with polyspherical predictors:\\ Uniform convergence rates and factor selection} \author{\large Chang Jun Im$^1$} \date{\normalsize $^1$The Institute for Data Innovation in Science, Seoul National University} \maketitle

\begin{abstract}
We study local Fr\'echet regression for responses taking values in a metric space and predictors on a product of spheres of possibly different dimensions. We propose local constant and local linear estimators with factor-specific bandwidths and establish uniform convergence rates. We select sphere factors by comparing cross-fitted prediction risks for the full and leave-one-factor-out models. The procedure consistently recovers the active factor set, and the selected local constant and local linear estimators attain the same uniform convergence rates as if the active factor set were known. For univariate Wasserstein responses supported on a common compact interval, the uniform rate holds without an entropy condition on the response space. Simulations with scalar and spherical responses and analyses of vessel speeds and drifter speed distributions illustrate the methods.
\end{abstract}

\small
\begin{quotation}
\end{quotation}
\normalsize

\section{Introduction}\label{sec:introduction}

Directional and periodic data arise in studies of movement, spatial orientation and temporal cycles; see \citet{MardiaJupp2000}. Geographic location can be represented by a point on a two-sphere, whereas a heading or periodic time variable can be represented by a point on a circle. Nonparametric methods for directional data include kernel density estimation on spheres \citep{HallWatsonCabrera1987} and local polynomial regression for toroidal predictors with scalar responses \citep{DiMarzioPanzeraTaylor2009} and spherical predictors with scalar and spherical responses \citep{DiMarzioPanzeraTaylor2014}. Joint observations lie on a product of spheres, or a polysphere, with possibly different factor dimensions. \citet{GarciaPortuguesMeilanVila2026} studied kernel density estimation on polyspheres.

Beyond scalar and spherical responses, regression with directional predictors may involve probability distributions, covariance matrices, networks, shapes and phylogenetic trees as responses. Such responses may take values in metric spaces without a natural vector space structure, so conventional regression based on conditional expectations may not be directly applicable. \citet{PetersenMuller2019} introduced Fr\'echet regression for responses in a general metric space with Euclidean predictors, defining the regression target by minimizing conditional expected squared distance. \citet{ChenMuller2022} established uniform convergence results for local Fr\'echet regression in that setting. For directional predictors, \citet{ImJeonPark2025} developed local constant and local linear estimators on a sphere, and \citet{ImJeon2026} extended the construction to a torus. \citet{ImJeon2026Riemannian} considered general finite-dimensional Riemannian predictors. Although a polysphere is a Riemannian manifold, its product structure permits separate bandwidths and selection of entire sphere factors.

When several sphere factors are included as predictors, irrelevant factors increase the dimension of nonparametric smoothing. \citet{JeonParkVanKeilegom2021} studied additive regression for non-Euclidean predictors and responses taking values in a Hilbert space, reducing model complexity through an additive structure. For Euclidean predictors, \citet{TuckerWuMuller2023} studied variable selection for global Fr\'echet regression, and \citet{DongEtAl2026} studied variable selection for Nadaraya--Watson--Fr\'echet regression. \citet{YangEtAl2026} studied variable selection for additive global Fr\'echet regression. These approaches do not directly address the selection of entire sphere factors for possibly nonadditive Fr\'echet regression on a polysphere. A method for this setting should accommodate dependent factors and nonadditive regression functions while remaining invariant under rotations within each sphere.

We propose local constant and local linear estimators based on product directional kernels with factor-specific bandwidths. We establish convergence rates uniform over predictor locations and factor sets for deterministic bandwidth vectors, extending the pointwise analyses for spherical and toroidal predictors. As in \citet{ImJeon2026}, no bounded bandwidth ratio condition is imposed. For univariate Wasserstein responses supported on a compact interval, the usual local entropy condition can fail even at regular distributions; see \citet{ImJeonPark2026Entropy}. Motivated by the quantile projection approach used for global Fr\'echet regression in \citet{MullerPetersen2026Correction}, we establish uniform convergence rates without imposing this condition.

We select active factors by comparing the increase in cross-fitted prediction risk when each sphere factor is omitted from the full model. Under the stated regularity and threshold conditions, the procedure consistently recovers the active factor set, and the selected local constant and local linear estimators attain the same uniform convergence rates as if the active factor set were known. Under quadratic curvature and bandwidths of the stated order, these rates depend on the active factors' total intrinsic dimension. The selection method allows dependent factors and interactions. 

We assess the finite-sample performance of the proposed methods through simulations with scalar and spherical responses and predictors on $\mathbb S^2\times(\mathbb S^1)^3$. We further illustrate the methods through analyses of Automatic Identification System (AIS) data with vessel speeds as scalar responses and Global Drifter Program (GDP) data with distributions of drifter speeds as Wasserstein responses.

The remainder of the paper is organized as follows. Section~\ref{sec:setting} introduces the estimators and factor-selection procedure. Section~\ref{sec:theory} presents their asymptotic properties. Sections~\ref{sec:simulation} and~\ref{sec:data} report the simulation study and real-data analysis, respectively. Section~\ref{sec:discussion} concludes. All proofs and additional simulation results are provided in the Supplementary Material.

\section{Problem setting and estimators}\label{sec:setting}
\setcounter{equation}{0}

\subsection{Problem setting}\label{sec:targets}

Let $Y$ be a response in a separable metric space $(\M,d_\M)$. Let $p,d_1,\ldots,d_p$ be positive integers and write $\md=(d_1,\ldots,d_p)\tran$ and $|\md|=\sum_{\ell=1}^p d_\ell$. Let $\mX=(\mX_1\tran,\ldots,\mX_p\tran)\tran$ be a predictor on $\Sph^{\md}:=\prod_{\ell=1}^p\Sph^{d_\ell}$, where $\mX_\ell\in\Sph^{d_\ell}$ for $\ell=1,\ldots,p$. For $\mx=\left(\mx_1\tran,\ldots,\mx_p\tran\right)\tran\in\Sph^{\md}$ and $y\in\M$, define the conditional Fr\'echet objective and the conditional Fr\'echet mean by
\begin{align*}
M_\oplus(\mx,y):=\E\left[d_\M^2(Y,y)\mid\mX=\mx\right],\quad m_\oplus(\mx):=\argminop_{y\in\M}M_\oplus(\mx,y).
\end{align*}
The function $m_\oplus:\Sph^{\md}\to\M$ is called the Fr\'echet regression function. For scalar responses with the Euclidean distance, $m_\oplus$ is the ordinary conditional mean. Existence and uniqueness of the conditional Fr\'echet mean are assumed in Section~\ref{sec:response_conditions}.

Write $[q]:=\{1,\ldots,q\}$ for any positive integer $q$. For $A\subseteq[p]$, define $\md_A:=(d_\ell:\ell\in A)\tran$, $\mX_A:=(\mX_\ell\tran:\ell\in A)\tran$, $\mx_A:=(\mx_\ell\tran:\ell\in A)\tran$, $\Sph^{\md_A}:=\prod_{\ell\in A}\Sph^{d_\ell}$ and $|\md_A|:=\sum_{\ell\in A}d_\ell$. For $\mx_A\in\Sph^{\md_A}$ and $y\in\M$, define the corresponding conditional Fr\'echet objective and conditional Fr\'echet mean by
\begin{align*}
M_{\oplus,A}(\mx_A,y):=\E[d_\M^2(Y,y)\mid\mX_A=\mx_A],\quad m_{\oplus,A}(\mx_A):=\argminop_{y\in\M}M_{\oplus,A}(\mx_A,y).
\end{align*}
For $A=\varnothing$, the predictor space is a singleton, and the unconditional Fr\'echet objective and the unconditional Fr\'echet mean are given by
\begin{align*}
M_{\oplus,\varnothing}(y):=\E[d_\M^2(Y,y)],\quad m_{\oplus,\varnothing}:=\argminop_{y\in\M}M_{\oplus,\varnothing}(y).
\end{align*}
Expressions involving $\mx_\varnothing$ are understood as the corresponding unconditional quantities. For $A=[p]$, we have $M_{\oplus,[p]}=M_\oplus$ and $m_{\oplus,[p]}=m_\oplus$. We aim to estimate the Fr\'echet regression function $m_\oplus$ and identify the predictor factors that contribute to prediction.

\subsection{Local constant estimator}\label{sec:lc}

Let $\mathcal D_n:=\{(\mX^{(i)},Y^{(i)})\}_{i=1}^n$ denote an i.i.d. sample of $(\mX,Y)$. For $A=\varnothing$, define the local constant estimator as the empirical Fr\'echet mean
\begin{align*}
\widehat m_{\varnothing,0}\in\argminop_{y\in\M}\frac1n\sum_{i=1}^n d_\M^2(Y^{(i)},y).
\end{align*}

Let $\mh=(h_\ell:\ell\in[p])$ be a deterministic bandwidth vector with positive components, possibly depending on $n$. For each $A\subseteq[p]$, write $\mh_A=(h_\ell:\ell\in A)$, with $\mh_\varnothing$ denoting the empty bandwidth vector. Let $L:[0,\infty)\to[0,\infty)$ be a directional kernel. For $\mx_A,\mz_A\in\Sph^{\md_A}$, define the product directional kernel
\begin{align*}
\cL_{\mx_A,\mh_A}(\mz_A):=\prod_{\ell\in A}\left[h_\ell^{-d_\ell}L\left(\frac{1-\mz_\ell\tran\mx_\ell}{h_\ell^2}\right)\right].
\end{align*}
Define
\begin{align*}
\widehat\mu_{\mh_A}(\mx_A):=\frac1n\sum_{i=1}^n\cL_{\mx_A,\mh_A}\left(\mX_A^{(i)}\right).
\end{align*}
When $\widehat\mu_{\mh_A}(\mx_A)>0$, define the local constant weight by
\begin{align*}
\widehat W_{\mx_A,\mh_A,0}(\mz_A):=\frac{\cL_{\mx_A,\mh_A}(\mz_A)}{\widehat\mu_{\mh_A}(\mx_A)}.
\end{align*}
Under the conditions of \Cref{lem:empirical_design}, $\widehat\mu_{\mh_A}(\mx_A)>0$ uniformly in $\mx_A$ with probability tending to one. The empirical local Fr\'echet objective and the local constant estimator are
\begin{align*}
\widehat M_{\mh_A,0}(\mx_A,y):=\frac1n\sum_{i=1}^n\widehat W_{\mx_A,\mh_A,0}\left(\mX_A^{(i)}\right)d_\M^2(Y^{(i)},y),\quad \widehat m_{\mh_A,0}(\mx_A)\in\argminop_{y\in\M}\widehat M_{\mh_A,0}(\mx_A,y).
\end{align*}
If $\widehat\mu_{\mh_A}(\mx_A)=0$, set $\widehat m_{\mh_A,0}(\mx_A):=\widehat m_{\varnothing,0}$.

\subsection{Local linear estimator}\label{sec:ll}

For $A=\varnothing$, define the local linear estimator as $\widehat m_{\varnothing,1}:=\widehat m_{\varnothing,0}$.

For $A\neq\varnothing$, following \citet{ImJeonPark2025}, we obtain local linear equivalent weights from a first-order approximation for scalar responses and use them to construct an estimator for general metric-space responses. For $\ell\in A$ and $\mx_\ell,\mz_\ell\in\Sph^{d_\ell}$, define
\begin{align*}
\mH_{\mx_\ell}(\mz_\ell):=\begin{cases}\displaystyle\frac{\arccos(\mz_\ell\tran\mx_\ell)}{\sqrt{1-(\mz_\ell\tran\mx_\ell)^2}}\,\mz_\ell,&\mz_\ell\notin\{\mx_\ell,-\mx_\ell\},\\[1ex]
\mz_\ell,&\mz_\ell\in\{\mx_\ell,-\mx_\ell\}.\end{cases}
\end{align*}
For $\ell\in A$, let $T_{\mx_\ell}:=\{\mathbf v\in\R^{d_\ell+1}:\mx_\ell\tran\mathbf v=0\}$ denote the tangent space at $\mx_\ell$. For $\mx_A\in\Sph^{\md_A}$, define
\begin{align*}
\cT_{\mx_A}:=\bigoplus_{\ell\in A}T_{\mx_\ell},\quad \overline{\cT}_{\mx_A}:=\R\oplus\cT_{\mx_A},
\end{align*}
where the latter includes the intercept coordinate. Let $\mI_q$ denote the $q\times q$ identity matrix for any positive integer $q$. For $\mz_A\in\Sph^{\md_A}$, define
\begin{align*}
\bm\Phi_{\mx_A,\mh_A}(\mz_A)&:=\left(\left[h_\ell^{-1}(\mI_{d_\ell+1}-\mx_\ell\mx_\ell\tran)\mH_{\mx_\ell}(\mz_\ell)\right]\tran:\ell\in A\right)\tran\in\cT_{\mx_A},\\
\overline{\bm\Phi}_{\mx_A,\mh_A}(\mz_A)&:=\left(1,\bm\Phi_{\mx_A,\mh_A}(\mz_A)\tran\right)\tran\in\overline{\cT}_{\mx_A}.
\end{align*}
With $\mZ_{\mx_A,\mh_A}^{(i)}:=\overline{\bm\Phi}_{\mx_A,\mh_A}\left(\mX_A^{(i)}\right)$, define the local design matrix
\begin{align*}
\widehat\mGamma_{\mh_A}(\mx_A):=\frac1n\sum_{i=1}^n\cL_{\mx_A,\mh_A}\left(\mX_A^{(i)}\right)\mZ_{\mx_A,\mh_A}^{(i)}\left(\mZ_{\mx_A,\mh_A}^{(i)}\right)\tran.
\end{align*}
Under the conditions of \Cref{lem:empirical_design}, $\widehat\mGamma_{\mh_A}(\mx_A)$ has rank $1+|\md_A|$ uniformly in $\mx_A$ with probability tending to one. Let $\overline{\mathbf{e}}_{\zero}=(1,0,\ldots,0)\tran\in\R^{1+|\md_A|+|A|}$, and let $\dagger$ denote the Moore--Penrose inverse. When $\widehat\mGamma_{\mh_A}(\mx_A)$ has rank $1+|\md_A|$, define the local linear weight by
\begin{align}
\widehat W_{\mx_A,\mh_A,1}(\mz_A):=\cL_{\mx_A,\mh_A}(\mz_A)\overline{\mathbf{e}}_{\zero}\tran\widehat\mGamma_{\mh_A}(\mx_A)^\dagger\overline{\bm\Phi}_{\mx_A,\mh_A}(\mz_A).\label{eq:ll_weight}
\end{align}
The empirical local Fr\'echet objective and local linear estimator are
\begin{align*}
\widehat M_{\mh_A,1}(\mx_A,y):=\frac1n\sum_{i=1}^n\widehat W_{\mx_A,\mh_A,1}\left(\mX_A^{(i)}\right)d_\M^2(Y^{(i)},y),\quad \widehat m_{\mh_A,1}(\mx_A)\in\argminop_{y\in\M}\widehat M_{\mh_A,1}(\mx_A,y).
\end{align*}
If $\widehat\mGamma_{\mh_A}(\mx_A)$ has rank less than $1+|\md_A|$, set $\widehat m_{\mh_A,1}(\mx_A):=\widehat m_{\varnothing,0}$.

For computation of $\widehat W_{\mx_A,\mh_A,1}(\mz_A)$, let $\mB_{\mx_\ell}\in\R^{(d_\ell+1)\times d_\ell}$ have orthonormal columns perpendicular to $\mx_\ell$ for each $\ell\in A$, and set $\mB_{\mx_A}:=\operatorname{diag}(\mB_{\mx_\ell}:\ell\in A)$ and $\overline\mB_{\mx_A}:=\operatorname{diag}(1,\mB_{\mx_A})$. For $\mz_A\in\Sph^{\md_A}$, the corresponding coordinate vector satisfies
\begin{align*}
\overline\mB_{\mx_A}\tran\overline{\bm\Phi}_{\mx_A,\mh_A}(\mz_A)&=\left(1,\left[\mB_{\mx_A}\tran\bm\Phi_{\mx_A,\mh_A}(\mz_A)\right]\tran\right)\tran,\\
\mB_{\mx_A}\tran\bm\Phi_{\mx_A,\mh_A}(\mz_A)&=\left(\left[h_\ell^{-1}\mB_{\mx_\ell}\tran\mH_{\mx_\ell}(\mz_\ell)\right]\tran:\ell\in A\right)\tran.
\end{align*}
Define the coordinate design matrix by
\begin{align*}
\widehat\mQ_{\mh_A}(\mx_A)&:=\frac1n\sum_{i=1}^n\cL_{\mx_A,\mh_A}(\mX_A^{(i)})\left(\overline\mB_{\mx_A}\tran\mZ_{\mx_A,\mh_A}^{(i)}\right)\left(\overline\mB_{\mx_A}\tran\mZ_{\mx_A,\mh_A}^{(i)}\right)\tran\\
&=\overline\mB_{\mx_A}\tran\widehat\mGamma_{\mh_A}(\mx_A)\overline\mB_{\mx_A}.
\end{align*}
When $\widehat\mGamma_{\mh_A}(\mx_A)$ has rank $1+|\md_A|$, $\widehat\mQ_{\mh_A}(\mx_A)$ is invertible, and
\begin{align*}
\widehat\mGamma_{\mh_A}(\mx_A)^\dagger=\overline\mB_{\mx_A}\left(\widehat\mQ_{\mh_A}(\mx_A)\right)^{-1}\overline\mB_{\mx_A}\tran.
\end{align*}
Thus, the local linear weight can be computed as
\begin{align}\label{eq:ll_weight_coord}
\widehat W_{\mx_A,\mh_A,1}(\mz_A)=\cL_{\mx_A,\mh_A}(\mz_A)\overline{\mathbf{e}}_{\zero}\tran\overline\mB_{\mx_A}\left(\widehat\mQ_{\mh_A}(\mx_A)\right)^{-1}\overline\mB_{\mx_A}\tran\overline{\bm\Phi}_{\mx_A,\mh_A}(\mz_A).
\end{align}
By \Cref{prop:supp_basis_invariance}, \eqref{eq:ll_weight_coord} agrees with \eqref{eq:ll_weight} and is independent of the choice of $\{\mB_{\mx_\ell}:\ell\in A\}$.

\begin{remark}\label{rem:invariance_reduction}
For $p=1$, the weights in \eqref{eq:ll_weight} coincide with the spherical local linear weights of \citet{ImJeonPark2025} under the same antipodal convention. For $d_1=\cdots=d_p=1$, they coincide with the joint toroidal local linear weights of \citet{ImJeon2026} away from factorwise antipodes, up to an invertible diagonal rescaling of the slope coordinates.
\end{remark}

\subsection{Factor selection}\label{sec:selection_method}

For $A\subseteq[p]$, define the population prediction risk by
\begin{align*}
R(A):=\E\left[d_\M^2(Y,m_{\oplus,A}(\mX_A))\right].
\end{align*}
When it exists, let $A_0\subseteq[p]$ be the smallest factor set under set inclusion satisfying $R(A_0)=R([p])$. We call $A_0$ the active factor set. For each $\ell\in[p]$, define the leave-one-factor-out (LOFO) risk gap by
\begin{align*}
\Delta_\ell:=R([p]\setminus\{\ell\})-R([p]).
\end{align*}
This risk comparison is related to leave-one-covariate-out methods for assessing predictive importance; see \citet{LeiEtAl2018}. By \Cref{prop:lofo}, $A_0$ exists uniquely and $\Delta_\ell>0$ if and only if $\ell\in A_0$. Thus, factor selection requires only the full model and the $p$ models obtained by omitting one factor at a time, rather than all $2^p$ subsets.

For factor selection, let $\mb=(b_\ell:\ell\in[p])$ be a deterministic bandwidth vector with positive components, possibly depending on $n$, and write $\mb_A=(b_\ell:\ell\in A)$ for each $A\subseteq[p]$, with $\mb_\varnothing$ denoting the empty bandwidth vector. Fix an integer $K\ge2$ and partition the observation indices into $K$ predetermined folds $I_1,\ldots,I_K$, each of size proportional to $n$. For each candidate model $A\subseteq[p]$ and fold $k$, let $\widehat m_{\mb_A,0}^{(-k)}$ denote the local constant estimator fitted using observations outside $I_k$, with the empty-model estimator understood as the empirical Fr\'echet mean of these observations. Define the cross-fitted (CF) prediction risk estimator by
\begin{align*}
\widehat R_{\mathrm{CF}}(A):=\frac1n\sum_{k=1}^K\sum_{i\in I_k}d_\M^2\left(Y^{(i)},\widehat m_{\mb_A,0}^{(-k)}\left(\mX_A^{(i)}\right)\right).
\end{align*}

We now estimate the active factor set by thresholding the CF risk gaps. For an independent copy $Y'$ of $Y$, define
\begin{align}
V_Y:=\frac12\E[d_\M^2(Y,Y')],\quad \widehat V_Y:=\frac{1}{n(n-1)}\sum_{1\le i<j\le n}d_\M^2(Y^{(i)},Y^{(j)}).\label{eq:empirical_response_dispersion}
\end{align}
For scalar responses with the Euclidean distance, $\widehat V_Y$ is the unbiased sample variance. For a fixed tuning constant $c>0$ and a deterministic sequence $a_n>0$, define the selection threshold by
\begin{align*}
\widehat\tau_n:=c\widehat V_Ya_n.
\end{align*}
Sufficient conditions on $a_n$ for selection consistency are given in Section~\ref{sec:selection_theory}. The estimated LOFO risk gaps and selected factor set are
\begin{align*}
\widehat\Delta_\ell&:=\widehat R_{\mathrm{CF}}([p]\setminus\{\ell\})-\widehat R_{\mathrm{CF}}([p]),\quad \ell\in[p],\\
\widehat A_{\mb}&:=\{\ell\in[p]:\widehat\Delta_\ell>\widehat\tau_n\}.
\end{align*}
Finally, for $s\in\{0,1\}$ and $\mx\in\Sph^{\md}$, define the selected estimator by
\begin{align*}
\widehat m_{\mb,\mh,s}(\mx):=\widehat m_{\mh_{\widehat A_{\mb}},s}(\mx_{\widehat A_{\mb}}),
\end{align*}
where the right-hand side is understood as $\widehat m_{\varnothing,s}$ when $\widehat A_{\mb}=\varnothing$.

\section{Asymptotic theory}\label{sec:theory}
\setcounter{equation}{0}

We first establish uniform convergence rates for Wasserstein and general metric-space responses, and then study factor-selection consistency and post-selection rates. For $\mathbf t=(t_\ell:\ell\in[p])\in\{\mb,\mh\}$ and $A\subseteq[p]$, define
\begin{align*}
v_{\mathbf t}(A):=\prod_{\ell\in A}t_\ell^{d_\ell},\quad b_{\mathbf t}(A):=\sum_{\ell\in A}t_\ell^2,
\end{align*}
with $v_{\mathbf t}(\varnothing)=1$ and $b_{\mathbf t}(\varnothing)=0$. For a nonempty $A\subseteq[p]$, define the population versions of $\widehat\mu_{\mh_A}(\mx_A)$ and $\widehat\mGamma_{\mh_A}(\mx_A)$ by
\begin{align*}
\widetilde\mu_{\mh_A}(\mx_A):=\E\left[\cL_{\mx_A,\mh_A}(\mX_A)\right],\quad
\widetilde\mGamma_{\mh_A}(\mx_A):=\E\left[\cL_{\mx_A,\mh_A}(\mX_A)\overline{\bm\Phi}_{\mx_A,\mh_A}(\mX_A)\overline{\bm\Phi}_{\mx_A,\mh_A}(\mX_A)\tran\right].
\end{align*}
Under the conditions stated below, \Cref{lem:population_design} ensures that $\widetilde\mu_{\mh_A}(\mx_A)>0$ and the restriction of $\widetilde\mGamma_{\mh_A}(\mx_A)$ to $\overline{\cT}_{\mx_A}$ is positive definite for all sufficiently small bandwidths. Define the corresponding population equivalent weights by
\begin{align*}
\widetilde W_{\mx_A,\mh_A,0}(\mz_A):=\frac{\cL_{\mx_A,\mh_A}(\mz_A)}{\widetilde\mu_{\mh_A}(\mx_A)},\quad
\widetilde W_{\mx_A,\mh_A,1}(\mz_A):=\cL_{\mx_A,\mh_A}(\mz_A)\overline{\mathbf{e}}_{\zero}\tran\widetilde\mGamma_{\mh_A}(\mx_A)^\dagger\overline{\bm\Phi}_{\mx_A,\mh_A}(\mz_A).
\end{align*}
For $s\in\{0,1\}$, define the population local Fr\'echet objective and its minimizer by
\begin{align*}
\widetilde M_{\mh_A,s}(\mx_A,y):=\E\left[\widetilde W_{\mx_A,\mh_A,s}(\mX_A)d_\M^2(Y,y)\right],\quad \widetilde m_{\mh_A,s}(\mx_A):=\argminop_{y\in\M}\widetilde M_{\mh_A,s}(\mx_A,y).
\end{align*}
For $A=\varnothing$, both population weights equal one, so $\widetilde M_{\varnothing,s}=M_{\oplus,\varnothing}$ and $\widetilde m_{\varnothing,s}=m_{\oplus,\varnothing}$ for $s\in\{0,1\}$.

\subsection{Common regularity conditions}\label{sec:kernel_conditions}

\begin{customcon}{L1}\label{con:L1}
The kernel $L:[0,\infty)\to[0,\infty)$ is nonincreasing and satisfies
\begin{align*}
0<\int_0^\infty L^k(r^2)r^m\dd r<\infty,\quad k\in\{1,2\},\ m=0,1,\ldots.
\end{align*}
\end{customcon}

Similar kernel conditions were used in \citet{ImJeonPark2025} and \citet{ImJeon2026}. Condition~\ref{con:L1} is satisfied by the uniform, Epanechnikov, exponential and von Mises kernels, among others.

\begin{customcon}{B1}\label{con:B1}
For each $\mathbf t=(t_\ell:\ell\in[p])\in\{\mb,\mh\}$, as $n\to\infty$,
\begin{align*}
\max_{\ell\in[p]}t_\ell\to0,\quad \frac{nv_{\mathbf t}([p])}{\log n}\to\infty.
\end{align*}
\end{customcon}

Condition~\ref{con:B1} requires the bandwidths to shrink while the effective local sample sizes diverge. For a single common bandwidth $h$ in the full model, it reduces to $h\to0$ and $nh^{|\md|}/\log n\to\infty$, as in \citet{ImJeon2026Riemannian}. As in \citet{ImJeon2026}, no bounded bandwidth-ratio condition is imposed, unlike \citet{RuppertWand1994}.

Let $\omega_d$ denote surface-area measure on $\Sph^d$, and let $\omega_{\md_A}$ denote the corresponding product measure on $\Sph^{\md_A}$. Let $f$ be the density of $\mX$ with respect to $\omega_{\md}$. Denote $\zero_q\in\R^q$ as the zero vector for any positive integer $q$. For a function $\phi:\Sph^{\md}\to\R$, define its factorwise homogeneous extension $\bar\phi:\mathcal U_{\md}\to\R$, as in \citet{ImJeon2026}, by
\begin{align}\label{eq:homogeneous_extension}
\bar\phi(\mz):=\phi\left(\frac{\mz_\ell}{\|\mz_\ell\|}:\ell\in[p]\right),\quad \mz\in\mathcal U_{\md},
\end{align}
where $\mathcal U_{\md}:=\prod_{\ell=1}^p\left(\R^{d_\ell+1}\setminus\{\zero_{d_\ell+1}\}\right)$ and $\mz=(\mz_1\tran,\ldots,\mz_p\tran)\tran$. Let $\bar f:\mathcal U_{\md}\to\R$ denote the homogeneous extension of $f$.

\begin{customcon}{D1}\label{con:D1}
$f$ is strictly positive, and $\bar f$ is twice continuously differentiable.
\end{customcon}

Let $P_Y$ denote the distribution of $Y$ and $P_{Y\mid\mX=\mx}$ the conditional distribution of $Y$ given $\mX=\mx$. Suppose that $P_{Y\mid\mX=\mx}$ is absolutely continuous with respect to $P_Y$ for every $\mx\in\Sph^{\md}$. For $y\in\M$, using a jointly measurable version of the Radon--Nikodym derivative, define
\begin{align}\label{eq:conditional_density_ratio}
g_y(\mx):=\frac{\mathrm dP_{Y\mid\mX=\mx}}{\mathrm dP_Y}(y).
\end{align}
Let $\bar g_y$ denote the homogeneous extension of $g_y$ defined by \eqref{eq:homogeneous_extension}. Whenever $\bar g_y$ is twice continuously differentiable on $\mathcal U_{\md}$, let $\nabla\bar g_y:\mathcal U_{\md}\to\R^{|\md|+p}$ and $\nabla^2\bar g_y:\mathcal U_{\md}\to\R^{(|\md|+p)\times(|\md|+p)}$ denote its Euclidean gradient and Hessian, respectively. In this case, define
\begin{align}\label{eq:conditional_density_envelope}
G(y):=\sup_{\mx\in\Sph^{\md}}\left(g_y(\mx)+\|\nabla\bar g_y(\mx)\|+\|\nabla^2\bar g_y(\mx)\|_{\mathrm{op}}\right),
\end{align}
and set $G(y)=\infty$ otherwise. Also, define
\begin{align*}
D_\M(y):=\max\left\{1,d_\M(y,m_{\oplus,\varnothing})\right\}.
\end{align*}

\begin{customcon}{D2}\label{con:D2}
For every $y\in\M$, $\bar g_y$ is twice continuously differentiable, and $\E[D_\M(Y)^2G(Y)]<\infty$.
\end{customcon}

By compactness, Condition~\ref{con:D1} implies that $f$ is bounded away from zero and has bounded first and second derivatives on the polysphere. Condition~\ref{con:D2} provides an integrable bound on the conditional density ratios and their derivatives. Similar smoothness conditions were used in \citet{ImJeonPark2025} and \citet{ImJeon2026}. A Gaussian location model is considered in \Cref{ex:gaussian_density_envelope}, and the corresponding conditions for reduced models follow from \Cref{lem:marginalization}.

\subsection{Uniform convergence rates for Wasserstein responses}\label{sec:wasserstein}

Let $\mathcal W_2([a,b])$ be the space of probability distributions supported on a compact interval $[a,b]$, equipped with the 2-Wasserstein distance $W_2$, where $-\infty<a<b<\infty$. For $y\in\mathcal W_2([a,b])$, let $F_y$ denote its distribution function and $Q_y$ its quantile function, defined by
\begin{align*}
Q_y(u):=\inf\{t\in[a,b]:F_y(t)\ge u\},\quad u\in(0,1).
\end{align*}
Then $W_2(y,z)=\|Q_y-Q_z\|_{L^2(0,1)}$. The quantile functions form a closed convex subset $\mathcal Q_{[a,b]}$ of $L^2(0,1)$; see \citet{PetersenMuller2019}.

For $s\in\{0,1\}$, a nonempty $A\subseteq[p]$ and $\mx_A\in\Sph^{\md_A}$, suppose that $\widehat\mu_{\mh_A}(\mx_A)>0$ if $s=0$, and that $\widehat\mGamma_{\mh_A}(\mx_A)$ is invertible on $\overline{\cT}_{\mx_A}$ if $s=1$. Then the local constant ($s=0$) or local linear ($s=1$) estimator satisfies
\begin{align}
Q_{\widehat m_{\mh_A,s}(\mx_A)}=\Pi_{\mathcal Q_{[a,b]}}\left(\frac1n\sum_{i=1}^n\widehat W_{\mx_A,\mh_A,s}\left(\mX_A^{(i)}\right)Q_{Y^{(i)}}\right),\label{eq:wasserstein_projection}
\end{align}
where $\Pi_{\mathcal Q_{[a,b]}}$ denotes the projection onto $\mathcal Q_{[a,b]}$ in $L^2(0,1)$. For $A=\varnothing$, both estimators are the empirical Fr\'echet mean, whose quantile function is $n^{-1}\sum_{i=1}^nQ_{Y^{(i)}}$.

\begin{theorem}\label{thm:wasserstein}
Suppose $\M=\mathcal W_2([a,b])$ and Conditions~\ref{con:L1}, \ref{con:B1} and \ref{con:D1}--\ref{con:D2} hold. Then, for each $s\in\{0,1\}$ and $\mathbf t\in\{\mb,\mh\}$, uniformly over $A\subseteq[p]$,
\begin{align*}
\sup_{\mx_A\in\Sph^{\md_A}}W_2\left(\widehat m_{\mathbf t_A,s}(\mx_A),m_{\oplus,A}(\mx_A)\right)=O_\P\left(b_{\mathbf t}(A)+\sqrt{\frac{\log n}{nv_{\mathbf t}(A)}}\right).
\end{align*}
\end{theorem}

For a fixed nonempty $A\subseteq[p]$, choosing $h_\ell\asymp(\log n/n)^{1/(|\md_A|+4)}$ for all $\ell\in A$ yields the rate $O_\P\left((\log n/n)^{2/(|\md_A|+4)}\right)$, provided Condition~\ref{con:B1} holds. This matches the minimax-optimal sup-norm rate for twice-smooth scalar regression with a $|\md_A|$-dimensional Euclidean predictor established by \citet{Stone1982}. For $|\md_A|=1$, this rate is faster in order than the uniform upper bound $O_\P(n^{-1/(3+\epsilon)})$ obtained by \citet{ChenMuller2022} under quadratic curvature, although the predictor spaces and regularity assumptions differ.

\subsection{Uniform convergence rates for responses in general metric spaces}\label{sec:response_conditions}\label{sec:uniform_rates}

For responses in a general metric space, we impose an additional condition on the response distribution.

\begin{customcon}{B2}\label{con:B2}
For each $\mathbf t\in\{\mb,\mh\}$, as $n\to\infty$,
\begin{align*}
\max_{1\le i\le n}D_\M(Y^{(i)})=O_\P\left(\sqrt{\frac{nv_{\mathbf t}([p])}{\log n}}\right).
\end{align*}
\end{customcon}

Condition~\ref{con:B2} controls the largest observed response distance in the uniform stochastic bounds. It holds automatically for bounded responses under Condition~\ref{con:B1}. If $\E[D_\M(Y)^q]<\infty$ for some $q>2$, a sufficient condition is
\begin{align*}
\liminf_{n\to\infty}\frac{n^{1-2/q}}{\log n}\min\{v_{\mb}([p]),v_{\mh}([p])\}>0.
\end{align*}
Further sufficient conditions are given in \Cref{prop:TB}.

We next impose conditions on the Fr\'echet minimizers and the response space $\M$. In Conditions~\ref{con:M1}--\ref{con:M4}, (P), (L) and (E) refer to population, localized population and empirical conditions, respectively.

\begin{customcon}{M1}\label{con:M1}
The following conditions hold.
\begin{itemize}[leftmargin=2.2em]
\item[(P)] $m_{\oplus,\varnothing}$ exists and is unique, with $\E[d_\M^2(Y,m_{\oplus,\varnothing})]<\infty$.
\item[(E)] $\widehat m_{\varnothing,0}$ exists almost surely and can be chosen as a measurable function of $\mathcal D_n$.
\end{itemize}
\end{customcon}

\begin{customcon}{M2}\label{con:M2}
The following condition holds.
\begin{itemize}[leftmargin=2.2em]
\item[(P)] For every $A\subseteq[p]$, $m_{\oplus,A}(\mx_A)$ exists and is unique for all $\mx_A\in\Sph^{\md_A}$, and $m_{\oplus,A}:\Sph^{\md_A}\to\M$ is measurable. Furthermore, for every $\epsilon>0$,
\begin{align*}
\inf_{A\subseteq[p]}\inf_{\mx_A\in\Sph^{\md_A}}\inf_{y\in\M:\,d_\M\left(y,m_{\oplus,A}(\mx_A)\right)>\epsilon}\left[M_{\oplus,A}(\mx_A,y)-M_{\oplus,A}\left(\mx_A,m_{\oplus,A}(\mx_A)\right)\right]>0.
\end{align*}
\end{itemize}
\end{customcon}

\begin{customcon}{M3}\label{con:M3}
For a fixed $s\in\{0,1\}$, there exists $\widetilde h>0$ such that, for every nonempty $A\subseteq[p]$, $\mh_A$ satisfying $\max_{\ell\in A}h_\ell\le\widetilde h$ and $\mx_A\in\Sph^{\md_A}$, the following conditions hold.
\begin{itemize}[leftmargin=2.2em]
\item[(L)] $\widetilde m_{\mh_A,s}(\mx_A)$ exists and is unique, and $\widetilde m_{\mh_A,s}:\Sph^{\md_A}\to\M$ is measurable.
\item[(E)] $\widehat m_{\mh_A,s}(\mx_A)$ exists simultaneously for all $\mx_A\in\Sph^{\md_A}$ satisfying $\widehat\mu_{\mh_A}(\mx_A)>0$ if $s=0$, or $\operatorname{rank}(\widehat\mGamma_{\mh_A}(\mx_A))=1+|\md_A|$ if $s=1$, with probability one. The estimator $\widehat m_{\mh_A,s}$ can be chosen jointly measurable in $(\mathcal D_n,\mx_A)$.
\end{itemize}
\end{customcon}

\begin{customcon}{M4}\label{con:M4}
For a fixed $s\in\{0,1\}$, there exist constants $\beta_\oplus>1$, $c_\oplus>0$, $\eta_\oplus>0$ and $0<h_\oplus\le\widetilde h$ such that, for every $A\subseteq[p]$, $\mx_A\in\Sph^{\md_A}$ and $y\in\M$, the following conditions hold.
\begin{itemize}[leftmargin=2.2em]
\item[(P)] If $d_\M\left(y,m_{\oplus,A}(\mx_A)\right)\le\eta_\oplus$, then
\begin{align}
M_{\oplus,A}(\mx_A,y)-M_{\oplus,A}\left(\mx_A,m_{\oplus,A}(\mx_A)\right)\ge c_\oplus d_\M\left(y,m_{\oplus,A}(\mx_A)\right)^{\beta_\oplus}.\label{eq:margin}
\end{align}
\item[(L)] For $A\neq\varnothing$ and every $\mh_A$ satisfying $\max_{\ell\in A}h_\ell\le h_\oplus$, if $d_\M\left(y,\widetilde m_{\mh_A,s}(\mx_A)\right)\le\eta_\oplus$, then
\begin{align}
\widetilde M_{\mh_A,s}(\mx_A,y)-\widetilde M_{\mh_A,s}\left(\mx_A,\widetilde m_{\mh_A,s}(\mx_A)\right)\ge c_\oplus d_\M\left(y,\widetilde m_{\mh_A,s}(\mx_A)\right)^{\beta_\oplus}.\label{eq:localized_margin}
\end{align}
\end{itemize}
\end{customcon}

For $z\in\M$ and $r>0$, let $B_\M(z,r):=\{y\in\M:d_\M(y,z)<r\}$, and let $N(\epsilon,B,d_\M)$ denote the $\epsilon$-covering number of $B\subseteq\M$. For each $A\subseteq[p]$, define $\cK_A:=\{m_{\oplus,A}(\mx_A):\mx_A\in\Sph^{\md_A}\}$ and $d_\M(z,\cK_A):=\inf_{u\in\cK_A}d_\M(z,u)$.

\begin{customcon}{M5}\label{con:M5}
There exist $r_\M>0$, $\delta_\M>0$ and $C_\M<\infty$ such that for every $A\subseteq[p]$ and $z\in\M$,
\begin{align*}
\int_0^{1/2}\sqrt{1+\log N\left(\delta\epsilon,B_\M(z,\delta),d_\M\right)}\dd\epsilon\le C_\M,
\end{align*}
whenever $d_\M(z,\cK_A)<r_\M$ and $0<\delta\le\delta_\M$.
\end{customcon}

\begin{customcon}{M6}\label{con:M6}
Every bounded ball in $(\M,d_\M)$ is totally bounded.
\end{customcon}

Conditions~\ref{con:M1}--\ref{con:M3} concern the existence, uniqueness, measurability and separation of Fr\'echet minimizers, while Condition~\ref{con:M4} imposes local curvature. Condition~\ref{con:M5} is a local metric entropy condition. Similar curvature and entropy conditions were used in \citet{PetersenMuller2019}, \citet{ImJeonPark2025} and \citet{ImJeon2026}. Condition~\ref{con:M6} is used to establish preliminary uniform consistency and holds when closed bounded subsets of $\M$ are compact. Thus, the response space need not be bounded. An alternative to Condition~\ref{con:M6} is given in \Cref{rmk:localization_alternative}.

For responses in general metric spaces, let $\beta_\oplus$ be the exponent in Condition~\ref{con:M4}. For Wasserstein responses, set $\beta_\oplus=2$. For $\mathbf t\in\{\mb,\mh\}$, define
\begin{align*}
r_{n,\mathbf t}(A):=b_{\mathbf t}(A)^{1/(\beta_\oplus-1)}+\left(\frac{\log n}{nv_{\mathbf t}(A)}\right)^{1/[2(\beta_\oplus-1)]}.
\end{align*}

\begin{theorem}\label{thm:grid_uniform}
Suppose Conditions~\ref{con:L1}, \ref{con:B1}--\ref{con:B2}, \ref{con:D1}--\ref{con:D2} and \ref{con:M1}--\ref{con:M6} hold, where Conditions~\ref{con:M3}--\ref{con:M4} are imposed for a fixed $s\in\{0,1\}$. Then, for each $\mathbf t\in\{\mb,\mh\}$,
\begin{align}
\max_{A\subseteq[p]}\frac{1}{r_{n,\mathbf t}(A)}\sup_{\mx_A\in\Sph^{\md_A}}d_\M\left(\widehat m_{\mathbf t_A,s}(\mx_A),m_{\oplus,A}(\mx_A)\right)=O_\P(1).\label{eq:grid_uniform_rate}
\end{align}
\end{theorem}

The bounds hold simultaneously over all factor sets because $p$ is fixed. For a fixed nonempty $A\subseteq[p]$, if $\beta_\oplus=2$ and $h_\ell\asymp(\log n/n)^{1/(|\md_A|+4)}$ for every $\ell\in A$, then, under the stated conditions,
\begin{align*}
\sup_{\mx_A\in\Sph^{\md_A}}d_\M\left(\widehat m_{\mh_A,s}(\mx_A),m_{\oplus,A}(\mx_A)\right)=O_\P\left((\log n/n)^{2/(|\md_A|+4)}\right).
\end{align*}
This agrees with the uniform convergence rate of \citet{ImJeon2026Riemannian} for Riemannian predictors of intrinsic dimension $|\md_A|$ under quadratic curvature. Condition~\ref{con:B2} must hold at this bandwidth order; a finite fourth response moment alone need not suffice.

\subsection{Factor selection and post-selection estimation}\label{sec:selection_theory}

Let $\mathfrak C_p:=\{[p]\}\cup\{[p]\setminus\{\ell\}:\ell\in[p]\}$ denote the full and LOFO models. For responses in general metric spaces, suppose Condition~\ref{con:M4} holds for both $s=0$ and $s=1$ with the same exponent $\beta_\oplus$. Define
\begin{align*}
e_n:=r_{n,\mb}([p])+n^{-1/2}.
\end{align*}

\begin{theorem}\label{thm:selection}
Suppose the conditions of \Cref{thm:wasserstein} or \Cref{thm:grid_uniform} hold for the local constant estimators with bandwidth $\mb$. Suppose also that $0<V_Y<\infty$, $a_n=o(1)$ and $e_n=o(a_n)$. Then
\begin{align*}
\P\left(\widehat A_{\mb}=A_0\right)\to1.
\end{align*}
\end{theorem}

The conditions on $a_n$ make the threshold vanish more slowly than the estimation error in the LOFO risk gaps. Since each gap compares squared losses on the same observations, its stochastic error can be controlled without an additional fourth-moment assumption.

\begin{theorem}\label{thm:post_selection}
Suppose the conditions of \Cref{thm:selection} hold. For a fixed $s\in\{0,1\}$, suppose also that the conditions of \Cref{thm:wasserstein} or \Cref{thm:grid_uniform} hold for the local constant ($s=0$) or local linear ($s=1$) estimator with bandwidth $\mh$. Then
\begin{align*}
\sup_{\mx\in\Sph^{\md}}d_\M\left(\widehat m_{\mb,\mh,s}(\mx),m_\oplus(\mx)\right)=O_\P\left(r_{n,\mh}(A_0)\right).
\end{align*}
\end{theorem}

For $A_0\neq\varnothing$ and $\beta_\oplus=2$, if $h_\ell\asymp(\log n/n)^{1/(|\md_{A_0}|+4)}$ for every $\ell\in A_0$, then, under the conditions of \Cref{thm:post_selection},
\begin{align*}
\sup_{\mx\in\Sph^{\md}}d_\M\left(\widehat m_{\mb,\mh,s}(\mx),m_\oplus(\mx)\right)=O_\P\left((\log n/n)^{2/(|\md_{A_0}|+4)}\right).
\end{align*}
Thus, the selected estimator attains the same uniform convergence rate as if the active factor set $A_0$ were known. This rate agrees with the classical uniform rate for twice-smooth Euclidean nonparametric regression in dimension $|\md_{A_0}|$.

\section{Simulation study}\label{sec:simulation}

We conduct two simulation studies with predictors on $\Sph^2\times(\Sph^1)^3$. The first uses scalar responses and the second uses spherical responses. Throughout the simulation and real-data analyses, LC and LL denote the local constant and local linear estimators, respectively. In both studies, we consider $n\in\{500,1000,2000\}$, with strong and weak signal settings under low and high noise levels.

Factor selection uses $5$-fold cross-fitted LC prediction risks for the full and LOFO models, with $3$-fold cross validation for bandwidth selection within each outer training sample. In both simulation studies, $\widehat A$ denotes the factor set selected using bandwidths chosen by cross validation. We use
\begin{align}
\widehat\tau_n=0.005\,\widehat V_Y\left(\frac{\log n}{n}\right)^{1/6},\label{eq:simulation_threshold}
\end{align}
where $\widehat V_Y$ is defined in \eqref{eq:empirical_response_dispersion} and equals the sample variance for scalar responses. The multiplier $0.005$ was fixed before the reported experiments using a calibration study. Final bandwidths are selected by $5$-fold cross validation using common folds across methods. Threshold sensitivity results for both simulations are reported in Section~\ref{supp:auxiliary} of the Supplementary Material.

We use $R=200$ Monte Carlo replications per setting and evaluate the performance of the estimators by
\begin{align*}
\mathrm{MISE}:=\frac1R\sum_{r=1}^R\int_{\Sph^2\times(\Sph^1)^3}d_\M^2\left(\widehat m_\oplus^{[r]}(\mx),m_\oplus(\mx)\right)\dd P_{\mX}(\mx),
\end{align*}
where $\widehat m_\oplus^{[r]}$ denotes the fitted regression function in replication $r\in[R]$ for each method, and $d_\M$ is the Euclidean distance for scalar responses and the geodesic distance for spherical responses. Each integral is approximated using $5000$ independently generated test predictors, shared across methods within each replication. 

\subsection{Case of scalar responses}\label{sec:simulation_scalar}

We generate the predictors as follows. For $\kappa>0$, define the circular density
\begin{align*}
q_\kappa(\mathbf z\mid \mathbf a):=\frac{\exp(\kappa \mathbf a\tran\mathbf z)}{\int_{\Sph^1}\exp(\kappa \mathbf a\tran\mathbf u)\dd\omega_1(\mathbf u)},\quad \mathbf z\in\Sph^1,\ \mathbf a\in\R^2,\ \|\mathbf a\|\le1.
\end{align*}
For $i=1,\ldots,n$, independently generate $\Phi_i\sim\mathrm{Unif}(-\pi,\pi)$ and $U_i\sim\mathrm{Unif}(-1,1)$, and set
\begin{align*}
\mX_{i1}:=\left(\sqrt{1-U_i^2}\cos(\Phi_i),\sqrt{1-U_i^2}\sin(\Phi_i),U_i\right)\tran\in\Sph^2.
\end{align*}
For $\mx_1=\left(x_{11},x_{12},x_{13}\right)\tran\in\Sph^2$, let $\mathbf a(\mx_1):=\left(x_{11},x_{12}\right)\tran\in\R^2$. Given $\mX_{i1}$, draw $\mX_{i2}$ and $\mX_{i4}$ independently from $q_{1.5}(\cdot\mid \mathbf a(\mX_{i1}))$, and then draw $\mX_{i3}$ from $q_1(\cdot\mid \mX_{i2})$, independently of $(\mX_{i1},\mX_{i4})$ given $\mX_{i2}$. The resulting joint density is strictly positive and smooth on the polysphere.

Given the predictors, we generate scalar responses as follows. For $\mx_1=\left(x_{11},x_{12},x_{13}\right)\tran\in\Sph^2$ and $\mx_2=\left(\cos(\theta_2),\sin(\theta_2)\right)\tran\in\Sph^1$, define
\begin{align*}
g_1(\mx_1)&:=1.2x_{11}+0.8x_{12}-0.65x_{13}+0.55x_{12}^2,\quad
g_2(\theta_2):=\sin(\theta_2)+0.5\cos(2\theta_2),\\
g_{12}(\mx_1,\theta_2)&:=x_{12}\cos(\theta_2).
\end{align*}
Writing $\mX_i=(\mX_{i1}\tran,\ldots,\mX_{i4}\tran)\tran$ and $\mx=(\mx_1\tran,\ldots,\mx_4\tran)\tran\in\Sph^2\times(\Sph^1)^3$, where $\mx_2=(\cos(\theta_2),\sin(\theta_2))\tran$, we generate
\begin{align*}
m_\oplus(\mx)&:=1.8g_1(\mx_1)+\beta_2g_2(\theta_2)+\beta_{12}g_{12}(\mx_1,\theta_2),\\
Y_i&:=m_\oplus(\mX_i)+\sigma\epsilon_i,\quad \epsilon_i\sim N(0,1),
\end{align*}
where the errors are independent across observations and independent of the predictors.

We take $(\beta_2,\beta_{12})=(1.2,0.8)$ and $(0.6,0.4)$ for the strong and weak signal settings, respectively, and $\sigma=0.25$ and $1$ for the low and high noise settings. The regression function is nonadditive, with $A_0=\{1,2\}$ in all settings. Hence the active and full intrinsic dimensions are $3$ and $5$, respectively, while the inactive factors remain dependent on the active factors. 

We compare five classes of estimators:
\begin{itemize}
\item Full: the proposed LC and LL estimators using all four factors.
\item Selected: the proposed LC and LL estimators using the selected factors.
\item Oracle: the proposed LC and LL estimators using the true active factors as infeasible benchmarks.
\item Euclidean: LC and LL estimators based on the Euclidean-predictor Fr\'echet regression framework of \citet{PetersenMuller2019}, applied to $\left(\Phi_i,\arcsin(U_i),\Theta_{i2},\Theta_{i3},\Theta_{i4}\right)\tran$, where $\mX_{i\ell}=\left(\cos(\Theta_{i\ell}),\sin(\Theta_{i\ell})\right)\tran$ for $\ell=2,3,4$. They use a common bandwidth for longitude and latitude and a separate bandwidth for each remaining angle.
\item Riemannian: the LC and LL estimators of \citet{ImJeon2026Riemannian}, using product geodesic distance and normal-coordinate volume correction.
\end{itemize}
The proposed estimators use the product von Mises kernel. The Euclidean estimators use a product Epanechnikov kernel, whereas the Riemannian estimators use an isotropic Epanechnikov kernel. For $0<a<b$ and an integer $m\ge2$, let $\mathcal G_m(a,b):=\left\{a(b/a)^{(j-1)/(m-1)}:j=1,\ldots,m\right\}$ denote the grid of $m$ logarithmically spaced points between $a$ and $b$. The candidate bandwidth grids are
\begin{align*}
\mathcal H_{[4]}^{\mathrm{sel}}&:=\mathcal G_{20}(0.002,0.4)\times\mathcal G_{20}(0.05,6)^3,\quad \mathcal H_{[4]}^{\mathrm{est}}:=\mathcal G_{20}(0.002,0.8)\times\mathcal G_{20}(0.05,6)^3,\\
\mathcal H^{\mathrm{Euc}}&:=\mathcal G_{20}(0.05,7.5)\times\mathcal G_{20}(0.05,20)^3,\quad \mathcal H^{\mathrm{Rim}}:=\mathcal G_{50}(0.2,3),
\end{align*}
where $\mathcal H_{[4]}^{\mathrm{sel}}$ and $\mathcal H_{[4]}^{\mathrm{est}}$ denote the proposed bandwidth grids for factor selection and final estimation, respectively, and $\mathcal H^{\mathrm{Euc}}$ and $\mathcal H^{\mathrm{Rim}}$ denote the grids for the Euclidean and Riemannian estimators, respectively. For $A\subseteq[4]$, the proposed grids retain the corresponding coordinates.

Table~\ref{tab:simulation_selection} reports factor-selection frequencies. Both active factors are selected in every replication. Exact selection is at least $0.915$ at $n=500$ and equals $1.000$ in all settings at $n=2000$. All remaining errors arise from including $\mX_3$ or $\mX_4$.

\begin{table}[!t]
\centering
\caption{Factor-selection frequencies for scalar responses based on $200$ Monte Carlo replications per setting. Exact gives the frequency of $\widehat A=A_0$ for $A_0=\{1,2\}$, and $\mX_1,\ldots,\mX_4$ give individual inclusion frequencies. Here $\sigma$ denotes the noise standard deviation.}
\label{tab:simulation_selection}
\begin{tabular}{cccccccc}
\toprule
Signal & $\sigma$ & $n$ & Exact & $\mX_1$ & $\mX_2$ & $\mX_3$ & $\mX_4$\\
\midrule
Strong & 0.25 & 500 & 0.990 & 1.000 & 1.000 & 0.010 & 0.005\\
Strong & 0.25 & 1000 & 1.000 & 1.000 & 1.000 & 0.000 & 0.000\\
Strong & 0.25 & 2000 & 1.000 & 1.000 & 1.000 & 0.000 & 0.000\\
\midrule
Strong & 1 & 500 & 0.930 & 1.000 & 1.000 & 0.035 & 0.035\\
Strong & 1 & 1000 & 1.000 & 1.000 & 1.000 & 0.000 & 0.000\\
Strong & 1 & 2000 & 1.000 & 1.000 & 1.000 & 0.000 & 0.000\\
\midrule
Weak & 0.25 & 500 & 1.000 & 1.000 & 1.000 & 0.000 & 0.000\\
Weak & 0.25 & 1000 & 1.000 & 1.000 & 1.000 & 0.000 & 0.000\\
Weak & 0.25 & 2000 & 1.000 & 1.000 & 1.000 & 0.000 & 0.000\\
\midrule
Weak & 1 & 500 & 0.915 & 1.000 & 1.000 & 0.035 & 0.055\\
Weak & 1 & 1000 & 0.985 & 1.000 & 1.000 & 0.015 & 0.000\\
Weak & 1 & 2000 & 1.000 & 1.000 & 1.000 & 0.000 & 0.000\\
\bottomrule
\end{tabular}
\end{table}

Table~\ref{tab:simulation_mise} reports the MISEs. Selected LL has the smallest MISE among the data-driven estimators in every setting and is close to Oracle LL. Factor selection has a smaller effect on LC estimation. Full LC has a slightly smaller MISE than Selected LC only in the weak signal setting with low noise at $n=1000$. The proposed estimators also have smaller MISEs than the Euclidean and Riemannian estimators in every setting.

\begin{table}[!t]
\centering
\caption{MISEs for scalar responses based on $200$ Monte Carlo replications per setting. All values are multiplied by $10$. Full, Selected and Oracle use all four, the selected and the true active factors, respectively. Euclidean and Riemannian use all four factors. Bold values indicate the smallest MISE among the data-driven methods, excluding Oracle.}
\label{tab:simulation_mise}
\begingroup
\small
\resizebox{\linewidth}{!}{%
\begin{tabular}{ccc*{10}{c}}
\toprule
& & & \multicolumn{2}{c}{Full} & \multicolumn{2}{c}{Selected} & \multicolumn{2}{c}{Oracle} & \multicolumn{2}{c}{Riemannian} & \multicolumn{2}{c}{Euclidean}\\
\cmidrule(lr){4-5}\cmidrule(lr){6-7}\cmidrule(lr){8-9}\cmidrule(lr){10-11}\cmidrule(lr){12-13}
Signal & $\sigma$ & $n$ & LC & LL & LC & LL & LC & LL & LC & LL & LC & LL\\
\midrule
Strong & 0.25 & 500 & 1.0170 & 0.6090 & 1.0130 & \textbf{0.4370} & 1.0130 & 0.4360 & 7.1750 & 6.2080 & 2.6710 & 2.8730\\
Strong & 0.25 & 1000 & 0.6140 & 0.3580 & 0.6130 & \textbf{0.2730} & 0.6130 & 0.2730 & 5.6550 & 4.7740 & 1.6400 & 1.7510\\
Strong & 0.25 & 2000 & 0.3770 & 0.2220 & 0.3760 & \textbf{0.1680} & 0.3760 & 0.1680 & 4.4240 & 3.6360 & 0.9630 & 1.0070\\
\midrule
Strong & 1 & 500 & 1.9560 & 1.9060 & 1.9330 & \textbf{1.6120} & 1.9250 & 1.5970 & 8.0240 & 7.0790 & 3.4880 & 4.6260\\
Strong & 1 & 1000 & 1.3610 & 1.2310 & 1.3480 & \textbf{1.0540} & 1.3480 & 1.0540 & 6.2650 & 5.3490 & 2.2700 & 2.9650\\
Strong & 1 & 2000 & 0.9250 & 0.8190 & 0.9130 & \textbf{0.7240} & 0.9130 & 0.7240 & 4.8400 & 4.0250 & 1.4340 & 1.8600\\
\midrule
Weak & 0.25 & 500 & 0.6300 & 0.4150 & 0.6270 & \textbf{0.3150} & 0.6270 & 0.3150 & 6.0960 & 4.9690 & 1.5350 & 1.8540\\
Weak & 0.25 & 1000 & 0.4140 & 0.2630 & 0.4160 & \textbf{0.2130} & 0.4160 & 0.2130 & 4.8190 & 3.8330 & 0.9950 & 1.1210\\
Weak & 0.25 & 2000 & 0.2460 & 0.1610 & 0.2450 & \textbf{0.1260} & 0.2450 & 0.1260 & 3.7780 & 2.9290 & 0.6090 & 0.6230\\
\midrule
Weak & 1 & 500 & 1.3460 & 1.3340 & 1.3260 & \textbf{1.1950} & 1.3210 & 1.1830 & 6.9420 & 5.8130 & 2.4150 & 3.1900\\
Weak & 1 & 1000 & 0.9970 & 0.9470 & 0.9840 & \textbf{0.8370} & 0.9830 & 0.8360 & 5.4270 & 4.4080 & 1.6480 & 2.2440\\
Weak & 1 & 2000 & 0.6430 & 0.6540 & 0.6350 & \textbf{0.5790} & 0.6350 & 0.5790 & 4.1920 & 3.3140 & 1.0690 & 1.4380\\
\bottomrule
\end{tabular}%
}
\endgroup
\end{table}

\subsection{Case of spherical responses}\label{sec:simulation_sphere}

We next take $\M=\Sph^2$ with geodesic distance $d_{\Sph^2}(y,z):=\arccos(y\tran z)$ for $y,z\in\Sph^2$. The predictor distribution, functions $g_1$, $g_2$ and $g_{12}$, and strong and weak signal settings are the same as in Section~\ref{sec:simulation_scalar}. To generate spherical responses, for $\mx=(\mx_1\tran,\ldots,\mx_4\tran)\tran\in\Sph^2\times(\Sph^1)^3$ with $\mx_2=(\cos(\theta_2),\sin(\theta_2))\tran$, let $\mathbf e_3=(0,0,1)\tran\in\Sph^2$ and define
\begin{align*}
u_1(\mx):=0.38\tanh\left(\frac{1.8g_1(\mx_1)+\beta_2g_2(\theta_2)}{3}\right),\quad
u_2(\mx):=0.32\tanh\left(\frac{0.55g_1(\mx_1)+\beta_{12}g_{12}(\mx_1,\theta_2)}{1.8}\right).
\end{align*}
Writing $\mathbf v(\mx):=\left(u_1(\mx),u_2(\mx),0\right)\tran\in T_{\mathbf e_3}\Sph^2$, define
\begin{align*}
m_\oplus(\mx):=\cos\left(\|\mathbf v(\mx)\|\right)\mathbf e_3+\frac{\sin\left(\|\mathbf v(\mx)\|\right)}{\|\mathbf v(\mx)\|}\mathbf v(\mx),
\end{align*}
where $\sin t/t$ is defined as $1$ at $t=0$. Given $\mX_i$, we independently generate
\begin{align*}
Y_i\mid\mX_i\sim\operatorname{vMF}\left(m_\oplus(\mX_i),\kappa\right),\quad i=1,\ldots,n,
\end{align*}
where $\operatorname{vMF}(\mu,\kappa)$ denotes the von Mises--Fisher distribution on $\Sph^2$ with mean direction $\mu\in\Sph^2$ and concentration parameter $\kappa>0$.

We use $\kappa=400$ and $100$ for the low and high noise settings, respectively. The regression function is nonadditive, with $A_0=\{1,2\}$ in all settings. We use the same ten estimators, kernels, bandwidth grids and cross validation procedures as in Section~\ref{sec:simulation_scalar}. 

Table~\ref{tab:simulation_sphere_selection} reports factor selection frequencies. Both active factors are selected in every replication. Exact selection is at least $0.940$ across all settings and equals $1.000$ in every setting at $n=2000$. All remaining errors arise from including $\mX_3$ or $\mX_4$.

\begin{table}[!t]
\centering
\caption{Factor-selection frequencies for spherical responses based on $200$ Monte Carlo replications per setting. Exact gives the frequency of $\widehat A=A_0$ for $A_0=\{1,2\}$, and $\mX_1,\ldots,\mX_4$ give individual inclusion frequencies. Here $\kappa=400$ and $100$ correspond to low and high noise, respectively.}
\label{tab:simulation_sphere_selection}
\begin{tabular}{cccccccc}
\toprule
Signal & $\kappa$ & $n$ & Exact & $\mX_1$ & $\mX_2$ & $\mX_3$ & $\mX_4$\\
\midrule
Strong & 400 & 500 & 1.000 & 1.000 & 1.000 & 0.000 & 0.000\\
Strong & 400 & 1000 & 1.000 & 1.000 & 1.000 & 0.000 & 0.000\\
Strong & 400 & 2000 & 1.000 & 1.000 & 1.000 & 0.000 & 0.000\\
\midrule
Strong & 100 & 500 & 0.970 & 1.000 & 1.000 & 0.015 & 0.015\\
Strong & 100 & 1000 & 1.000 & 1.000 & 1.000 & 0.000 & 0.000\\
Strong & 100 & 2000 & 1.000 & 1.000 & 1.000 & 0.000 & 0.000\\
\midrule
Weak & 400 & 500 & 0.995 & 1.000 & 1.000 & 0.000 & 0.005\\
Weak & 400 & 1000 & 1.000 & 1.000 & 1.000 & 0.000 & 0.000\\
Weak & 400 & 2000 & 1.000 & 1.000 & 1.000 & 0.000 & 0.000\\
\midrule
Weak & 100 & 500 & 0.940 & 1.000 & 1.000 & 0.030 & 0.030\\
Weak & 100 & 1000 & 0.995 & 1.000 & 1.000 & 0.005 & 0.000\\
Weak & 100 & 2000 & 1.000 & 1.000 & 1.000 & 0.000 & 0.000\\
\bottomrule
\end{tabular}
\end{table}

Table~\ref{tab:simulation_sphere_mise} reports the MISEs. Selected LL has the smallest MISE among the data-driven estimators in every setting and is close to Oracle LL. Factor selection has a smaller effect on LC estimation. The proposed estimators also have smaller MISEs than the Euclidean and Riemannian estimators in every setting.

\begin{table}[!t]
\centering
\caption{MISEs for spherical responses based on $200$ Monte Carlo replications per setting. All values are multiplied by $10^3$. Full, Selected and Oracle use all four, the selected and the true active factors, respectively. Euclidean and Riemannian use all four factors. Bold values indicate the smallest MISE among the data-driven methods, excluding Oracle.}
\label{tab:simulation_sphere_mise}
\begingroup
\small
\resizebox{\linewidth}{!}{%
\begin{tabular}{ccc*{10}{c}}
\toprule
& & & \multicolumn{2}{c}{Full} & \multicolumn{2}{c}{Selected} & \multicolumn{2}{c}{Oracle} & \multicolumn{2}{c}{Riemannian} & \multicolumn{2}{c}{Euclidean}\\
\cmidrule(lr){4-5}\cmidrule(lr){6-7}\cmidrule(lr){8-9}\cmidrule(lr){10-11}\cmidrule(lr){12-13}
Signal & $\kappa$ & $n$ & LC & LL & LC & LL & LC & LL & LC & LL & LC & LL\\
\midrule
Strong & 400 & 500 & 1.901 & 1.540 & 1.894 & \textbf{1.268} & 1.894 & 1.268 & 9.464 & 8.493 & 3.695 & 4.536\\
Strong & 400 & 1000 & 1.200 & 1.054 & 1.196 & \textbf{0.885} & 1.196 & 0.885 & 7.475 & 6.484 & 2.393 & 2.873\\
Strong & 400 & 2000 & 0.776 & 0.650 & 0.773 & \textbf{0.544} & 0.773 & 0.544 & 5.813 & 4.932 & 1.447 & 1.783\\
\midrule
Strong & 100 & 500 & 2.962 & 3.097 & 2.924 & \textbf{2.687} & 2.918 & 2.679 & 10.830 & 9.613 & 4.930 & 6.625\\
Strong & 100 & 1000 & 2.073 & 1.942 & 2.058 & \textbf{1.707} & 2.058 & 1.707 & 8.449 & 7.277 & 3.345 & 4.311\\
Strong & 100 & 2000 & 1.460 & 1.286 & 1.451 & \textbf{1.134} & 1.451 & 1.134 & 6.511 & 5.527 & 2.130 & 2.811\\
\midrule
Weak & 400 & 500 & 1.266 & 1.180 & 1.258 & \textbf{1.011} & 1.258 & 1.010 & 8.682 & 7.317 & 2.390 & 2.884\\
Weak & 400 & 1000 & 0.848 & 0.763 & 0.843 & \textbf{0.651} & 0.843 & 0.651 & 6.861 & 5.584 & 1.592 & 1.967\\
Weak & 400 & 2000 & 0.591 & 0.483 & 0.590 & \textbf{0.424} & 0.590 & 0.424 & 5.351 & 4.248 & 1.016 & 1.261\\
\midrule
Weak & 100 & 500 & 2.035 & 2.140 & 2.014 & \textbf{1.942} & 2.007 & 1.931 & 10.047 & 8.440 & 3.613 & 4.377\\
Weak & 100 & 1000 & 1.534 & 1.482 & 1.512 & \textbf{1.312} & 1.512 & 1.312 & 7.837 & 6.388 & 2.476 & 3.287\\
Weak & 100 & 2000 & 1.025 & 1.051 & 1.014 & \textbf{0.935} & 1.014 & 0.935 & 6.048 & 4.836 & 1.651 & 2.266\\
\bottomrule
\end{tabular}%
}
\endgroup
\end{table}

\FloatBarrier

\section{Real-data analysis}
\label{sec:data}

\subsection{AIS data analysis}\label{subsec:data_ais}

We analyze Automatic Identification System (AIS) records for cargo vessels in the California coastal region during 2021--2024, obtained from the NOAA Pacific Marine Environmental Laboratory ERDDAP server \citep{NOAAPMELAIS} at \url{https://data.pmel.noaa.gov/pmel/erddap/}. The response $Y_i$ is speed over ground (SOG) in knots, and the predictors are geographic location, course over ground, local time of day and annual phase. We retain records with vessel-type codes $70$--$79$, valid positions and courses, and $0.5<\mathrm{SOG}\le40$, and randomly select one observation per Maritime Mobile Service Identity (MMSI), yielding $4613$ observations. Data from 2021--2022, 2023 and 2024 are used for training, validation and testing, respectively.

We represent geographic location on $\Sph^2$ and the remaining predictors on $\Sph^1$. For longitude $\lambda_i$ and latitude $\varphi_i$ in radians, define
\begin{align*}
\mX_{i1}:=\left(\cos(\varphi_i)\cos(\lambda_i),\cos(\varphi_i)\sin(\lambda_i),\sin(\varphi_i)\right)\tran,\quad
\mX_{i\ell}:=\left(\cos(\theta_{i\ell}),\sin(\theta_{i\ell})\right)\tran,\quad \ell=2,3,4,
\end{align*}
where $\theta_{i2}$ is course over ground and $\theta_{i3}$ and $\theta_{i4}$ are the daily and annual phases. Thus $\mX_i=(\mX_{i1}\tran,\ldots,\mX_{i4}\tran)\tran\in\Sph^2\times(\Sph^1)^3$. For the response, we define $Y_i:=\mathrm{SOG}_i\in\R$.

We compare the proposed Full and Selected LC and LL estimators with the Euclidean and Riemannian estimators from Section~\ref{sec:simulation_scalar}. We use the bandwidth grids in Section~\ref{sec:simulation_scalar}. 

Factor selection uses the 2021--2022 training sample with $5$-fold cross-fitting and $3$-fold inner cross validation, using $c=0.005$. Final bandwidths are selected by minimizing prediction error on the 2023 validation sample for estimators fitted to the 2021--2022 training data. 

\begin{figure}[!t]
\centering
\includegraphics[width=\textwidth]{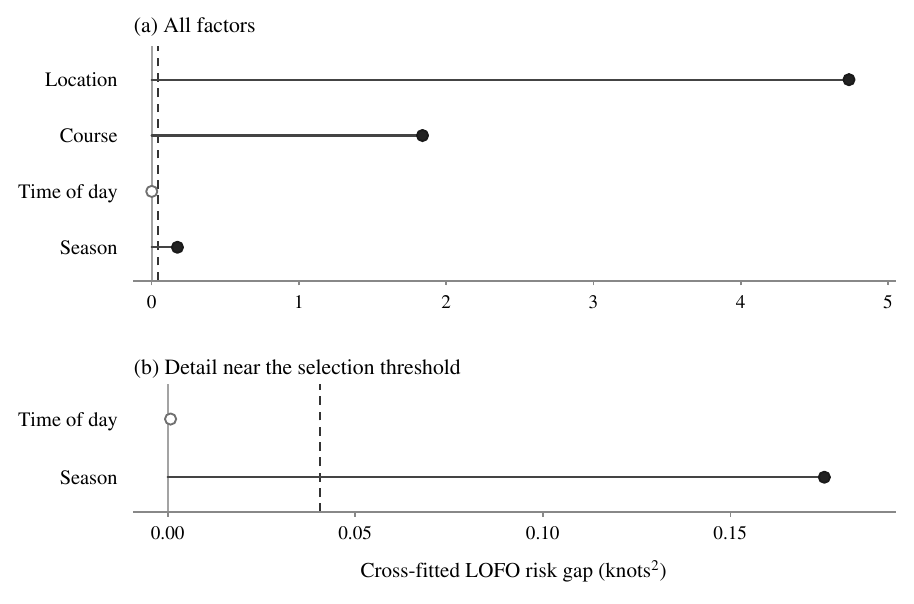}
\caption{Estimated LOFO risk gaps for the AIS training data. The dashed line indicates the selection threshold with $c=0.005$. Filled and open circles denote selected and excluded factors, respectively. Panels (a) and (b) show the gaps for all factors and an enlarged view of those for time of day and season, respectively.}
\label{fig:ais_lofo}
\end{figure}

Figure~\ref{fig:ais_lofo} shows the cross-fitted LOFO risk gaps. Location, course and season are selected, reducing the intrinsic predictor dimension from five to four. The selected set is unchanged over the eight equally spaced multipliers $c\in\{0.0025j:j=1,\ldots,8\}$. 

Table~\ref{tab:ais_prediction} reports RMSEs for the $1162$ observations in the 2024 test sample. Full LL has the smallest RMSE, and Selected LL has a similar RMSE with one fewer factor. The proposed estimators have smaller RMSEs than the corresponding Euclidean and Riemannian estimators.

\begin{table}[!t]
\centering
\caption{Root mean squared prediction errors (knots) for the AIS data in 2024. All values are multiplied by $10$. Full and Selected denote the proposed estimators using all four factors and the selected factors, respectively. Euclidean and Riemannian use all four factors. Bold indicates the smallest RMSE.}
\label{tab:ais_prediction}
\begingroup
\begin{tabular}{*{8}{c}}
\toprule
\multicolumn{2}{c}{Full} & \multicolumn{2}{c}{Selected} & \multicolumn{2}{c}{Riemannian} & \multicolumn{2}{c}{Euclidean}\\
\cmidrule(lr){1-2}\cmidrule(lr){3-4}\cmidrule(lr){5-6}\cmidrule(lr){7-8}
LC & LL & LC & LL & LC & LL & LC & LL\\
\midrule
34.383 & \textbf{34.026} & 34.380 & 34.107 & 39.019 & 36.688 & 35.807 & 34.721\\
\bottomrule
\end{tabular}
\endgroup
\end{table}

\begin{figure}[!t]
\centering
\makebox[\textwidth][c]{%
\includegraphics[width=1.5\textwidth]{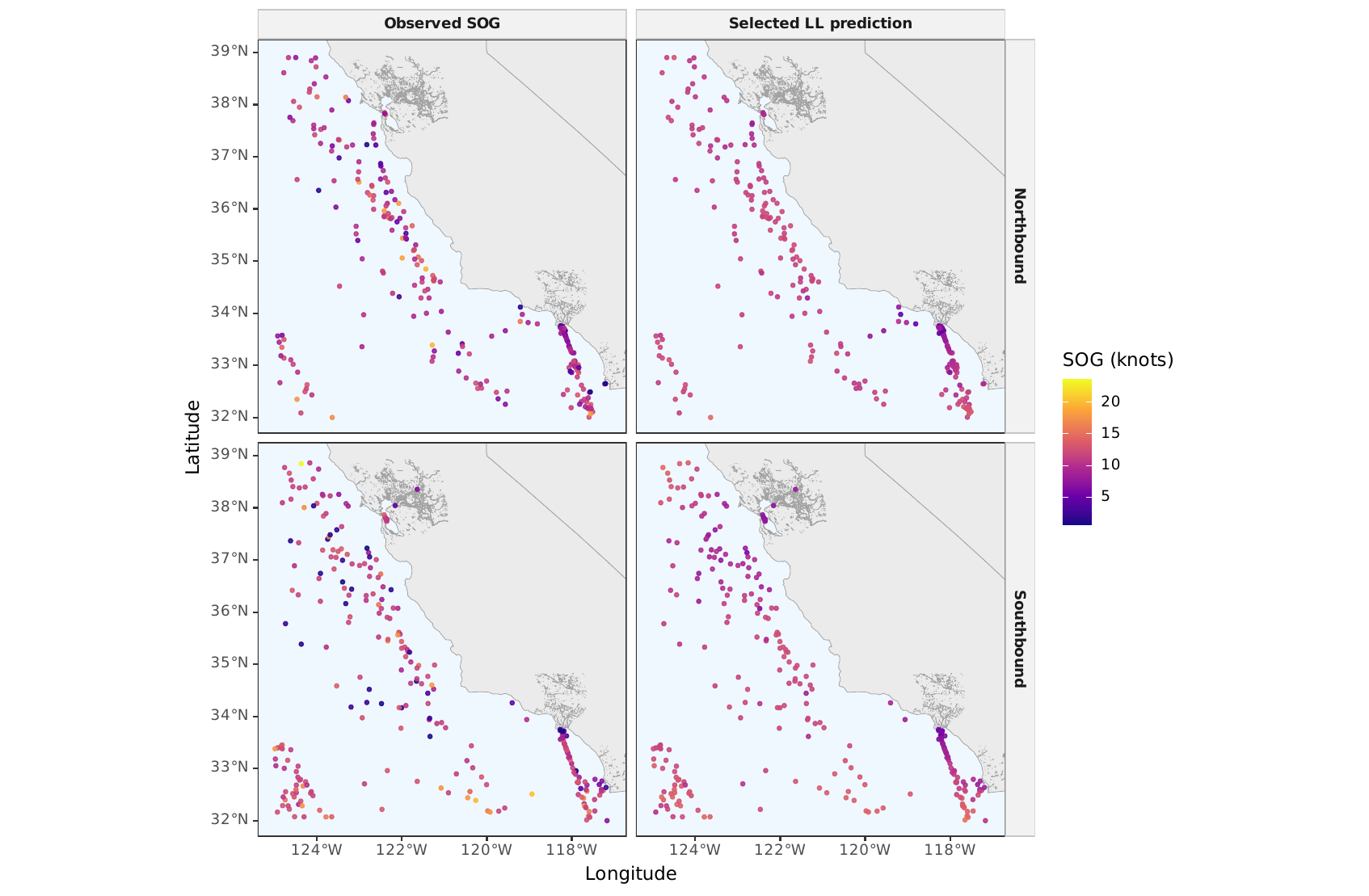}%
}
\caption{Observed SOG and Selected LL predictions (knots) at 2024 vessel locations. Northbound and southbound vessels have courses within $45^\circ$ of $0^\circ$ and $180^\circ$, respectively, giving $244$ and $263$ displayed observations. Each pair of panels uses the same locations, and all panels share a common color scale. The estimator uses the 2021--2022 training sample and bandwidths selected on the 2023 validation sample.}
\label{fig:ais_fitted}
\end{figure}

Figure~\ref{fig:ais_fitted} compares observed SOG values with Selected LL predictions at the same 2024 vessel locations for northbound and southbound vessels. The fitted values retain broad spatial variation while being less variable than the individual observed speeds.

\FloatBarrier

\subsection{GDP data analysis}\label{subsec:data_gdp}

We analyze ocean drifter observations from the Global Drifter Program (GDP) during 2021--2024, obtained from the UC San Diego Lagrangian Drifter Laboratory \citep{CenturioniBraaschSchonau2025} through \url{https://erddap.ldl.ucsd.edu/erddap/}. The response is the empirical distribution of displacement speeds over the following seven days, and the predictors are geographic location, recent heading, previous heading, annual phase and local solar time. We retain observations with two preceding and $28$ subsequent consecutive six-hour intervals, displacement speeds in $[0,3]$ meters per second, preceding speeds of at least $0.005$ meters per second, and valid heading positions within $85^\circ$ latitude. We randomly select one eligible observation per World Meteorological Organization (WMO) identifier and calendar year, yielding $4174$ observations from $2427$ distinct WMO identifiers in 2021--2023 and $1405$ observations in 2024. Data from 2021--2023 and 2024 are used for final estimation and testing, respectively.

We represent geographic location by $\mX_{i1}\in\Sph^2$ and recent heading, previous heading, annual phase and local solar time by $\mX_{i2},\ldots,\mX_{i5}\in\Sph^1$, respectively. The two headings are the bearings over the preceding six-hour intervals. The annual phase is the fraction of the UTC calendar year elapsed, multiplied by $2\pi$, and local solar time is UTC time in hours plus longitude in degrees divided by $15$, modulo $24$. Thus $\mX_i=(\mX_{i1}\tran,\ldots,\mX_{i5}\tran)\tran\in\Sph^2\times(\Sph^1)^4$, with intrinsic dimension six. To construct the response, let $s_{i,1},\ldots,s_{i,28}$ denote the subsequent displacement speeds and define
\begin{align*}
Y_i:=\frac1{28}\sum_{k=1}^{28}\delta_{s_{i,k}}.
\end{align*}

We compare the proposed Full and Selected LC and LL estimators with the Euclidean and Riemannian estimators from Section~\ref{sec:simulation_scalar} and two additional classes of estimators:
\begin{itemize}
\item Directional: LC and LL estimators using the proposed product von Mises kernel with a common bandwidth across factors.
\item Geodesic: LC and LL estimators using an isotropic Gaussian kernel based on the product geodesic distance, without volume correction.
\end{itemize}
Both additional LL estimators use the local linear coordinate construction in Section~\ref{sec:ll} with their respective kernels and bandwidths. The Euclidean estimators use a product Gaussian kernel applied to ordinary differences of longitude, latitude and angular coordinates. All competing estimators use all five factors. The candidate bandwidth grids are
\begin{align*}
\mathcal H_{[5]}^{\mathrm{GDP}}&:=\mathcal G_{20}(0.08,1.1)\times\mathcal G_{20}(0.25,5)^4,\quad \mathcal H_{\mathrm{GDP}}^{\mathrm{Euc}}:=\mathcal G_{20}(0.08,6)\times\mathcal G_{20}(0.25,6.5)^4,\\
\mathcal H_{\mathrm{GDP}}^{\mathrm{Rim}}&:=\mathcal G_{50}(0.7,3.12),\quad \mathcal H_{\mathrm{GDP}}^{\mathrm{Dir}}:=\mathcal G_{20}(0.3,5),\quad \mathcal H_{\mathrm{GDP}}^{\mathrm{Geo}}:=\mathcal G_{20}(0.4,5),
\end{align*}
where $\mathcal H_{[5]}^{\mathrm{GDP}}$ denotes the proposed bandwidth grid, and $\mathcal H_{\mathrm{GDP}}^{\mathrm{Euc}}$, $\mathcal H_{\mathrm{GDP}}^{\mathrm{Rim}}$, $\mathcal H_{\mathrm{GDP}}^{\mathrm{Dir}}$ and $\mathcal H_{\mathrm{GDP}}^{\mathrm{Geo}}$ denote the grids for the Euclidean, Riemannian, Directional and Geodesic estimators, respectively.

Factor selection uses the 2021--2023 estimation sample with $5$-fold cross-fitting and $3$-fold inner cross validation, using $c=0.005$ and $\widehat V_Y$ computed from pairwise squared Wasserstein distances as in \eqref{eq:empirical_response_dispersion}. Final bandwidths are selected by $5$-fold cross validation on the same estimation sample, with observations sharing a WMO identifier assigned to the same fold in each partition.

\begin{figure}[!t]
\centering
\includegraphics[width=\textwidth]{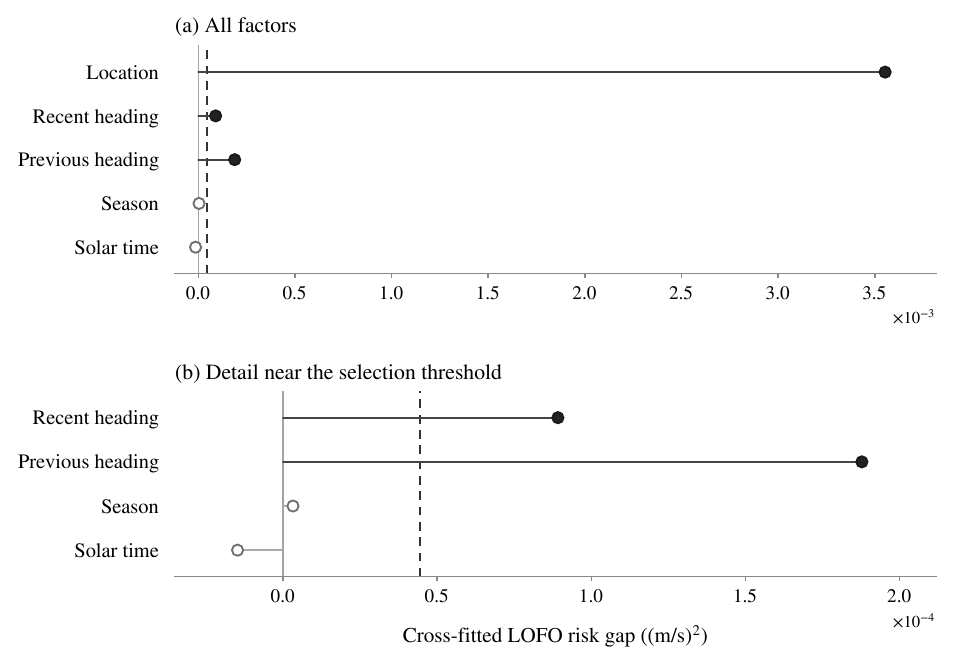}
\caption{Estimated LOFO risk gaps for the GDP estimation sample. The dashed line indicates the selection threshold with $c=0.005$. Filled and open circles denote selected and excluded factors, respectively. Panels (a) and (b) show the gaps for all factors and an enlarged view of those for the four circular predictors, respectively.}
\label{fig:gdp_lofo}
\end{figure}

Figure~\ref{fig:gdp_lofo} shows the cross-fitted LOFO risk gaps. Location, recent heading and previous heading are selected, reducing the intrinsic predictor dimension from six to four. The selected factors are unchanged from the original-grid analysis. All three factors remain selected for $c\le0.01$ on the eight-point grid used for AIS. For $c\ge0.0125$, recent heading is excluded, while location and previous heading remain selected.

Table~\ref{tab:gdp_prediction} reports Wasserstein prediction RMSEs. Full LL has the smallest RMSE, followed by Selected LL. Their corresponding MSEs differ by approximately $0.12\%$, while Selected LL uses three of the five factors. The proposed estimators also have smaller RMSEs than the competing estimators.

\begin{table}[!t]
\centering
\caption{Root mean squared Wasserstein prediction errors (meters per second, multiplied by $10^2$) for the GDP data in 2024. Full and Selected denote the proposed estimators using all five factors and the selected factors, respectively. Directional, Geodesic, Riemannian and Euclidean use all five factors. Bold indicates the smallest RMSE.}
\label{tab:gdp_prediction}
\begingroup
\small
\resizebox{\linewidth}{!}{%
\begin{tabular}{*{12}{c}}
\toprule
\multicolumn{2}{c}{Full} & \multicolumn{2}{c}{Selected} & \multicolumn{2}{c}{Directional} & \multicolumn{2}{c}{Geodesic} & \multicolumn{2}{c}{Riemannian} & \multicolumn{2}{c}{Euclidean}\\
\cmidrule(lr){1-2}\cmidrule(lr){3-4}\cmidrule(lr){5-6}\cmidrule(lr){7-8}\cmidrule(lr){9-10}\cmidrule(lr){11-12}
LC & LL & LC & LL & LC & LL & LC & LL & LC & LL & LC & LL\\
\midrule
14.756 & \textbf{14.721} & 14.779 & 14.730 & 15.932 & 15.856 & 15.942 & 15.900 & 16.032 & 15.979 & 14.980 & 15.187\\
\bottomrule
\end{tabular}%
}
\endgroup
\end{table}

\begin{figure}[!t]
\centering
\includegraphics[width=\textwidth]{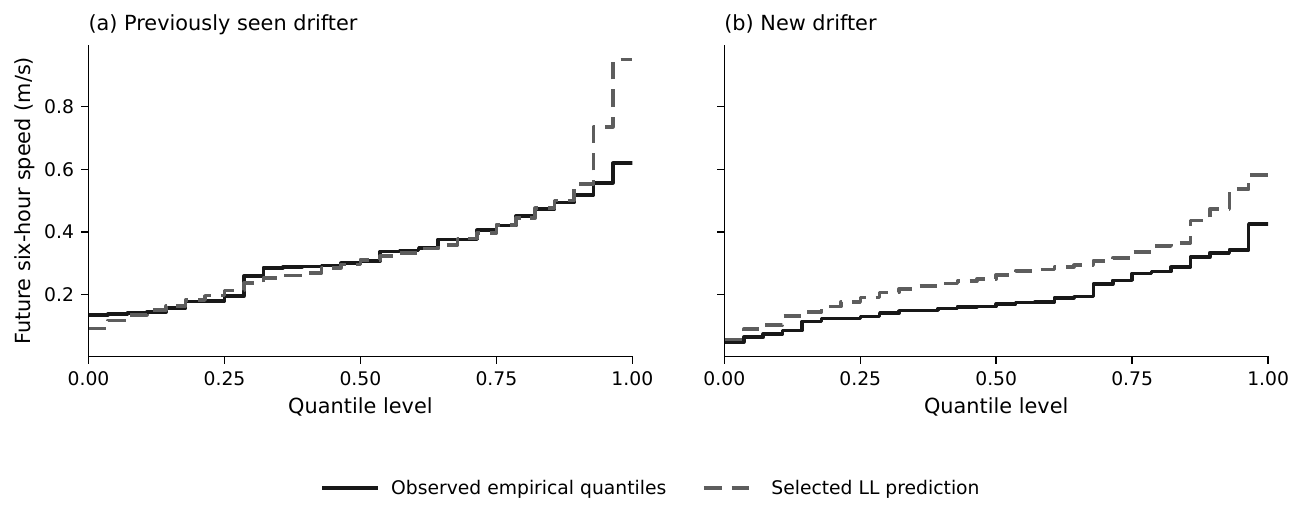}
\caption{Observed empirical quantile functions and Selected LL predictions for two 2024 drifters. WMO 1601735 occurs in the 2021--2023 estimation sample, whereas WMO 1801717 does not.}
\label{fig:gdp_fitted}
\end{figure}

Figure~\ref{fig:gdp_fitted} compares observed empirical quantile functions with Selected LL predictions for two 2024 drifters, one whose identifier occurs in the estimation sample and one whose identifier does not. The displayed observations have squared Wasserstein prediction losses closest to the respective group medians. The largest discrepancies occur in the upper quantiles.

\FloatBarrier

\section{Discussion}\label{sec:discussion}

We have developed local Fr\'echet regression with factor-specific bandwidths and factor selection for predictors on a product of spheres. The product structure gives rotation-invariant units for selection while allowing dependence and interactions among factors. The uniform convergence results are established for deterministic bandwidth vectors. Together with consistent factor selection, they yield uniform rates for the selected LC and LL estimators. Under the stated bandwidth and curvature conditions, these rates depend on the active factors' total intrinsic dimension. For Wasserstein responses on a compact interval, quantile projection gives a separate proof without an entropy condition on the response space.

In the simulations, factor selection improves local linear estimation and yields errors close to those of the oracle estimator. In the AIS and GDP analyses, fewer factors are retained while prediction errors remain close to those of the full estimator. The applications have restricted geographic support, and GDP observations from the same drifter may remain dependent across years. These features lie outside the independent-sampling theory with a positive density over the whole product sphere. Grouped cross validation prevents the same identifier from appearing in both parts of a tuning split but does not establish independence of the observations. The expanded-grid GDP comparison is also a post-test sensitivity analysis.

The active factor set is defined by preservation of optimal prediction risk, rather than by a collection of significance tests. The theory assumes a fixed number of factors and deterministic bandwidth sequences. Growing numbers of factors, factor-level inference, dependent observations and smoothness-adaptive bandwidth selection remain topics for further study.

\section*{Acknowledgements}
Chang Jun Im was supported by the National Research Foundation of Korea grant funded by the Korea government (MSIT) (No.\ RS-2025-00515381). Part of this research was conducted during the author's visit to the Institute for Mathematical and Statistical Innovation (IMSI), which is supported by the National Science Foundation under Grant DMS-2425650.

\clearpage
\setcounter{section}{0}
\setcounter{subsection}{0}
\setcounter{equation}{0}
\setcounter{table}{0}
\renewcommand{\thetable}{S.\arabic{table}}
\renewcommand{\theHtable}{supp.\arabic{table}}
\setcounter{theorem}{0}
\setcounter{lemma}{0}
\setcounter{corollary}{0}
\setcounter{proposition}{0}
\setcounter{definition}{0}
\setcounter{remark}{0}
\renewcommand{\thesection}{S.\arabic{section}}
\renewcommand{\theequation}{S.\arabic{section}.\arabic{equation}}
\renewcommand{\theHequation}{S.\arabic{section}.\arabic{equation}}
\renewcommand{\theHsection}{supp.\arabic{section}}
\renewcommand{\theHtheorem}{supp.\arabic{section}.\arabic{theorem}}
\renewcommand{\theHlemma}{supp.\arabic{section}.\arabic{lemma}}
\renewcommand{\theHcorollary}{supp.\arabic{section}.\arabic{corollary}}
\renewcommand{\theHproposition}{supp.\arabic{section}.\arabic{proposition}}
\renewcommand{\theHdefinition}{supp.\arabic{section}.\arabic{definition}}
\renewcommand{\theHremark}{supp.\arabic{section}.\arabic{remark}}

\begin{center}
{\Large\bf Supplement to\\[0.4em] ``Adaptive local Fr\'echet regression with polyspherical predictors''}\\[0.8em] {\large by Chang Jun Im}
\end{center}

\vspace{0.5em}

\noindent
This Supplement proves the results of the main paper. We first fix a factor set and a deterministic bandwidth sequence and obtain bounds uniform in the predictor. Since $p$ is fixed, a finite union gives simultaneous bounds over the factor sets and the two bandwidth vectors. Constants denoted by $C$ may change between displays. They may depend on the fixed factor dimensions, the kernel and the regularity constants, but not on the predictor or on sufficiently small bandwidths. For a measurable function $F$ of $(\mX,Y)$, write $PF:=\E[F(\mX,Y)]$ and $P_nF:=n^{-1}\sum_{i=1}^nF(\mX^{(i)},Y^{(i)})$. Predictor-only functions are regarded as functions constant in the response argument.

For $A=\varnothing$, the product surface measure is unit mass on the singleton predictor space, $v_{\mathbf t}(A)=1$ and $b_{\mathbf t}(A)=0$. The augmented regressor, design matrices and both equivalent weights are one. All empirical estimators on the empty model are the empirical Fr\'echet mean. This convention also applies to the fallback estimator, using only the training observations when a fold is omitted.

\section{Auxiliary results}\label{supp:auxiliary}
\setcounter{equation}{0}

In the proofs, $\mathbf t_A=(t_\ell:\ell\in A)$ may be any deterministic positive bandwidth sequence satisfying the displayed conditions for that model. For the main results, it is a subvector of $\mathbf t\in\{\mb,\mh\}$. For all sufficiently large $n$,
\begin{align}\label{eq:subset_bandwidth_orders}
v_{\mathbf t}(A)\ge v_{\mathbf t}([p]),\quad b_{\mathbf t}(A)\le b_{\mathbf t}([p]).
\end{align}
Thus Conditions~\ref{con:B1}--\ref{con:B2} imply the corresponding conditions with $v_{\mathbf t}(A)$ in place of $v_{\mathbf t}([p])$. Abbreviations such as $v=v_{\mathbf t}(A)$ and $\rho_n=(\log n/(nv))^{1/2}$ are local to the proofs in which they are introduced.

\subsection{Integration and directional moments}\label{supp:geometry}

For $\mx\in\Sph^d$, let $\mB_\mx\in\R^{(d+1)\times d}$ have orthonormal columns perpendicular to $\mx$, and put $\mP_\mx:=\mI_{d+1}-\mx\mx\tran$. For $\mz\notin\{\mx,-\mx\}$, write
\begin{align*}
\mz=\mx\cos\theta+\mB_\mx\bm\xi\sin\theta,\quad \theta\in(0,\pi),\quad \bm\xi\in\Sph^{d-1}.
\end{align*}
Then $\mP_\mx\mH_\mx(\mz)=\theta\mB_\mx\bm\xi$. At the two endpoints this projected vector is zero under the convention in Section~\ref{sec:ll}. Lemma~2 of \citet{GarciaPortuguesEtAl2013S}, followed by the substitution $u=\cos\theta$, gives
\begin{align}\label{eq:supp_sphere_integration}
\int_{\Sph^d}F(\mz)\dd\omega_d(\mz)=\int_0^\pi\int_{\Sph^{d-1}}F(\mx\cos\theta+\mB_\mx\bm\xi\sin\theta)\sin^{d-1}\theta\dd\omega_{d-1}(\bm\xi)\dd\theta
\end{align}
for every nonnegative measurable or integrable $F$. Surface measure on $\Sph^0=\{-1,1\}$ is counting measure, so the formula also holds for $d=1$. Applying \eqref{eq:supp_sphere_integration} factor by factor and using Tonelli's theorem gives the product formula. This is the same reduction used in the spherical and toroidal analyses of \citet{ImJeonPark2025S} and \citet{ImJeon2026S}.

For nonnegative integers $a$, $k\in\{1,2\}$ and positive integers $d$, define
\begin{align*}
c_{d,a,k}:=|\Sph^{d-1}|2^{(d+a)/2}\int_0^\infty L^k(r^2)r^{d+a-1}\dd r.
\end{align*}
For a nonempty $A$, write $\theta_\ell:=\arccos(\mx_\ell\tran\mz_\ell)$ and let $\bm\Phi_\ell$ denote the $\ell$th block of $\bm\Phi_{\mx_A,\mathbf t_A}(\mz_A)$.

\begin{lemma}\label{lem:moment_orders}
Suppose Condition~\ref{con:L1} holds. For every fixed multi-index $\mathbf j=(j_\ell:\ell\in A)$ of nonnegative integers and $k\in\{1,2\}$, as $\max_{\ell\in A}t_\ell\to0$,
\begin{align*}
\int\cL_{\mx_A,\mathbf t_A}(\mz_A)^k\prod_{\ell\in A}\theta_\ell^{j_\ell}\dd\omega_{\md_A}(\mz_A)=v_{\mathbf t}(A)^{1-k}\prod_{\ell\in A}t_\ell^{j_\ell}\left(\prod_{\ell\in A}c_{d_\ell,j_\ell,k}+o(1)\right).
\end{align*}
The remainder is independent of $\mx_A$. For a single factor and $a\in\{0,2\}$, the corresponding normalized remainder is $O(t_\ell^2)$. In particular,
\begin{align*}
\int\cL_{\mx_A,\mathbf t_A}\dd\omega_{\md_A}&=\prod_{\ell\in A}c_{d_\ell,0,1}+O(b_{\mathbf t}(A)),\\
\int\cL_{\mx_A,\mathbf t_A}\bm\Phi_{\mx_A,\mathbf t_A}\dd\omega_{\md_A}&=\zero,\\
\int\cL_{\mx_A,\mathbf t_A}\bm\Phi_\ell\bm\Phi_k\tran\dd\omega_{\md_A}&=\zero\quad(\ell\ne k),\\
\int\cL_{\mx_A,\mathbf t_A}\bm\Phi_\ell\bm\Phi_\ell\tran\dd\omega_{\md_A}&=\frac{c_{d_\ell,2,1}}{d_\ell}\prod_{j\in A\setminus\{\ell\}}c_{d_j,0,1}\mP_{\mx_\ell}+O(b_{\mathbf t}(A)).
\end{align*}
For fixed $\alpha\in(0,\pi)$ and fixed nonnegative integers $a,q$,
\begin{align*}
\sup_{\theta>\alpha}L^k\left(\frac{1-\cos\theta}{h^2}\right)(\theta/h)^a=o(h^{d+q}).
\end{align*}
The same order holds after integration with respect to surface measure on $\Sph^d$.
\end{lemma}

\begin{proof}
By \eqref{eq:supp_sphere_integration}, it suffices to calculate the one-factor integral. Fix $\alpha\in(0,\pi)$ and write $c_\alpha=1-\cos\alpha$. The substitution $u=(1-\cos\theta)/h^2$ gives
\begin{align}\label{eq:isotropic_change}
h^{-d-a}\int_0^\alpha L^k\left(\frac{1-\cos\theta}{h^2}\right)\theta^a\sin^{d-1}\theta\dd\theta=\int_0^{c_\alpha/h^2}L^k(u)\left(\frac{\arccos(1-h^2u)}h\right)^a(2u-h^2u^2)^{(d-2)/2}\dd u.
\end{align}
On this interval, $2-h^2u\ge1+\cos\alpha>0$ and $\arccos(1-h^2u)/h\le C\sqrt u$. The integrand is bounded by $CL^k(u)u^{(d+a-2)/2}$, including when $d=1$. The bound is integrable by Condition~\ref{con:L1}. Dominated convergence therefore gives the limit
\begin{align*}
2^{(d+a-2)/2}\int_0^\infty L^k(u)u^{(d+a-2)/2}\dd u=2^{(d+a)/2}\int_0^\infty L^k(r^2)r^{d+a-1}\dd r.
\end{align*}

To control the remaining interval, put $F_k(r)=L^k(r^2)$. For every positive integer $m$, monotonicity gives
\begin{align}\label{eq:isotropic_tail_decay}
R^mF_k(R)\le\frac{m}{1-2^{-m}}\int_{R/2}^{R}r^{m-1}F_k(r)\dd r=o(1).
\end{align}
Consequently, for any fixed nonnegative $a,q$,
\begin{align*}
\sup_{\theta>\alpha}L^k\left(\frac{1-\cos\theta}{h^2}\right)(\theta/h)^a\le\pi^ah^{-a}L^k(c_\alpha/h^2)=o(h^{d+q}),
\end{align*}
where \eqref{eq:isotropic_tail_decay} is used with $m=d+a+q$. The sphere has finite surface area, so the integrated assertion follows. Taking $q=0$ also makes the omitted part of the normalized moment $o(1)$.

For $a=0$ or $2$, on $0\le h^2u\le c_\alpha$,
\begin{align*}
\left|\left(\frac{\arccos(1-h^2u)}{h\sqrt{2u}}\right)^a(1-h^2u/2)^{(d-2)/2}-1\right|\le Ch^2u.
\end{align*}
The expression at $u=0$ is defined by continuity. Thus the difference between the integrand in \eqref{eq:isotropic_change} and its limit has integral at most
\begin{align*}
Ch^2\int_0^\infty L^k(u)u^{(d+a)/2}\dd u=O(h^2).
\end{align*}
The omitted part of the limiting integral is at most
\begin{align*}
\frac{h^2}{c_\alpha}\int_{c_\alpha/h^2}^\infty L^k(u)u^{(d+a)/2}\dd u=o(h^2).
\end{align*}
The fixed-angle remainder is $o(h^2)$ after normalization by the preceding tail bound with $q=2$. This proves the stated remainder. Multiplying the one-factor integrals, including their factors $t_\ell^{-kd_\ell}$, proves the product moment formula.

For the vector moments, rotational symmetry gives
\begin{align*}
\int_{\Sph^{d-1}}\bm\xi\dd\omega_{d-1}(\bm\xi)=\zero_d,\quad \int_{\Sph^{d-1}}\bm\xi\bm\xi\tran\dd\omega_{d-1}(\bm\xi)=\frac{|\Sph^{d-1}|}{d}\mI_d.
\end{align*}
Since $\bm\Phi_\ell=t_\ell^{-1}\theta_\ell\mB_{\mx_\ell}\bm\xi_\ell$ except on a surface-null set, each first moment is zero. Every cross block contains two separate first angular moments and is also zero. The diagonal block follows from $\mB_{\mx_\ell}\mB_{\mx_\ell}\tran=\mP_{\mx_\ell}$ and the moment with $a=2$. These calculations use product surface measure, not independence of the predictor factors.
\end{proof}

\subsection{Reduced conditional distributions}\label{supp:response}

For $A^c=[p]\setminus A$, define
\begin{align}\label{eq:marginal_conditional_density}
\begin{split}
f_A(\mx_A)&:=\int_{\Sph^{\md_{A^c}}}f(\mx_A,\mx_{A^c})\dd\omega_{\md_{A^c}}(\mx_{A^c}),\\
g_{A,y}(\mx_A)&:=\frac{1}{f_A(\mx_A)}\int_{\Sph^{\md_{A^c}}}g_y(\mx_A,\mx_{A^c})f(\mx_A,\mx_{A^c})\dd\omega_{\md_{A^c}}(\mx_{A^c}).
\end{split}
\end{align}
The arguments on the right retain the original factor order. We use these conditional-law versions for the reduced models. For $A=[p]$, these are $f$ and $g_y$. For $A=\varnothing$, set $f_A=1$ and $g_{A,y}=1$ for $P_Y$-almost every $y$. Homogeneous extensions for a subset normalize its retained factors. For $j=0,1,2$, write
\begin{align*}
J_j:=\E[D_\M(Y)^jG(Y)],\quad T_G(u):=\E[D_\M(Y)^2G(Y)\1\{D_\M(Y)>u\}],\quad u>0.
\end{align*}
These quantities are used only in this Supplement.

\begin{lemma}\label{lem:marginalization}
Under Conditions~\ref{con:D1}--\ref{con:D2}, the densities $f_A$ are bounded above and away from zero and have twice continuously differentiable homogeneous extensions. The functions $g_{A,y}$ are jointly measurable versions of the conditional density ratios, and
\begin{align*}
\sup_{A\subseteq[p]}\sup_{\mx_A}\left(g_{A,y}(\mx_A)+\|\nabla\bar g_{A,y}(\mx_A)\|+\|\nabla^2\bar g_{A,y}(\mx_A)\|_{\mathrm{op}}\right)\le CG(y)
\end{align*}
for $P_Y$-almost every $y$. Moreover, $J_j<\infty$, $T_G(u)\to0$ as $u\to\infty$, and
\begin{align*}
\sup_{A,\mx_A}\E[D_\M(Y)^j\mid\mX_A=\mx_A]&\le CJ_j,\quad j=0,1,2,\\
\sup_{A,\mx_A}\E[D_\M(Y)^2\1\{D_\M(Y)>u\}\mid\mX_A=\mx_A]&\le CT_G(u).
\end{align*}
The assertions involving $D_\M$ use the reference mean defined in the main paper.
\end{lemma}

\begin{proof}
Since $f$ is positive and continuous on the compact polysphere, $c\le f\le C$. Integration over the omitted factors gives the same bounds, with different positive constants, for $f_A$. Differentiation under this finite surface integral is justified by bounded derivatives of $\bar f$ on a product of compact annuli around the polysphere.

For a fixed $A$ and $y$, let $N_{A,y}$ be the numerator in \eqref{eq:marginal_conditional_density}. The product rule gives
\begin{align*}
\sup_{\mx_A}\left(|N_{A,y}|+\|\nabla\bar N_{A,y}\|+\|\nabla^2\bar N_{A,y}\|_{\mathrm{op}}\right)\le CG(y).
\end{align*}
For a scalar numerator $N$ and positive denominator $f_A$, the quotient derivatives are
\begin{align*}
\nabla(N/f_A)&=f_A^{-1}\nabla N-Nf_A^{-2}\nabla f_A,\\
\nabla^2(N/f_A)&=f_A^{-1}\nabla^2N-f_A^{-2}\left(\nabla N\nabla f_A\tran+\nabla f_A\nabla N\tran+N\nabla^2f_A\right)+2Nf_A^{-3}\nabla f_A\nabla f_A\tran.
\end{align*}
Applying these identities to the homogeneous extensions proves the derivative bound. All constants can be made common to the finitely many factor sets.

The joint law of $(\mX,Y)$ has density $f(\mx)g_y(\mx)$ with respect to $\omega_{\md}\otimes P_Y$. Integrating against a nonnegative product test function and then over the omitted factors proves the conditional-law identity in \eqref{eq:marginal_conditional_density}. Joint measurability follows from parameterized integration. Derivatives are measurable as limits of difference quotients, and the supremum defining $G$ can be taken over a fixed countable dense subset because the derivatives are continuous. In addition,
\begin{align*}
\int_{\Sph^{\md}}g_y(\mx)f(\mx)\dd\omega_{\md}(\mx)=1
\end{align*}
for $P_Y$-almost every $y$, and hence $G(y)\ge1$ on this set. This also proves the envelope assertion for $A=\varnothing$.

For a nonnegative measurable $\psi$, the conditional-law identity yields
\begin{align*}
\E[\psi(Y)\mid\mX_A=\mx_A]=\int_\M\psi(y)g_{A,y}(\mx_A)P_Y(\dd y)\le C\E[\psi(Y)G(Y)].
\end{align*}
Take $\psi=D_\M^j$ and $\psi=D_\M^2\1\{D_\M>u\}$. Since $D_\M\ge1$, Condition~\ref{con:D2} implies $J_j\le J_2<\infty$. Dominated convergence gives $T_G(u)\to0$.
\end{proof}

\begin{example}\label{ex:gaussian_density_envelope}
Consider $Y=\mu(\mX)+\sigma\varepsilon$, where $\varepsilon\sim N(0,1)$ is independent of $\mX$, $\sigma>0$, and $\mu$ has a twice continuously differentiable homogeneous extension. Suppose Condition~\ref{con:D1} holds. If $\varphi_\sigma$ is the $N(0,\sigma^2)$ density, then $p_Y(y)=\E[\varphi_\sigma(y-\mu(\mX))]>0$ and
\begin{align*}
g_y(\mx)=\frac{\varphi_\sigma(y-\mu(\mx))}{p_Y(y)}.
\end{align*}
On the polysphere,
\begin{align*}
\nabla\bar g_y&=g_y\frac{y-\mu}{\sigma^2}\nabla\bar\mu,\\
\nabla^2\bar g_y&=g_y\left(\left(\frac{(y-\mu)^2}{\sigma^4}-\frac1{\sigma^2}\right)\nabla\bar\mu(\nabla\bar\mu)\tran+\frac{y-\mu}{\sigma^2}\nabla^2\bar\mu\right).
\end{align*}
Let $M=\sup_\mx|\mu(\mx)|$. Bounded derivatives of $\bar\mu$ give
\begin{align*}
G(y)\le\frac{C(1+y^2)}{p_Y(y)}\exp\left(-\frac{((|y|-M)_+)^2}{2\sigma^2}\right).
\end{align*}
Here $m_{\oplus,\varnothing}=\E[Y]$ and $D_\M(y)\le C(1+|y|)$. Multiplication by $P_Y(\dd y)=p_Y(y)\dd y$ cancels the denominator, proving $\E[D_\M(Y)^2G(Y)]<\infty$. Thus Condition~\ref{con:D2} holds for every $y\in\R$ without a positive lower bound for $p_Y$.
\end{example}

\subsection{Taylor expansion and equivalent weights}\label{supp:design}\label{supp:reductions}

We use a Taylor expansion along product geodesics, as in the toroidal calculation of \citet{ImJeon2026S}. Let $F:\Sph^{\md_A}\to\R$ have a twice continuously differentiable homogeneous extension. For nonantipodal factor pairs, write $\mz_\ell=\mx_\ell\cos\theta_\ell+\mB_{\mx_\ell}\bm\xi_\ell\sin\theta_\ell$ and define
\begin{align*}
\bm\gamma_\ell(u):=\mx_\ell\cos(u\theta_\ell)+\mB_{\mx_\ell}\bm\xi_\ell\sin(u\theta_\ell),\quad 0\le u\le1.
\end{align*}
Homogeneity in each factor implies $\nabla_\ell\bar F(\bm\gamma(u))\tran\bm\gamma_\ell(u)=0$. Since $\bm\gamma_\ell''(u)=-\theta_\ell^2\bm\gamma_\ell(u)$, the chain rule gives
\begin{align*}
\frac{\dd^2}{\dd u^2}\bar F(\bm\gamma(u))&=\bm\gamma'(u)\tran\nabla^2\bar F(\bm\gamma(u))\bm\gamma'(u),\\
\|\bm\gamma'(u)\|^2&=\sum_{\ell\in A}\theta_\ell^2.
\end{align*}
Let $\bm\Phi_\ell(\mz_A)$ denote the $\ell$th block of $\bm\Phi_{\mx_A,\mathbf t_A}(\mz_A)$. Taylor's formula with integral remainder therefore yields
\begin{align}\label{eq:product_taylor_remainder}
\begin{split}
&F(\mz_A)-F(\mx_A)-\sum_{\ell\in A}t_\ell\nabla_\ell\bar F(\mx_A)\tran\bm\Phi_\ell(\mz_A)\\
&\qquad=\int_0^1(1-u)\bm\gamma'(u)\tran\nabla^2\bar F(\bm\gamma(u))\bm\gamma'(u)\dd u.
\end{split}
\end{align}
In particular, the absolute remainder is bounded by $\frac12\sup_{\mz_A}\|\nabla^2\bar F(\mz_A)\|_{\mathrm{op}}\sum_{\ell\in A}\theta_\ell^2$. The first-order bound is $|F(\mz_A)-F(\mx_A)|\le\sup\|\nabla\bar F\|(\sum_\ell\theta_\ell^2)^{1/2}$. For a Hilbert-valued $F$, the same calculation uses the first and second Fr\'echet derivatives and the Bochner integral remainder. The acceleration term is again zero by factorwise homogeneity. Factorwise antipodes have zero probability under Condition~\ref{con:D1}, so these formulas apply inside the population integrals below. The empirical kernel classes retain their prescribed endpoint values.

\begin{proposition}\label{prop:supp_basis_invariance}
For each predictor value, the coordinate design matrix is invertible if and only if the ambient design has rank $1+|\md_A|$. On this event the ambient and coordinate representations give the same local linear equivalent weights. These weights do not depend on the orthonormal bases of the factorwise tangent spaces. The assertions hold for both empirical and population designs.
\end{proposition}

\begin{proof}
Write $\mB=\mB_{\mx_A}$ and $\overline\mB=\operatorname{diag}(1,\mB)$. Then
\begin{align*}
\overline\mB\tran\overline\mB=\mI_{1+|\md_A|},\quad \overline\mB\overline\mB\tran=1\oplus\bigoplus_{\ell\in A}\mP_{\mx_\ell}.
\end{align*}
Every augmented regressor belongs to the image of $\overline\mB$. If $\mGamma$ is either design matrix and $\mQ=\overline\mB\tran\mGamma\overline\mB$, then
\begin{align*}
\mGamma=\overline\mB\mQ\overline\mB\tran,\quad \operatorname{rank}(\mGamma)=\operatorname{rank}(\mQ).
\end{align*}
When $\mQ$ is invertible, $\overline\mB\mQ^{-1}\overline\mB\tran$ satisfies the four Moore--Penrose identities, because $\mGamma$ and $\mQ$ are symmetric and $\overline\mB\tran\overline\mB=\mI$. Thus
\begin{align*}
\mGamma^\dagger=\overline\mB\mQ^{-1}\overline\mB\tran,\quad \overline\mB\tran\overline{\mathbf e}_{\zero}=(1,0,\ldots,0)\tran\in\R^{1+|\md_A|}.
\end{align*}
Substituting the first identity into the ambient weight proves \eqref{eq:ll_weight_coord}. Notice that the last vector is the coordinate intercept vector, not the ambient vector of dimension $1+|\md_A|+|A|$.

Any other choice of bases has the form $\mB_\ell\mathbf O_\ell$ with $\mathbf O_\ell$ orthogonal. Let $\mathbf O=\operatorname{diag}(1,\mathbf O_\ell:\ell\in A)$. The coordinate regressor changes from $\mathbf z$ to $\mathbf O\tran\mathbf z$, the coordinate matrix changes to $\mathbf O\tran\mQ\mathbf O$, and the coordinate intercept vector $\mathbf e$ satisfies $\mathbf O\mathbf e=\mathbf e$. Hence
\begin{align*}
\mathbf e\tran(\mathbf O\tran\mQ\mathbf O)^{-1}\mathbf O\tran\mathbf z=\mathbf e\tran\mathbf O\tran\mQ^{-1}\mathbf z=\mathbf e\tran\mQ^{-1}\mathbf z.
\end{align*}
This proves invariance without a globally continuous choice of bases.
\end{proof}

\begin{proof}[Verification of \Cref{rem:invariance_reduction}]
For $p=1$, the coordinate regressor is $(1,h^{-1}\mB_\mx\tran\mH_\mx(\mz))\tran$. It is obtained from the spherical regressor in \citet{ImJeonPark2025S} by an invertible change of slope coordinates. If $\mathbf S=\operatorname{diag}(1,\mathbf S_1)$ is invertible, replacing a regressor $\mathbf z$ by $\mathbf S\mathbf z$ changes its design $\mQ$ to $\mathbf S\mQ\mathbf S\tran$, and
\begin{align*}
\mathbf e\tran(\mathbf S\mQ\mathbf S\tran)^{-1}\mathbf S\mathbf z=\mathbf e\tran\mQ^{-1}\mathbf z.
\end{align*}
A common positive kernel normalization also cancels from the equivalent weight. When all $d_\ell=1$, choose $\mB_{\mx_\ell}=\mathbf R\mx_\ell$, where $\mathbf R$ rotates by $\pi/2$. Then $\mB_{\mx_\ell}\tran\mH_{\mx_\ell}(\mz_\ell)$ is the signed angular displacement in \citet{ImJeon2026S} away from antipodes. The same calculation proves the toroidal reduction.
\end{proof}

Put $c_{0,A}=\prod_{\ell\in A}c_{d_\ell,0,1}$, $c_{2,A,\ell}=(c_{d_\ell,2,1}/d_\ell)\prod_{j\in A\setminus\{\ell\}}c_{d_j,0,1}$ and
\begin{align*}
\mGamma_{A,0}(\mx_A):=c_{0,A}\oplus\bigoplus_{\ell\in A}c_{2,A,\ell}\mP_{\mx_\ell}.
\end{align*}

\begin{lemma}\label{lem:population_design}
Under Conditions~\ref{con:L1} and \ref{con:D1}, uniformly in $\mx_A$ as $\max_\ell t_\ell\to0$,
\begin{align*}
\widetilde\mu_{\mathbf t_A}(\mx_A)&=c_{0,A}f_A(\mx_A)+O(b_{\mathbf t}(A)),\\
\widetilde\mGamma_{\mathbf t_A}(\mx_A)&=f_A(\mx_A)\mGamma_{A,0}(\mx_A)+O\left(\sum_{\ell\in A}t_\ell\right).
\end{align*}
For all sufficiently small bandwidths, the population design has rank $1+|\md_A|$ and its nonzero eigenvalues are bounded away from zero uniformly in $\mx_A$. In addition, $\sup_{\mx_A}\|\widetilde\mGamma_{\mathbf t_A}(\mx_A)^\dagger\|_{\mathrm{op}}\le C$ and, for $s=0,1$,
\begin{align*}
\E[\widetilde W_{\mx_A,\mathbf t_A,s}(\mX_A)]=1,\quad \sup_{\mx_A}\E[|\widetilde W_{\mx_A,\mathbf t_A,s}(\mX_A)|]\le C.
\end{align*}
The local linear weights reproduce the first-order tangent terms. On their design events the empirical weights satisfy the corresponding sample-average identities.
\end{lemma}

\begin{proof}
The properties of $f_A$ used here follow from Condition~\ref{con:D1} by integration over the omitted factors and do not require Condition~\ref{con:D2}. By \eqref{eq:product_taylor_remainder}, the integral of the linear term in $f_A$ is zero under the symmetric product kernel. Lemma~\ref{lem:moment_orders} bounds the integral of the remainder by $Cb_{\mathbf t}(A)$, proving the denominator expansion.

For a design entry, replacing $f_A(\mz_A)$ by $f_A(\mx_A)$ introduces a remainder bounded by
\begin{align*}
C\int\cL_{\mx_A,\mathbf t_A}(\mz_A)(1+\|\bm\Phi_{\mx_A,\mathbf t_A}(\mz_A)\|^2)\left(\sum_{\ell\in A}\theta_\ell^2\right)^{1/2}\dd\omega_{\md_A}(\mz_A)\le C\sum_{\ell\in A}t_\ell.
\end{align*}
The last inequality follows by expanding the finite sum and applying Lemma~\ref{lem:moment_orders} to moments through degree three. Its constant does not involve a bandwidth ratio. The remaining integral has the block structure in that lemma, with remainder $O(b_{\mathbf t}(A))$. This proves the matrix expansion.

Let $\mathbf a\in\overline{\cT}_{\mx_A}$ satisfy $\|\mathbf a\|=1$. All coefficients in $\mGamma_{A,0}$ are positive, so
\begin{align*}
\mathbf a\tran\widetilde\mGamma_{\mathbf t_A}(\mx_A)\mathbf a\ge\inf_{\mx_A}f_A(\mx_A)\min(c_{0,A},\min_{\ell\in A}c_{2,A,\ell})-C\sum_{\ell\in A}t_\ell\ge c>0
\end{align*}
for small bandwidths. Both designs vanish on the orthogonal complement of $\overline{\cT}_{\mx_A}$, which consists of the factorwise isotropic directions. Proposition~\ref{prop:supp_basis_invariance} now gives the rank and inverse bounds.

Write $\mathbf Z=\overline{\bm\Phi}_{\mx_A,\mathbf t_A}(\mX_A)$ and $\mathbf e=\overline{\mathbf e}_{\zero}$. Since $\mathbf Z\tran\mathbf e=1$,
\begin{align*}
\widetilde\mGamma\mathbf e=\E[\cL\mathbf Z],\quad \E[\widetilde W_1\mathbf Z]=\widetilde\mGamma\widetilde\mGamma^\dagger\mathbf e=\mathbf e.
\end{align*}
Thus
\begin{align}\label{eq:population_reproduction}
\E[\widetilde W_1]=1,\quad \E[\widetilde W_1\bm\Phi_{\mx_A,\mathbf t_A}(\mX_A)]=\zero.
\end{align}
Replacing expectations by sample averages gives, on the rank event,
\begin{align*}
\frac1n\sum_{i=1}^n\widehat W_1(\mX_A^{(i)})=1,\quad \frac1n\sum_{i=1}^n\widehat W_1(\mX_A^{(i)})\bm\Phi_{\mx_A,\mathbf t_A}(\mX_A^{(i)})=\zero.
\end{align*}
For LC, normalization follows directly from the denominator.

The inverse bound also gives
\begin{align}\label{eq:population_weight_envelope}
|\widetilde W_{\mx_A,\mathbf t_A,1}(\mz_A)|\le C\cL_{\mx_A,\mathbf t_A}(\mz_A)(1+\|\bm\Phi_{\mx_A,\mathbf t_A}(\mz_A)\|).
\end{align}
Lemma~\ref{lem:moment_orders} and boundedness of $f_A$ give the absolute first-moment bound. The LC weight is nonnegative and has expectation one.
\end{proof}

\subsection{Kernel covering bounds}\label{supp:kernel_entropy}

We establish covering bounds for $v_{\mathbf t}(A)\cL_{\mx_A,\mathbf t_A}$ and its tangent multipliers. A class $\cF$ with envelope $F$ is of VC type if, for some constants $C,V<\infty$,
\begin{align*}
\sup_Q N(\epsilon\|F\|_{Q,2},\cF,L_2(Q))\le(C/\epsilon)^V,\quad 0<\epsilon\le1,
\end{align*}
where the supremum is over probability measures $Q$.

\begin{lemma}\label{lem:kernel_covering}
Under Condition~\ref{con:L1}, the classes consisting of $v_{\mathbf t}(A)\cL_{\mx_A,\mathbf t_A}$ multiplied by coordinate monomials of degree at most two in $\bm\Phi_{\mx_A,\mathbf t_A}$ have finite constant envelopes and are of VC type. The same assertion holds for $(v_{\mathbf t}(A)\cL_{\mx_A,\mathbf t_A})^2$ and coordinate monomials of degree at most four. Predictor values and all positive bandwidths may vary. The covering constants depend only on the kernel and fixed dimensions.
\end{lemma}

\begin{proof}
For $\mx\in\Sph^d$ and $h>0$, define $r_{\mx,h}:\Sph^d\to[0,\infty)$ by $r_{\mx,h}(\mz):=\|\mz-\mx\|/h$. Let $\mathbf e_1,\ldots,\mathbf e_{d+1}$ be the standard basis of $\R^{d+1}$. Define $\phi_0,\ldots,\phi_{d+1}:\Sph^d\to\R$ by $\phi_0(\mz):=1$ and $\phi_a(\mz):=\mathbf e_a\tran\mz$ for $a=1,\ldots,d+1$. Then
\begin{align*}
r_{\mx,h}^2(\mz)=\frac{2-2\mx\tran\mz}{h^2}=\frac{2}{h^2}\phi_0(\mz)-\frac{2}{h^2}\sum_{a=1}^{d+1}(\mathbf e_a\tran\mx)\phi_a(\mz),\quad \mz\in\Sph^d.
\end{align*}
Thus, the functions $r_{\mx,h}^2$ belong to the $(d+2)$-dimensional vector space spanned by $\phi_0,\ldots,\phi_{d+1}$. By Lemmas~2.6.15 and~2.6.18(viii) and Theorem~2.6.7 of \citet{vanDerVaartWellner1996S}, for any fixed bounded monotone function $g:[0,\infty)\to\R$, there exist constants $C,V<\infty$ such that
\begin{align*}
\sup_Q N\left(\epsilon U_g,\{g\circ r_{\mx,h}:\mx\in\Sph^d,\ h>0\},L_2(Q)\right)\le(C/\epsilon)^V,\quad 0<\epsilon\le1,
\end{align*}
where $U_g:=\max\{1,\sup_{r\ge0}|g(r)|\}$ and the supremum is taken over probability measures $Q$ on $\Sph^d$.

For $k\in\{1,2\}$ and $j=0,\ldots,4$, write $F_k(r):=L^k(r^2/2)$ and $\psi_{j,k}(r):=r^jF_k(r)$. The function $\psi_{0,k}$ is bounded and nonincreasing. For $j\ge1$, define
\begin{align*}
\psi_{j,k}^{+}(r):=j\int_0^r u^{j-1}F_k(u)\dd u,\quad \psi_{j,k}^{-}(r):=\psi_{j,k}^{+}(r)-\psi_{j,k}(r).
\end{align*}
Both functions are bounded and nondecreasing. Indeed, Condition~\ref{con:L1} gives
\begin{align*}
0\le\psi_{j,k}^{-}(r)\le\psi_{j,k}^{+}(r)\le j\int_0^\infty u^{j-1}F_k(u)\dd u<\infty,
\end{align*}
and, for $0\le r_1<r_2$,
\begin{align*}
\psi_{j,k}^{-}(r_2)-\psi_{j,k}^{-}(r_1)\ge r_1^j\left(F_k(r_1)-F_k(r_2)\right)\ge0.
\end{align*}
Thus $\psi_{j,k}=\psi_{j,k}^{+}-\psi_{j,k}^{-}$. By Theorem~2.6.7 of \citet{vanDerVaartWellner1996S} and the triangle inequality, the classes $\{\psi_{j,k}\circ r_{\mx,h}:\mx\in\Sph^d,\ h>0\}$ have polynomial $L_2(Q)$ covering bounds with constant envelopes.

Write $u=\mx\tran\mz$ and define
\begin{align*}
\boldsymbol\eta_{\mx}(\mz):=\begin{cases}\displaystyle\frac{\mz-u\mx}{\sqrt{1-u^2}},&|u|<1,\\\zero_{d+1},&|u|=1,\end{cases}\quad q(s):=\begin{cases}\displaystyle\frac{2\arcsin(s/2)}s,&0<s\le2,\\1,&s=0.\end{cases}
\end{align*}
The convexity of $\arcsin$ on $[0,1]$ implies that $q$ is nondecreasing, with $1\le q(s)\le\pi/2$. Hence $\{q(r_{\mx,1})^j:\mx\in\Sph^d\}$ is a bounded VC-subgraph class. Also, $\|\boldsymbol\eta_{\mx}(\mz)\|\le1$. For an ambient coordinate $a$, put $w_a=z_a-ux_a$. The subgraph inequality $\tau<[\boldsymbol\eta_{\mx}(\mz)]_a$ is determined by the signs of
\begin{align*}
1-u^2,\quad w_a,\quad \tau,\quad w_a^2-\tau^2(1-u^2),
\end{align*}
including the case $|u|=1$. These are polynomials of bounded degree in $(\mz,\tau)$. By Lemmas~2.6.15 and~2.6.17 of \citet{vanDerVaartWellner1996S}, the corresponding sign classes and their finite unions and intersections are VC. Thus the coordinate classes of $\boldsymbol\eta_{\mx}$ are bounded VC-subgraph classes. Theorem~2.6.7 of the same reference gives polynomial $L_2(Q)$ covering bounds for these classes and for the classes generated by $q^j$.

Since $\|\mz-\mx\|=2\sin(\arccos(u)/2)$,
\begin{align*}
&L^k\left(\frac{1-\mx\tran\mz}{h^2}\right)\prod_{\nu=1}^j[h^{-1}\mP_\mx\mH_\mx(\mz)]_{a_\nu}\\
&\quad=\psi_{j,k}(r_{\mx,h}(\mz))q(r_{\mx,1}(\mz))^j\prod_{\nu=1}^j[\boldsymbol\eta_\mx(\mz)]_{a_\nu},
\end{align*}
including $j=0$ and the prescribed values at $\mz=\pm\mx$. For uniformly bounded functions,
\begin{align*}
\left\|\prod_{\nu=1}^m f_\nu-\prod_{\nu=1}^m g_\nu\right\|_{Q,2}\le\sum_{\nu=1}^m\|f_\nu-g_\nu\|_{Q,2}\prod_{\mu\ne\nu}\max(\|f_\mu\|_\infty,\|g_\mu\|_\infty).
\end{align*}
Combining finite covers with this inequality preserves polynomial covering bounds. Apply this argument to the displayed factors and then to the sphere factors. Taking $k=1$ with total monomial degree at most two and $k=2$ with total monomial degree at most four proves the lemma.
\end{proof}

We use inequality (2.2) in Theorem~2.1 of \citet{GineGuillou2002S}. For a measurable VC-type class $\cF$ with constant envelope $U$, variance at most $\sigma^2$ and covering bound $(\mathcal A/\epsilon)^V$, where $0<\sigma\le U$ and $\mathcal A\ge3\sqrt e$,
\begin{align}\label{eq:vc_expectation_bound}
\E\left[\sup_{F\in\cF}|(P_n-P)F|\right]\le C\left(\sqrt{\frac{V\sigma^2\log(\mathcal A U/\sigma)}n}+\frac{VU\log(\mathcal A U/\sigma)}n\right).
\end{align}
The kernel classes above satisfy the required measurability condition.

\subsection{Empirical designs and response moments}

\begin{lemma}\label{lem:empirical_design}
Suppose Conditions~\ref{con:L1}, \ref{con:D1} and the fixed-model version of Condition~\ref{con:B1} hold. Put $\rho_n=(\log n/(nv_{\mathbf t}(A)))^{1/2}$. Then
\begin{align*}
\sup_{\mx_A}\|\widehat\mGamma_{\mathbf t_A}(\mx_A)-\widetilde\mGamma_{\mathbf t_A}(\mx_A)\|_{\mathrm{op}}=O_\P(\rho_n).
\end{align*}
With probability tending to one, all LC denominators are positive and all LL designs have rank $1+|\md_A|$. On this event,
\begin{align*}
\sup_{\mx_A}\|\widehat\mGamma_{\mathbf t_A}(\mx_A)^\dagger-\widetilde\mGamma_{\mathbf t_A}(\mx_A)^\dagger\|_{\mathrm{op}}&=O_\P(\rho_n),\\
|\widehat W_{\mx_A,\mathbf t_A,1}(\mz_A)-\widetilde W_{\mx_A,\mathbf t_A,1}(\mz_A)|&\le C_n\rho_n\cL_{\mx_A,\mathbf t_A}(\mz_A)(1+\|\bm\Phi_{\mx_A,\mathbf t_A}(\mz_A)\|),
\end{align*}
where $C_n=O_\P(1)$ is common to both predictor arguments. For LC the same bound holds without the tangent factor. Also $\sup_{\mx_A}P_n|\widehat W_{\mx_A,\mathbf t_A,s}|=O_\P(1)$ for $s=0,1$.
\end{lemma}

\begin{proof}
Write $v=v_{\mathbf t}(A)$ and $\mathbf Z=\overline{\bm\Phi}_{\mx_A,\mathbf t_A}(\mX_A)$. Each design entry is the empirical average of $\cL$ times an augmented coordinate monomial of degree at most two. Lemmas~\ref{lem:moment_orders} and \ref{lem:kernel_covering} give a constant envelope $U=C_0/v$ and a second-moment bound $\sigma^2=C_1/v$ for the corresponding class. Enlarge the constants to ensure $\sigma\le U$. Since $nv/\log n\to\infty$, $v^{-1}=O(n)$ and
\begin{align*}
\log(\mathcal A U/\sigma)=O(\log n).
\end{align*}
By \eqref{eq:vc_expectation_bound},
\begin{align*}
\E\left[\sup_{\mx_A}|(P_n-P)(\cL Z_aZ_b)|\right]\le C\left(\sqrt{\frac{\log n}{nv}}+\frac{\log n}{nv}\right)=O(\rho_n).
\end{align*}
Markov's inequality and the fixed matrix dimension prove the operator-norm bound.

By Lemma~\ref{lem:population_design}, the restricted population eigenvalues have a common positive lower bound. The empirical perturbation tends to zero uniformly, and the two matrices vanish on the same isotropic space. Their restrictions to $\overline{\cT}_{\mx_A}$ are therefore simultaneously positive definite with probability tending to one, with inverse norms uniformly bounded. On this event the usual inverse identity on that subspace gives
\begin{align*}
\widehat\mGamma^\dagger-\widetilde\mGamma^\dagger=\widehat\mGamma^\dagger(\widetilde\mGamma-\widehat\mGamma)\widetilde\mGamma^\dagger.
\end{align*}
It follows that the inverse difference has order $O_\P(\rho_n)$. The intercept entry of the design is $\widehat\mu_{\mathbf t_A}$, which proves the LC assertion.

Subtracting the LL weight formulas and using $\|\overline{\bm\Phi}\|\le1+\|\bm\Phi\|$ gives the asserted weight difference. For LC use
\begin{align*}
\widehat\mu^{-1}-\widetilde\mu^{-1}=\frac{\widetilde\mu-\widehat\mu}{\widehat\mu\widetilde\mu}.
\end{align*}
The empirical zeroth moment and the trace of the tangent block give $\sup_{\mx_A}P_n\cL=O_\P(1)$ and $\sup_{\mx_A}P_n(\cL\|\bm\Phi\|^2)=O_\P(1)$. Hence
\begin{align*}
P_n(\cL\|\bm\Phi\|)\le\left(P_n\cL\,P_n(\cL\|\bm\Phi\|^2)\right)^{1/2}=O_\P(1)
\end{align*}
uniformly in the predictor. Combining this inequality with the empirical inverse bound proves the last assertion. LC has the exact sample-average normalization $P_n\widehat W_0=1$.
\end{proof}

\begin{proposition}\label{prop:TB}
Suppose Condition~\ref{con:B1} holds and put $v_n=\min(v_{\mb}([p]),v_{\mh}([p]))$. Each of the following is sufficient for Condition~\ref{con:B2}.
\begin{enumerate}[label=(\roman*),leftmargin=2.2em]
\item $D_\M(Y)$ is almost surely bounded.
\item For some $q>2$, $\E[D_\M(Y)^q]<\infty$ and $\liminf_n n^{1-2/q}v_n/\log n>0$.
\item For some $q>2$ and $C<\infty$, $\P(D_\M(Y)>u)\le Cu^{-q}$ for all sufficiently large $u$, and the bandwidth condition in (ii) holds.
\item For some $C,c,\alpha>0$, $\P(D_\M(Y)>u)\le C\exp(-cu^\alpha)$ for all sufficiently large $u$, and
\begin{align*}
\liminf_{n\to\infty}\frac{(nv_n/\log n)^{\alpha/2}}{\log n}>0.
\end{align*}
\end{enumerate}
For a fixed model of intrinsic dimension $d_A=|\md_A|$ and bandwidths $t_\ell\asymp(\log n/n)^{1/(d_A+4)}$, a moment of order $q>\max(2,(d_A+4)/2)$ suffices for its fixed-model tail condition. This last assertion concerns $v_{\mathbf t}(A)$, not the product over omitted factors.
\end{proposition}

\begin{proof}
Let $T_n=\sqrt{nv_n/\log n}$ and $p_n(C)=\P(D_\M(Y)>CT_n)$. Independence gives
\begin{align*}
1-\exp(-np_n(C))\le\P\left(\max_{1\le i\le n}D_\M(Y^{(i)})>CT_n\right)=1-(1-p_n(C))^n\le np_n(C).
\end{align*}
Thus the maximal-tail condition is equivalent to $\lim_{C\to\infty}\limsup_{n\to\infty}np_n(C)=0$. Since $T_n\to\infty$, bounded responses satisfy it. Under (ii), Markov's inequality gives
\begin{align*}
np_n(C)\le C^{-q}\E[D_\M(Y)^q]\left(\frac{\log n}{n^{1-2/q}v_n}\right)^{q/2}.
\end{align*}
First taking the upper limit in $n$ and then sending $C\to\infty$ proves (ii). The same calculation with the polynomial-tail bound proves (iii). If the ratio in (ii) tends to infinity, the bound vanishes for every fixed $C>0$, giving the stronger order $o_\P(T_n)$. Under (iv), $np_n(C)\le C_1n\exp(-cC^\alpha T_n^\alpha)$, which tends to zero for sufficiently large $C$.

For the stated fixed-model order, $v_{\mathbf t}(A)\asymp(\log n/n)^{d_A/(d_A+4)}$. Therefore
\begin{align*}
\frac{n^{1-2/q}v_{\mathbf t}(A)}{\log n}\asymp n^{4/(d_A+4)-2/q}(\log n)^{-4/(d_A+4)}\longrightarrow\infty
\end{align*}
when $q>(d_A+4)/2$. The condition $q>2$ ensures applicability of (ii). These moment restrictions are sufficient, not necessary.
\end{proof}

\begin{lemma}\label{lem:local_weighted_moments}
Suppose Conditions~\ref{con:L1}, \ref{con:D1}--\ref{con:D2}, \ref{con:M1} and the fixed-model versions of Conditions~\ref{con:B1}--\ref{con:B2} hold. Write $\cL_i=\cL_{\mx_A,\mathbf t_A}(\mX_A^{(i)})$, $\bm\Phi_i=\bm\Phi_{\mx_A,\mathbf t_A}(\mX_A^{(i)})$ and $D_i=D_\M(Y^{(i)})$. For $j=0,1,2$,
\begin{align*}
\sup_{\mx_A}\frac1n\sum_{i=1}^n\cL_i(1+\|\bm\Phi_i\|^j)D_i^2&=O_\P(1),\\
\sup_{\mx_A}\frac{v_{\mathbf t}(A)}n\sum_{i=1}^n\cL_i^2(1+\|\bm\Phi_i\|^{2j})D_i^2&=O_\P(1).
\end{align*}
These statements also hold with $D_i$ in place of $D_i^2$. For $s=0,1$ and $T\ge1$,
\begin{align*}
\sup_{\mx_A}\E[\widetilde W_{\mx_A,\mathbf t_A,s}(\mX_A)^2D_\M(Y)^2]&\le C/v_{\mathbf t}(A),\\
\sup_{\mx_A}\E[|\widetilde W_{\mx_A,\mathbf t_A,s}(\mX_A)|D_\M(Y)\1\{D_\M(Y)>T\}]&\le CT_G(T)/T.
\end{align*}
\end{lemma}

\begin{proof}
Put $v=v_{\mathbf t}(A)$ and $\rho_n=(\log n/(nv))^{1/2}$. It suffices to use nonnegative finite sums of absolute coordinate monomials that dominate $1+\|\bm\Phi\|^j$ and $1+\|\bm\Phi\|^{2j}$. Their orders and constant envelopes are those in Lemmas~\ref{lem:moment_orders} and \ref{lem:kernel_covering}. Let $w$ be either resulting kernel weight, including the factor $v$ for the squared kernel. Then
\begin{align*}
0\le w\le C/v,\quad Pw\le C,\quad w^2\le(C/v)w,\quad P[wD_\M^2]\le CJ_2.
\end{align*}
The last inequality follows by conditioning on $\mX_A$ and using Lemma~\ref{lem:marginalization}. The same conditioning gives $P[wD_\M^2\1\{D_\M>T\}]\le CT_G(T)$.

Fix $C_0>0$ and truncate at $T_n=C_0/\rho_n$. The class $wD_\M^2\1\{D_\M\le T_n\}$ has a constant envelope $U=CT_n^2/v$ and second moment at most $CT_n^2/v$, since
\begin{align*}
P[w^2D_\M^4\1\{D_\M\le T_n\}]\le\frac{CT_n^2}{v}P[wD_\M^2]\le\frac{C'T_n^2}{v}.
\end{align*}
Multiplication by $D_\M^2\1\{D_\M\le T_n\}$ preserves the kernel covering bound under the probability measure proportional to $D_\M^4\1\{D_\M\le T_n\}\,dQ$. A constant-envelope covering bound follows because the multiplier is at most $T_n^2$. The envelope-to-variance ratio is of order $T_n/\sqrt v=C_0\sqrt{n/\log n}$, so its logarithm is $O(\log n)$ for fixed $C_0$. By \eqref{eq:vc_expectation_bound},
\begin{align*}
\E\left[\sup_{\mx_A}|(P_n-P)(wD_\M^2\1\{D_\M\le T_n\})|\right]\le C(T_n\rho_n+T_n^2\rho_n^2)=O(C_0+C_0^2).
\end{align*}
The means are bounded, so Markov's inequality gives boundedness in probability for the truncated averages. Given a probability tolerance, Condition~\ref{con:B2} allows $C_0$ to be fixed so that $\P(\max_iD_i>T_n)$ is smaller than that tolerance for all sufficiently large $n$. On the complementary event the truncated and original averages agree. This proves both empirical assertions. Since $D_i\le D_i^2$, it also proves their first-moment versions.

By \eqref{eq:population_weight_envelope} and Lemma~\ref{lem:marginalization},
\begin{align*}
P[\widetilde W_1^2D_\M^2]&\le CJ_2\int\cL^2(1+\|\bm\Phi\|^2)f_A\dd\omega_{\md_A}\le C/v,\\
P[|\widetilde W_1|D_\M\1\{D_\M>T\}]&\le\frac{CT_G(T)}{T}\int\cL(1+\|\bm\Phi\|)f_A\dd\omega_{\md_A}\le CT_G(T)/T.
\end{align*}
For LC omit the tangent factor and use the positive lower bound for its denominator. Only the first power of $G$ is integrated in these calculations.
\end{proof}

\section{Uniform convergence for Wasserstein responses}\label{supp:wasserstein}
\setcounter{equation}{0}

Throughout this section, $\M=\mathcal W_2([a,b])$, $\mathbb H=L^2(0,1)$ and $B=\max\{|a|,|b|\}$. The quantile set $\mathcal Q_{[a,b]}$ consists of the elements of $\mathbb H$ admitting a nondecreasing representative with values in $[a,b]$. Quantiles are identified up to Lebesgue-null sets. In particular, $\|Q_y\|_{\mathbb H}\le B$, including for distributions with atoms.

\begin{lemma}\label{lem:quantile_projection}
The set $\mathcal Q_{[a,b]}$ is nonempty, closed and convex. If a real-valued weight $W$ is integrable and $\E[W]=1$, then $\E[WW_2^2(Y,y)]$ has a unique minimizer whose quantile is $\Pi_{\mathcal Q_{[a,b]}}\E[WQ_Y]$. The same conclusion holds for an empirical weighted objective with weights summing to one, whether or not the weights are nonnegative. These objectives satisfy a global quadratic curvature inequality with constant one. The empirical estimators, including the fallback, have jointly measurable versions on the entire predictor domain.
\end{lemma}

\begin{proof}
The constant function $a$ belongs to $\mathcal Q_{[a,b]}$, and convex combinations preserve monotonicity and the interval constraint. To prove closedness, let $q_j\to q$ in $L^2(0,1)$ and take a subsequence converging almost everywhere on a full-measure set $E$. For $u<v$ in $E$, monotone representatives satisfy $q_j(u)\le q_j(v)$, so $q(u)\le q(v)$. Extending this ordered function from $E$ by its one-sided limits gives a nondecreasing representative with values in $[a,b]$. Thus $q\in\mathcal Q_{[a,b]}$.

The Bochner integral $Z=\E[WQ_Y]$ exists because $\E[|W|]\,B<\infty$. For $q\in\mathcal Q_{[a,b]}$, expansion of the square gives
\begin{align*}
\E[W\|Q_Y-q\|_{\mathbb H}^2]&=\E[W\|Q_Y\|_{\mathbb H}^2]-2\langle Z,q\rangle_{\mathbb H}+\|q\|_{\mathbb H}^2\\
&=\|q-Z\|_{\mathbb H}^2+\E[W\|Q_Y\|_{\mathbb H}^2]-\|Z\|_{\mathbb H}^2.
\end{align*}
Projection onto a nonempty closed convex subset of a Hilbert space gives the unique minimizer $q^*=\Pi_{\mathcal Q_{[a,b]}}Z$. For $q\in\mathcal Q_{[a,b]}$, $\langle Z-q^*,q-q^*\rangle_{\mathbb H}\le0$, and hence
\begin{align*}
\|q-Z\|_{\mathbb H}^2-\|q^*-Z\|_{\mathbb H}^2=\|q-q^*\|_{\mathbb H}^2+2\langle q-q^*,q^*-Z\rangle_{\mathbb H}\ge\|q-q^*\|_{\mathbb H}^2.
\end{align*}
This proves curvature without assuming $W\ge0$.

For empirical weights take $\alpha_i=n^{-1}\widehat W_{\mx_A,\mathbf t_A,s}(\mX_A^{(i)})$. By Lemma~\ref{lem:population_design}, $\sum_i\alpha_i=1$ on the appropriate design event. Substituting $Z=\sum_i\alpha_iQ_{Y^{(i)}}$ in the preceding calculation proves \eqref{eq:wasserstein_projection}. If the design fails, the estimator is the empirical Fr\'echet mean, whose quantile is the ordinary average. Conditional mean quantiles belong to the closed convex set, so the same argument gives the conditional Fr\'echet mean and its curvature.

The map $y\mapsto Q_y$ is the isometric embedding of $\mathcal W_2([a,b])$ into the separable space $\mathbb H$. Kernel evaluation and the Moore--Penrose inverse on its rank strata are Borel. Projection onto $\mathcal Q_{[a,b]}$ is continuous. Their composition and the specified fallback therefore give a jointly measurable estimator for every data set and predictor value. This also supplies simultaneous existence of the minimizers in this response space.
\end{proof}

For the remainder of this section, write
\begin{align*}
\mu_A(\mx_A)&:=\E[Q_Y\mid\mX_A=\mx_A],\\
\widetilde Z_{\mathbf t_A,s}(\mx_A)&:=\E[\widetilde W_{\mx_A,\mathbf t_A,s}(\mX_A)Q_Y],\\
\widehat Z_{\mathbf t_A,s}(\mx_A)&:=\frac1n\sum_{i=1}^n\widehat W_{\mx_A,\mathbf t_A,s}(\mX_A^{(i)})Q_{Y^{(i)}}.
\end{align*}
The last expression is defined on the design domain and is set equal to the empirical mean quantile elsewhere. By Lemma~\ref{lem:quantile_projection}, $Q_{m_{\oplus,A}(\mx_A)}=\mu_A(\mx_A)$ and the local fitted quantiles are the corresponding projections of $\widetilde Z$ and $\widehat Z$.

\begin{lemma}\label{lem:quantile_bias}
Under Conditions~\ref{con:L1} and \ref{con:D1}--\ref{con:D2}, for $s=0,1$,
\begin{align*}
\sup_{\mx_A}\|\widetilde Z_{\mathbf t_A,s}(\mx_A)-\mu_A(\mx_A)\|_{\mathbb H}=O(b_{\mathbf t}(A)).
\end{align*}
The same bound holds for $\sup_{\mx_A}W_2(\widetilde m_{\mathbf t_A,s}(\mx_A),m_{\oplus,A}(\mx_A))$. The functions $\mu_A$ have twice continuously differentiable homogeneous extensions with uniformly bounded derivatives.
\end{lemma}

\begin{proof}
By Lemma~\ref{lem:marginalization},
\begin{align*}
\mu_A(\mx_A)=\int_\M Q_y g_{A,y}(\mx_A)P_Y(\dd y).
\end{align*}
For derivative orders zero, one and two, the Hilbert norm of the integrand is bounded by $BCG(y)$. This is integrable by Condition~\ref{con:D2}. Difference quotients and dominated convergence therefore justify differentiation of the Bochner integral and give bounded continuous derivatives of the homogeneous extension. The product $f_A\mu_A$ has the same property.

For LC, subtracting the target from the normalized kernel average gives
\begin{align*}
\widetilde Z_{\mathbf t_A,0}(\mx_A)-\mu_A(\mx_A)=\frac{\int\cL_{\mx_A,\mathbf t_A}(\mz_A)f_A(\mz_A)(\mu_A(\mz_A)-\mu_A(\mx_A))\dd\omega_{\md_A}(\mz_A)}{\widetilde\mu_{\mathbf t_A}(\mx_A)}.
\end{align*}
Apply the Hilbert-valued form of \eqref{eq:product_taylor_remainder} to $F(\mz_A)=f_A(\mz_A)(\mu_A(\mz_A)-\mu_A(\mx_A))$, holding $\mx_A$ fixed. The linear term integrates to zero by Lemma~\ref{lem:moment_orders}. The remainder has norm at most $C\sum_\ell\theta_\ell^2$, whose kernel integral is at most $Cb_{\mathbf t}(A)$. The denominator is bounded away from zero by Lemma~\ref{lem:population_design}.

For LL, define the first-order polynomial in the tangent displacement by
\begin{align*}
(P_1\mu_A)(\mz_A):=\mu_A(\mx_A)+\sum_{\ell\in A}t_\ell\,\mathrm D_\ell\bar\mu_A(\mx_A)[\bm\Phi_\ell(\mz_A)].
\end{align*}
Here $\mathrm D_\ell$ is the derivative as a linear map into $\mathbb H$, and $\bm\Phi_\ell$ is the $\ell$th block of the scaled tangent vector. By \eqref{eq:population_reproduction},
\begin{align*}
\widetilde Z_{\mathbf t_A,1}(\mx_A)-\mu_A(\mx_A)=\E[\widetilde W_1(\mu_A(\mX_A)-(P_1\mu_A)(\mX_A))].
\end{align*}
The Taylor remainder, \eqref{eq:population_weight_envelope} and Lemma~\ref{lem:moment_orders} yield
\begin{align*}
\|\widetilde Z_{\mathbf t_A,1}(\mx_A)-\mu_A(\mx_A)\|_{\mathbb H}&\le C\sum_{\ell\in A}t_\ell^2\E[\cL(1+\|\bm\Phi\|)\|\bm\Phi_\ell\|^2]\\
&\le C'\sum_{\ell\in A}t_\ell^2.
\end{align*}
All derivatives and moment constants are uniform in $\mx_A$. Since $\mu_A(\mx_A)\in\mathcal Q_{[a,b]}$, nonexpansiveness of projection proves the asserted Wasserstein bias bound. The empty-model bias is zero.
\end{proof}

\begin{lemma}\label{lem:quantile_stochastic}
Under Conditions~\ref{con:L1}, \ref{con:D1} and the fixed-model version of Condition~\ref{con:B1}, for $s=0,1$,
\begin{align*}
\sup_{\mx_A}\|\widehat Z_{\mathbf t_A,s}(\mx_A)-\widetilde Z_{\mathbf t_A,s}(\mx_A)\|_{\mathbb H}=O_\P\left(\sqrt{\frac{\log n}{nv_{\mathbf t}(A)}}\right).
\end{align*}
\end{lemma}

\begin{proof}
Put $v=v_{\mathbf t}(A)$ and $\rho_n=(\log n/(nv))^{1/2}$. First consider
\begin{align*}
S_n(\mx_A,u):=(P_n-P)(\widetilde W_{\mx_A,\mathbf t_A,s}Q_Y(u)),\quad 0<u<1.
\end{align*}
For fixed $u$, $Q_Y(u)$ is a measurable multiplier bounded by $B$. By Lemma~\ref{lem:population_design}, $v\widetilde W_s$ is a finite linear combination of $v\cL$ and its first tangent multipliers with uniformly bounded coefficients. Enlarge these coefficients to a fixed Euclidean ball. Lemma~\ref{lem:kernel_covering} then gives common covering constants, independent of $u$ and of the bandwidths. Lemma~\ref{lem:moment_orders} gives a constant envelope $U=C_0/v$ and a second-moment proxy $\sigma^2=C_1/v$ for the scalar class.

Let $\lambda_n=\log(\mathcal A U/\sigma)$. Then $\lambda_n=O(\log n)$ and
\begin{align*}
\frac{U^2\lambda_n}{n\sigma^2}=O\left(\frac{\log n}{nv}\right)=o(1).
\end{align*}
The probability inequality (2.3), with its lower-range condition (2.4), in \citet[Theorem~2.1]{GineGuillou2002S} therefore implies
\begin{align*}
\P\left(n\sup_{\mx_A}|S_n(\mx_A,u)|>w\right)\le C\exp\left(-\frac{c w^2}{n/v+w/v}\right)
\end{align*}
for $w\ge C(\lambda_n/v+\sqrt{n\lambda_n/v})$. To see the denominator, the square $(\sqrt n\sigma+U\sqrt{\lambda_n})^2$ in that inequality is at most $Cn/v$, and $\log(1+x)\ge x/(1+x)$ gives the displayed form. Thus, after enlarging the constants, for $z\ge0$,
\begin{align}\label{eq:quantile_process_tail}
\P\left(\sup_{\mx_A}|S_n(\mx_A,u)|>C\left(\sqrt{\frac{\log n+z}{nv}}+\frac{\log n+z}{nv}\right)\right)\le C'e^{-z}.
\end{align}
The constants are independent of $u$.

Let $a(z)=C(\sqrt{(\log n+z)/(nv)}+(\log n+z)/(nv))$. Integrating \eqref{eq:quantile_process_tail} gives
\begin{align*}
\sup_{0<u<1}\E\left[\sup_{\mx_A}|S_n(\mx_A,u)|^2\right]&\le a(0)^2+C'\int_0^\infty e^{-z}\frac{\dd}{\dd z}a(z)^2\dd z\\
&\le C''(\rho_n^2+\rho_n^4).
\end{align*}
Quantile evaluation is jointly measurable in $(Y,u)$. The measurable projection theorem, applied over the compact predictor index, makes the scalar supremum measurable on the completion of the product of the sampling space and $(0,1)$ with Lebesgue measure. The same argument applies to its sections; completion does not change their integrals. Thus Tonelli's theorem gives
\begin{align*}
\E\left[\sup_{\mx_A}\|S_n(\mx_A,\cdot)\|_{\mathbb H}^2\right]&\le\int_0^1\E\left[\sup_{\mx_A}|S_n(\mx_A,u)|^2\right]\dd u\\
&\le C(\rho_n^2+\rho_n^4)=O(\rho_n^2).
\end{align*}
Markov's inequality proves the required order for the population-weighted process. No covering number of a response-indexed quantile class is used.

On the simultaneous design event, the estimated-weight term satisfies
\begin{align*}
\sup_{\mx_A}\|P_n((\widehat W_s-\widetilde W_s)Q_Y)\|_{\mathbb H}\le B\sup_{\mx_A}P_n|\widehat W_s-\widetilde W_s|=O_\P(\rho_n)
\end{align*}
by Lemma~\ref{lem:empirical_design}. The exact decomposition
\begin{align*}
\widehat Z-\widetilde Z=(P_n-P)(\widetilde W_sQ_Y)+P_n((\widehat W_s-\widetilde W_s)Q_Y)
\end{align*}
completes the proof on that event. Its complement has probability tending to zero, so the fallback does not affect the asserted stochastic order.
\end{proof}

\begin{proof}[Proof of \Cref{thm:wasserstein}]
Fix a nonempty $A$, $s\in\{0,1\}$ and $\mathbf t\in\{\mb,\mh\}$. By \eqref{eq:subset_bandwidth_orders}, the fixed-model bandwidth conditions hold. Lemma~\ref{lem:quantile_projection} and nonexpansiveness give, on the design event,
\begin{align*}
\sup_{\mx_A}W_2(\widehat m_{\mathbf t_A,s}(\mx_A),m_{\oplus,A}(\mx_A))&\le\sup_{\mx_A}\|\widehat Z_{\mathbf t_A,s}(\mx_A)-\widetilde Z_{\mathbf t_A,s}(\mx_A)\|_{\mathbb H}\\
&\quad+\sup_{\mx_A}\|\widetilde Z_{\mathbf t_A,s}(\mx_A)-\mu_A(\mx_A)\|_{\mathbb H}\\
&=O_\P\left(b_{\mathbf t}(A)+\sqrt{\frac{\log n}{nv_{\mathbf t}(A)}}\right)
\end{align*}
by Lemmas~\ref{lem:quantile_bias}--\ref{lem:quantile_stochastic}. The design-failure probability tends to zero. For $A=\varnothing$, the empirical mean quantile satisfies
\begin{align*}
\E\left[\left\|\frac1n\sum_{i=1}^n(Q_{Y^{(i)}}-\E[Q_Y])\right\|_{\mathbb H}^2\right]=\frac1n\E[\|Q_Y-\E[Q_Y]\|_{\mathbb H}^2]\le\frac{4B^2}{n},
\end{align*}
which is stronger than the required empty-model bound. There are finitely many factor sets and two bandwidth vectors, so a finite union gives simultaneous probability constants for their individual rates.

For later selection arguments, Lemma~\ref{lem:quantile_projection} gives existence, uniqueness and quadratic separation of the conditional targets. Lemma~\ref{lem:quantile_bias} gives their continuity, and all response distances are at most $b-a$. Thus the target and second-moment properties needed for selection hold in this branch without Conditions~\ref{con:M5}--\ref{con:M6} or the tail assumption~\ref{con:B2}.
\end{proof}

\section{Uniform convergence for responses in general metric spaces}\label{supp:stochastic}
\setcounter{equation}{0}

\subsection{Conditional targets and preliminary consistency}

In this section write $m_0=m_{\oplus,\varnothing}$. For $y,m,w\in\M$, the triangle inequality gives
\begin{align}\label{eq:metric_increment}
|d_\M^2(w,y)-d_\M^2(w,m)|\le d_\M(y,m)(2d_\M(w,m)+d_\M(y,m)).
\end{align}
If $d_\M(m,m_0)\le R$ and $d_\M(y,m)\le\delta\le R'$, the right-hand side is at most $C_{R,R'}\delta D_\M(w)$. This bound uses $D_\M(w)=\max\{1,d_\M(w,m_0)\}$.

\begin{lemma}\label{lem:target_continuity}
Under Conditions~\ref{con:D1}--\ref{con:D2} and \ref{con:M1}--\ref{con:M2}, the conditional Fr\'echet means belong to a common bounded ball centered at $m_0$. For every $A\subseteq[p]$, $m_{\oplus,A}$ is continuous and its image $\cK_A$ is compact. For every fixed $R<\infty$,
\begin{align*}
\sup_{y:\,d_\M(y,m_0)\le R}|M_{\oplus,A}(\mx_A',y)-M_{\oplus,A}(\mx_A,y)|\le C_R\|\mx_A'-\mx_A\|.
\end{align*}
The constants can be chosen uniformly over $A$.
\end{lemma}

\begin{proof}
Write $m=m_{\oplus,A}(\mx_A)$. Conditional optimality gives
\begin{align*}
\E[d_\M^2(Y,m)\mid\mX_A=\mx_A]&\le\E[d_\M^2(Y,m_0)\mid\mX_A=\mx_A],\\
d_\M(m,m_0)^2&\le2\E[d_\M^2(Y,m)+d_\M^2(Y,m_0)\mid\mX_A=\mx_A]\le4CJ_2,
\end{align*}
where Lemma~\ref{lem:marginalization} gives the last bound.

For $d_\M(y,m_0)\le R$, $d_\M^2(w,y)\le(1+R)^2D_\M(w)^2$. The conditional density representation and Lemma~\ref{lem:marginalization} consequently give uniformly bounded first and second predictor derivatives of $M_{\oplus,A}(\cdot,y)$ on the polysphere. To bound the first difference, integrate the first derivative along a product geodesic between $\mx_A$ and $\mx_A'$. Since $\arccos(\mx_\ell\tran\mx_\ell')\le(\pi/2)\|\mx_\ell-\mx_\ell'\|$, this yields the stated bound with a constant depending only on $R$ and the fixed regularity quantities.

Let $\mx_{A,j}\to\mx_A$, $m_j=m_{\oplus,A}(\mx_{A,j})$ and $m=m_{\oplus,A}(\mx_A)$. All these targets belong to one fixed bounded ball. By their respective optimality properties,
\begin{align*}
0&\le M_{\oplus,A}(\mx_A,m_j)-M_{\oplus,A}(\mx_A,m)\\
&\le2\sup_{y:\,d_\M(y,m_0)\le R}|M_{\oplus,A}(\mx_{A,j},y)-M_{\oplus,A}(\mx_A,y)|\longrightarrow0.
\end{align*}
Condition~\ref{con:M2} forces $d_\M(m_j,m)\to0$. This proves continuity without compactness of bounded subsets of $\M$. The image $\cK_A$ is compact because the predictor space is compact. The empty-model target is constant, so its assertions follow immediately.
\end{proof}

\begin{lemma}\label{lem:population_prelocalization}
Fix $s\in\{0,1\}$ and suppose Conditions~\ref{con:L1}, \ref{con:D1}--\ref{con:D2} and \ref{con:M1}--\ref{con:M3} hold for this $s$. For sufficiently small bandwidths, the population local minimizers belong to a common bounded ball centered at $m_0$. For every fixed $R<\infty$,
\begin{align*}
\sup_{\mx_A}\sup_{y:\,d_\M(y,m_0)\le R}|\widetilde M_{\mathbf t_A,s}(\mx_A,y)-M_{\oplus,A}(\mx_A,y)|=O(b_{\mathbf t}(A)).
\end{align*}
Consequently $\sup_{\mx_A}d_\M(\widetilde m_{\mathbf t_A,s}(\mx_A),m_{\oplus,A}(\mx_A))\to0$ as $\max_\ell t_\ell\to0$.
\end{lemma}

\begin{proof}
For any integrable signed weight with $\int W=1$, write $a(w)=d_\M(w,m_0)$ and $r=d_\M(y,m_0)$. The triangle inequality implies
\begin{align*}
|d_\M^2(w,y)-r^2|\le2ra(w)+a(w)^2.
\end{align*}
Putting $K_j=\int|W|a^j$ gives
\begin{align*}
\int Wd_\M^2(w,y)\ge r^2-2K_1r-K_2,\quad \int Wd_\M^2(w,m_0)\le K_2.
\end{align*}
At any existing minimizer, comparison with $m_0$ gives
\begin{align}\label{eq:minimizer_radius_bound}
r^2-2K_1r-2K_2\le0,\quad r\le K_1+\sqrt{K_1^2+2K_2}.
\end{align}
Lemmas~\ref{lem:marginalization}, \ref{lem:moment_orders} and \ref{lem:population_design} bound the population values of $K_1,K_2$ uniformly. Condition~\ref{con:M3}(L) supplies attainment. The preceding calculation bounds existing minimizers and does not replace that existence condition.

Fix $R<\infty$ and let $F_y(\mz_A)=M_{\oplus,A}(\mz_A,y)$. For derivative orders zero, one and two, the density-ratio integral is dominated by $C_R D_\M(w)^2G(w)$. Thus the homogeneous extensions of $F_y$ and $f_AF_y$ have derivatives bounded uniformly over $d_\M(y,m_0)\le R$. For LC,
\begin{align*}
\widetilde M_{\mathbf t_A,0}(\mx_A,y)-F_y(\mx_A)=\frac{\int\cL_{\mx_A,\mathbf t_A}(\mz_A)(F_y(\mz_A)-F_y(\mx_A))f_A(\mz_A)\dd\omega_{\md_A}(\mz_A)}{\widetilde\mu_{\mathbf t_A}(\mx_A)}.
\end{align*}
By \eqref{eq:product_taylor_remainder}, the numerator is the sum of a linear tangent term and a remainder bounded by $C_R\int\cL\sum_\ell\theta_\ell^2\dd\omega_{\md_A}$. The linear integral is zero and the remainder is $O(b_{\mathbf t}(A))$ by Lemma~\ref{lem:moment_orders}. The denominator is bounded away from zero.

For LL, let $(P_1F_y)(\mz_A)=F_y(\mx_A)+\sum_\ell t_\ell\nabla_\ell\bar F_y(\mx_A)\tran\bm\Phi_\ell(\mz_A)$. By \eqref{eq:population_reproduction},
\begin{align*}
\widetilde M_{\mathbf t_A,1}(\mx_A,y)-F_y(\mx_A)=\E[\widetilde W_1(F_y(\mX_A)-(P_1F_y)(\mX_A))].
\end{align*}
By \eqref{eq:product_taylor_remainder} and \eqref{eq:population_weight_envelope}, the absolute value is bounded by $C_R\sum_\ell t_\ell^2\E[\cL(1+\|\bm\Phi\|)\|\bm\Phi_\ell\|^2]$, hence by $C_R'b_{\mathbf t}(A)$. This proves the bounded-region approximation for both estimators.

The conditional and population local targets lie in one such region. Let $m=m_{\oplus,A}(\mx_A)$ and $\widetilde m=\widetilde m_{\mathbf t_A,s}(\mx_A)$. Their optimality gives
\begin{align*}
0\le M_{\oplus,A}(\mx_A,\widetilde m)-M_{\oplus,A}(\mx_A,m)\le2\sup_{y:\,d_\M(y,m_0)\le R}|\widetilde M_{\mathbf t_A,s}(\mx_A,y)-M_{\oplus,A}(\mx_A,y)|.
\end{align*}
The right-hand side tends to zero uniformly. Condition~\ref{con:M2} proves the asserted population consistency.
\end{proof}

\begin{proposition}\label{prop:empirical_prelocalization}
Fix $s\in\{0,1\}$ and assume Conditions~\ref{con:L1}, \ref{con:D1}--\ref{con:D2}, \ref{con:M1}--\ref{con:M3}, \ref{con:M6} and the fixed-model versions of Conditions~\ref{con:B1}--\ref{con:B2}. Then
\begin{align*}
\sup_{\mx_A}d_\M(\widehat m_{\mathbf t_A,s}(\mx_A),m_{\oplus,A}(\mx_A))=o_\P(1).
\end{align*}
The same conclusion holds for the distance from $\widetilde m_{\mathbf t_A,s}$.
\end{proposition}

\begin{proof}
Fix a nonempty $A$ and write $v=v_{\mathbf t}(A)$, $b=b_{\mathbf t}(A)$ and $\rho_n=(\log n/(nv))^{1/2}$. Condition~\ref{con:M3}(E) supplies empirical minimizers at every valid predictor value almost surely. On the design event, apply \eqref{eq:minimizer_radius_bound} to the empirical measure. Lemmas~\ref{lem:empirical_design} and \ref{lem:local_weighted_moments} give
\begin{align*}
\sup_{\mx_A}P_n(|\widehat W_s|d_\M(Y,m_0)^j)=O_\P(1),\quad j=1,2.
\end{align*}
It follows that the supremum of the empirical minimizer radii is $O_\P(1)$.

For fixed $R$, put $B_R=\{y:d_\M(y,m_0)\le R\}$ and center each objective at $m_0$. For a fixed $y\in B_R$, let $\Delta_y(Y)=d_\M^2(Y,y)-d_\M^2(Y,m_0)$. By \eqref{eq:metric_increment}, $|\Delta_y(Y)|\le C_RD_\M(Y)$. The empirical-to-population-local contrast is
\begin{align*}
(P_n-P)(\widetilde W_s\Delta_y)+P_n((\widehat W_s-\widetilde W_s)\Delta_y).
\end{align*}
Truncate the first term at $D_\M\le T_n=C_0/\rho_n$. Its envelope is $C_RT_n/v$ and its second moment is at most $C_R/v$ by Lemma~\ref{lem:local_weighted_moments}. Kernel covering and \eqref{eq:vc_expectation_bound} give order $O_\P(\rho_n)$ for fixed $C_0$. The empirical truncation changes nothing on $\max_iD_\M(Y^{(i)})\le T_n$. The change in population centering is at most
\begin{align*}
C_RP[|\widetilde W_s|D_\M\1\{D_\M>T_n\}]\le C_RT_G(T_n)/T_n=o(\rho_n).
\end{align*}
Condition~\ref{con:B2} allows $C_0$ to be chosen for any prescribed probability tolerance. The second term is bounded by
\begin{align*}
C_nC_R\rho_nP_n[\cL(1+\|\bm\Phi\|)D_\M]=O_\P(\rho_n)
\end{align*}
by Lemmas~\ref{lem:empirical_design} and \ref{lem:local_weighted_moments}. The same argument applies at every point of any fixed finite response net.

Condition~\ref{con:M6} supplies a finite $\epsilon$-net $y_1,\ldots,y_{N_\epsilon}$ of $B_R$. By \eqref{eq:metric_increment}, whenever $d_\M(y,y_j)\le\epsilon$,
\begin{align*}
|\widehat M(\mx_A,y)-\widehat M(\mx_A,y_j)|\le C_R\epsilon P_n(|\widehat W_s|D_\M)=\epsilon O_\P(1)
\end{align*}
uniformly in $\mx_A$. The corresponding conditional difference is at most $C_R\epsilon J_1$. Lemma~\ref{lem:population_prelocalization} bounds the population-local approximation. Thus the supremum over $\mx_A$ and $y\in B_R$ of the difference of centered empirical and conditional objectives tends to zero in probability, by first sending $n\to\infty$ and then $\epsilon\downarrow0$ at the fixed $R$.

Fix $\eta>0$ and $\zeta>0$. By Lemma~\ref{lem:target_continuity} and the empirical radius bound, choose $R$ containing all conditional targets and such that the limiting upper probability of any empirical minimizer lying outside $B_R$ is below $\zeta$. Condition~\ref{con:M2} gives $c_\eta>0$ for the objective gap at distance greater than $\eta$. If the centered-objective approximation error on $B_R$ is below $c_\eta/3$, empirical optimality gives a conditional gap at the empirical minimizer below $2c_\eta/3$. Hence
\begin{align*}
\limsup_{n\to\infty}\P\left(\sup_{\mx_A}d_\M(\widehat m_{\mathbf t_A,s}(\mx_A),m_{\oplus,A}(\mx_A))>\eta\right)\le\zeta.
\end{align*}
The design event has probability tending to one. Since $\zeta$ is arbitrary, this proves the first conclusion. Lemma~\ref{lem:population_prelocalization} and the triangle inequality prove the second. For $A=\varnothing$, use weights one, the law of large numbers for $D_\M^2$, the same finite-net argument and Condition~\ref{con:M1}(E). There is only one evaluation point, so no additional simultaneous assertion is needed.
\end{proof}

\begin{remark}\label{rmk:localization_alternative}
Condition~\ref{con:M6} is used for the finite response net in Proposition~\ref{prop:empirical_prelocalization}. If bounded balls are not totally bounded, the subsequent rate argument remains valid with that preliminary consistency step replaced by
\begin{align*}
\sup_{\mx_A}d_\M(\widehat m_{\mathbf t_A,s}(\mx_A),\widetilde m_{\mathbf t_A,s}(\mx_A))\le\eta_\star
\end{align*}
with probability tending to one, for the fixed radius $\eta_\star$ used below, and with simultaneous empirical optimality on that event. This is an alternative assumption, not a consequence of local curvature or entropy alone.
\end{remark}

\subsection{Smoothing bias}\label{supp:bias}

\begin{lemma}\label{lem:localized_bias}
Fix $s\in\{0,1\}$ and assume Conditions~\ref{con:L1}, \ref{con:D1}--\ref{con:D2} and \ref{con:M1}--\ref{con:M4} for this $s$. For $m=m_{\oplus,A}(\mx_A)$ and $d_\M(y,m)\le\eta$ with any fixed sufficiently small $\eta>0$,
\begin{align*}
&|\widetilde M_{\mathbf t_A,s}(\mx_A,y)-\widetilde M_{\mathbf t_A,s}(\mx_A,m)-M_{\oplus,A}(\mx_A,y)+M_{\oplus,A}(\mx_A,m)|\\
&\qquad\le Cb_{\mathbf t}(A)d_\M(y,m).
\end{align*}
Consequently,
\begin{align*}
\sup_{\mx_A}d_\M(\widetilde m_{\mathbf t_A,s}(\mx_A),m_{\oplus,A}(\mx_A))=O\left(b_{\mathbf t}(A)^{1/(\beta_\oplus-1)}\right).
\end{align*}
\end{lemma}

\begin{proof}
Fix $\mx_A$, let $m=m_{\oplus,A}(\mx_A)$, and put $\delta=d_\M(y,m)$. When $\delta=0$ the increment is zero. For $\delta>0$, define within this proof
\begin{align*}
\Delta(\mz_A):=M_{\oplus,A}(\mz_A,y)-M_{\oplus,A}(\mz_A,m).
\end{align*}
Keep $m$ and $y$ fixed while differentiating in $\mz_A$. By \eqref{eq:metric_increment} and Lemma~\ref{lem:target_continuity},
\begin{align*}
|d_\M^2(w,y)-d_\M^2(w,m)|\le C\delta D_\M(w).
\end{align*}
Lemma~\ref{lem:marginalization} therefore gives the integrable derivative majorant $C\delta D_\M(w)G(w)$. For $j=0,1,2$,
\begin{align*}
\partial^j\bar\Delta(\mz_A)=\int_\M(d_\M^2(w,y)-d_\M^2(w,m))\partial^j\bar g_{A,w}(\mz_A)P_Y(\dd w),\quad \sup_{\mz_A}\|\partial^j\bar\Delta(\mz_A)\|\le C\delta J_1.
\end{align*}
The norm is absolute value for $j=0$, Euclidean norm for $j=1$ and operator norm for $j=2$. The product rule gives the same order for the derivatives of $f_A\Delta$. No derivative of the Fr\'echet regression function is used.

For LC, the increment of the smoothed objective satisfies
\begin{align*}
&\widetilde M_{\mathbf t_A,0}(\mx_A,y)-\widetilde M_{\mathbf t_A,0}(\mx_A,m)-\Delta(\mx_A)\\
&\qquad=\frac{\int\cL_{\mx_A,\mathbf t_A}(\mz_A)(\Delta(\mz_A)-\Delta(\mx_A))f_A(\mz_A)\dd\omega_{\md_A}(\mz_A)}{\widetilde\mu_{\mathbf t_A}(\mx_A)}.
\end{align*}
Apply \eqref{eq:product_taylor_remainder} to $(\Delta(\mz_A)-\Delta(\mx_A))f_A(\mz_A)$. Its linear term integrates to zero by Lemma~\ref{lem:moment_orders}, and its remainder has absolute value at most $C\delta\sum_\ell\theta_\ell^2$. The kernel integral is at most $C\delta b_{\mathbf t}(A)$. Dividing by the denominator bound in Lemma~\ref{lem:population_design} proves the LC assertion.

For LL, put
\begin{align*}
(P_1\Delta)(\mz_A):=\Delta(\mx_A)+\sum_{\ell\in A}t_\ell\nabla_\ell\bar\Delta(\mx_A)\tran\bm\Phi_\ell(\mz_A).
\end{align*}
This polynomial includes the constant term. By \eqref{eq:population_reproduction}, its weighted expectation is $\Delta(\mx_A)$, and hence
\begin{align*}
&\widetilde M_{\mathbf t_A,1}(\mx_A,y)-\widetilde M_{\mathbf t_A,1}(\mx_A,m)-\Delta(\mx_A)\\
&\qquad=\E[\widetilde W_1(\Delta(\mX_A)-(P_1\Delta)(\mX_A))].
\end{align*}
By \eqref{eq:product_taylor_remainder}, the remainder is at most $C\delta\sum_\ell t_\ell^2\|\bm\Phi_\ell\|^2$ away from a probability-null set of antipodes. Thus, by \eqref{eq:population_weight_envelope} and Lemma~\ref{lem:moment_orders}, the absolute weighted expectation is at most
\begin{align*}
C\delta\sum_{\ell\in A}t_\ell^2\E[\cL(1+\|\bm\Phi\|)\|\bm\Phi_\ell\|^2]\le C'\delta b_{\mathbf t}(A).
\end{align*}
This proves the first assertion for both estimators without a bandwidth-ratio restriction.

By Lemma~\ref{lem:population_prelocalization}, the population local target lies in the neighborhood in Condition~\ref{con:M4}(P) for all sufficiently small bandwidths. Let $d=d_\M(\widetilde m_{\mathbf t_A,s}(\mx_A),m)$. Population optimality and \eqref{eq:margin} give
\begin{align*}
c_\oplus d^{\beta_\oplus}&\le M_{\oplus,A}(\mx_A,\widetilde m_{\mathbf t_A,s}(\mx_A))-M_{\oplus,A}(\mx_A,m)\\
&\le Cb_{\mathbf t}(A)d.
\end{align*}
Divide by $d$ when $d>0$. The constants are uniform in $\mx_A$, which proves the bias rate. For $A=\varnothing$ the bias is zero.
\end{proof}

\subsection{Localized stochastic bounds}

For a small fixed $\eta>0$, write $\cN_A:=\{z\in\M:d_\M(z,\cK_A)<\eta\}$. The constants and this neighborhood in the following lemma do not depend on the bandwidth.

\begin{lemma}\label{lem:moving_anchor_entropy}
Under the conditions of Lemma~\ref{lem:target_continuity} and Condition~\ref{con:M5}, $\eta$ can be chosen so that $\cN_A$ is totally bounded and
\begin{align*}
N(\epsilon,\cN_A,d_\M)\le C\epsilon^{-q},\quad 0<\epsilon<\eta,
\end{align*}
for fixed $C,q<\infty$. For $0<\delta\le\eta/4$, let
\begin{align*}
\cH_\delta:=\{w\mapsto d_\M^2(w,y)-d_\M^2(w,a):a,y\in\cN_A,\ d_\M(y,a)\le\delta\}.
\end{align*}
For every probability measure $Q$ with $\|D_\M\|_{Q,2}<\infty$,
\begin{align*}
N(\epsilon\delta\|D_\M\|_{Q,2},\cH_\delta,L_2(Q))\le C\left(\frac{C}{\delta\epsilon}\right)^{2q},\quad 0<\epsilon\le1.
\end{align*}
If the assumptions of Lemma~\ref{lem:population_prelocalization} also hold for the chosen $s$, then for sufficiently small bandwidths the class includes all increments centered at $a=\widetilde m_{\mathbf t_A,s}(\mx_A)$ with $d_\M(y,a)\le\delta$, uniformly in $\mx_A$.
\end{lemma}

\begin{proof}
Choose $\eta>0$ so that $4\eta<\min(r_\M,\delta_\M)$. Monotonicity of covering numbers and Condition~\ref{con:M5}, integrated over $[1/16,1/8]$, give
\begin{align*}
\frac1{16}\sqrt{1+\log N(r/8,B_\M(z,r),d_\M)}\le C_\M
\end{align*}
whenever $r\le\delta_\M$ and $d_\M(z,\cK_A)<r_\M$. Thus every such ball is covered by at most $N_0=\lceil\exp(256C_\M^2)\rceil$ balls of radius $r/8$.

By Lemma~\ref{lem:target_continuity}, a finite $\eta$-net of $\cK_A$ covers $\cN_A$ by balls of radius $2\eta$. Cover each such ball at radius $r/8$, discard covering balls not meeting $\cN_A$, and choose a new center in each nonempty intersection. The retained part is contained in the ball of radius $r/2$ about its new center, since the recentering distance is at most $r/8$ and the enlarged radius allows strict containment. Every new center belongs to $\cN_A$ and remains within distance $\eta<r_\M$ of $\cK_A$. This argument can be iterated at every scale. After $j$ steps there are at most $N_1N_0^j$ balls of radius $2\eta\,2^{-j}$. Choosing $j$ according to $\epsilon$ proves the polynomial bound with $q=\log_2N_0$ and a different constant $C$. The same bound, up to a constant and a factor two in radius, holds for the closure in the metric completion.

All response arguments in $\cN_A$ lie in a fixed bounded ball about $m_0$. Twice applying \eqref{eq:metric_increment} yields
\begin{align*}
&|[d_\M^2(w,y)-d_\M^2(w,a)]-[d_\M^2(w,y')-d_\M^2(w,a')]|\\
&\qquad\le CD_\M(w)(d_\M(y,y')+d_\M(a,a')).
\end{align*}
Cover each response argument by a net of radius $\delta\epsilon/(4C)$. There are at most $C(C/(\delta\epsilon))^{2q}$ pairs of net cells. For every occupied cell pair choose an actual admissible pair $(a,y)$ as its representative. The preceding inequality shows that the corresponding functions form the asserted $L_2(Q)$ cover. This choice of representatives keeps the centers of the function cover in the class even if the two original net centers are not at distance at most $\delta$.

By Lemma~\ref{lem:population_prelocalization}, the population local targets have distance at most $\eta/4$ from $\cK_A$ for small bandwidths. Every $y$ at distance at most $\delta\le\eta/4$ from such a target belongs to $\cN_A$, proving the final assertion.
\end{proof}

Fix $0<\eta_\star<\min(1,\eta_\oplus,\eta/4)/4$. For a fixed $A,\mathbf t_A,s$, put $a(\mx_A)=\widetilde m_{\mathbf t_A,s}(\mx_A)$ and define
\begin{align*}
\mathbb Z_n(\mx_A,y):=\widehat M_{\mathbf t_A,s}(\mx_A,y)-\widehat M_{\mathbf t_A,s}(\mx_A,a(\mx_A))-\widetilde M_{\mathbf t_A,s}(\mx_A,y)+\widetilde M_{\mathbf t_A,s}(\mx_A,a(\mx_A)).
\end{align*}
The empirical objectives in this expression are used on their simultaneous design event.

\begin{lemma}\label{lem:complete_modulus}
Assume Conditions~\ref{con:L1}, \ref{con:D1}--\ref{con:D2}, \ref{con:M1}--\ref{con:M5} for the fixed $s$ and the fixed-model versions of Conditions~\ref{con:B1}--\ref{con:B2}. Let $\rho_n=(\log n/(nv_{\mathbf t}(A)))^{1/2}$. For deterministic $0<\delta_n\le\eta_\star$ with $\log(e/\delta_n)=O(\log n)$,
\begin{align*}
\sup_{\mx_A}\sup_{d_\M(y,a(\mx_A))\le\delta_n}|\mathbb Z_n(\mx_A,y)|=O_\P(\delta_n\rho_n).
\end{align*}
Put $r_n=\rho_n^{1/(\beta_\oplus-1)}$ and $\delta_{j,n}=\min(2^jr_n,\eta_\star)$ for $0\le j\le J_n$, where $J_n$ is the first integer for which $2^{J_n}r_n\ge\eta_\star$. Then, on the simultaneous design event,
\begin{align*}
\max_{0\le j\le J_n}\sup_{\mx_A}\frac1{\delta_{j,n}}\sup_{d_\M(y,a(\mx_A))\le\delta_{j,n}}|\mathbb Z_n(\mx_A,y)|=O_\P(\rho_n).
\end{align*}
\end{lemma}

\begin{proof}
Abbreviate $v=v_{\mathbf t}(A)$ and $\Delta_{y,a}(Y)=d_\M^2(Y,y)-d_\M^2(Y,a(\mx_A))$. Direct subtraction of the two objectives gives
\begin{align}\label{eq:complete_increment_decomposition}
\mathbb Z_n(\mx_A,y)=(P_n-P)(\widetilde W_s\Delta_{y,a})+P_n((\widehat W_s-\widetilde W_s)\Delta_{y,a}).
\end{align}
For $d_\M(y,a(\mx_A))\le\delta_n$, Lemma~\ref{lem:population_prelocalization} and \eqref{eq:metric_increment} give $|\Delta_{y,a}(Y)|\le C\delta_nD_\M(Y)$.

Consider the first term. The normalized weight $v\widetilde W_s$ belongs to a finite linear combination of the kernel classes in Lemma~\ref{lem:kernel_covering}, with coefficients in a fixed Euclidean ball. Lemma~\ref{lem:moving_anchor_entropy} bounds the loss-increment class, including both moving response arguments. Product covering, using the measure proportional to $D_\M^2\,dQ$ for the kernel part, yields a fixed exponent and a prefactor whose logarithm is at most $C\log(C/\delta_n)$ for the normalized product class.

Fix $C_0>0$ and truncate at $D_\M\le T_n=C_0/\rho_n$. By Lemma~\ref{lem:local_weighted_moments}, the truncated class has envelope and second-moment proxies
\begin{align*}
U_n=C\delta_nT_n/v,\quad \sigma_n^2=C'\delta_n^2/v.
\end{align*}
Its covering logarithm at these proxies satisfies
\begin{align*}
\log(\mathcal A_nU_n/\sigma_n)\le C\left(\log(C/\delta_n)+\log(C+T_n/\sqrt v)\right)=O(\log n),
\end{align*}
because $T_n/\sqrt v=C_0\sqrt{n/\log n}$. By \eqref{eq:vc_expectation_bound}, with a fixed covering exponent,
\begin{align*}
\E\left[\sup_{\mx_A,y}|(P_n-P)(\widetilde W_s\Delta_{y,a}\1\{D_\M\le T_n\})|\right]&\le C\left(\delta_n\sqrt{\frac{\log n}{nv}}+\frac{\delta_nT_n\log n}{nv}\right)\\
&\le C'(1+C_0)\delta_n\rho_n.
\end{align*}
The dependence of $\mathcal A_n$ on $\delta_n$ is included in the displayed logarithm, rather than in an unspecified maximal-inequality constant.

On $\max_iD_\M(Y^{(i)})\le T_n$ the empirical sums before and after truncation agree. The population centerings differ by at most
\begin{align*}
C\delta_nP[|\widetilde W_s|D_\M\1\{D_\M>T_n\}]\le C\delta_nT_G(T_n)/T_n=o(\delta_n\rho_n).
\end{align*}
Condition~\ref{con:B2} makes the complementary empirical event arbitrarily unlikely after choosing $C_0$. Markov's inequality proves the $O_\P(\delta_n\rho_n)$ bound for the first term in \eqref{eq:complete_increment_decomposition}.

For the second term, Lemmas~\ref{lem:empirical_design} and \ref{lem:local_weighted_moments} give
\begin{align*}
\sup_{\mx_A,y}P_n(|\widehat W_s-\widetilde W_s||\Delta_{y,a}|)&\le C_n\delta_n\rho_n\sup_{\mx_A}P_n[\cL(1+\|\bm\Phi\|)D_\M]\\
&=O_\P(\delta_n\rho_n).
\end{align*}
Its random multiplicative bound is independent of the radius. Combining the two terms proves the first conclusion.

For the second conclusion, divide the increment class at radius $\delta_{j,n}$ by $\delta_{j,n}$ and take the union over $j$. Since $v\le1$ eventually and $nv/\log n\to\infty$,
\begin{align*}
\rho_n\to0,\quad \log(1/r_n)=O(\log n),\quad J_n=O(\log n).
\end{align*}
The union adds $\log(J_n+1)=O(\log\log n)$ to the covering logarithm. Its radius-normalized envelope is at most $CT_n/v$, its second moment is at most $C/v$, and its centering tail is at most $CT_G(T_n)/T_n=o(\rho_n)$. The smallest radius has $\log(e/r_n)=O(\log n)$, so the preceding maximal-inequality calculation applies to this union with the same order $O_\P(\rho_n)$. The estimated-weight bound is already linear in each radius. This supplies one random bound for every predictor and every radius, proving the second conclusion.
\end{proof}

\subsection{Proof of the uniform convergence theorem}\label{supp:uniform_proofs}

\begin{proof}[Proof of \Cref{thm:grid_uniform}]
Fix a nonempty $A$, $\mathbf t\in\{\mb,\mh\}$ and $s$. By \eqref{eq:subset_bandwidth_orders}, the fixed-model versions of Conditions~\ref{con:B1}--\ref{con:B2} hold. Proposition~\ref{prop:empirical_prelocalization} puts the empirical minimizer within distance $\eta_\star$ of its population local target with probability tending to one. Let $\rho_n=(\log n/(nv_{\mathbf t}(A)))^{1/2}$ and $r_n=\rho_n^{1/(\beta_\oplus-1)}$. Lemma~\ref{lem:complete_modulus} provides $C_n=O_\P(1)$ such that all its radius-normalized bounds are at most $C_n\rho_n$.

By Condition~\ref{con:M3}(E), empirical optimality holds simultaneously at all valid predictor points almost surely. On the localization and design events, write
\begin{align*}
d=d_\M(\widehat m_{\mathbf t_A,s}(\mx_A),\widetilde m_{\mathbf t_A,s}(\mx_A)).
\end{align*}
If $d\le r_n$, the desired bound already holds. Otherwise choose $j$ with $\delta_{j-1,n}<d\le\delta_{j,n}$, so $\delta_{j,n}\le2d$. Empirical optimality and \eqref{eq:localized_margin} imply
\begin{align*}
c_\oplus d^{\beta_\oplus}&\le\widetilde M(\mx_A,\widehat m)-\widetilde M(\mx_A,\widetilde m)\\
&\le|\mathbb Z_n(\mx_A,\widehat m)|\le C_n\delta_{j,n}\rho_n\le2C_nd\rho_n.
\end{align*}
Division by $d$ proves a uniform $O_\P(\rho_n^{1/(\beta_\oplus-1)})$ bound for the empirical-to-population-local distance. Lemma~\ref{lem:localized_bias} and the triangle inequality give
\begin{align*}
\sup_{\mx_A}d_\M(\widehat m_{\mathbf t_A,s}(\mx_A),m_{\oplus,A}(\mx_A))=O_\P\left(b_{\mathbf t}(A)^{1/(\beta_\oplus-1)}+\left(\frac{\log n}{nv_{\mathbf t}(A)}\right)^{1/[2(\beta_\oplus-1)]}\right).
\end{align*}
All exceptional-event probabilities tend to zero, including the design failure on which the empirical mean fallback is used.

For $A=\varnothing$, the population local target is $m_0$, its bias is zero and the weights are one. Condition~\ref{con:M1}(E) supplies the single empirical minimizer. The finite-net proof of Proposition~\ref{prop:empirical_prelocalization} applies with no predictor index. Lemma~\ref{lem:moving_anchor_entropy} applies around $\cK_\varnothing=\{m_0\}$. Repeating the increment calculation with $v=1$ and using \eqref{eq:margin} in place of \eqref{eq:localized_margin} gives the required rate $O_\P((\log n/n)^{1/[2(\beta_\oplus-1)]})$.

For each fixed $A$ the ratio of its error to its own $r_{n,\mathbf t}(A)$ is bounded in probability. Given $\epsilon>0$, choose a probability constant for each of the finitely many factor sets with error probability at most $\epsilon/2^p$, and take their largest constant and sample-size threshold. A union bound proves \eqref{eq:grid_uniform_rate}, and gives the equivalent formulation with the model-specific rate stated separately for each $A$.
\end{proof}

\section{Factor selection and post-selection estimation}\label{supp:lofo}
\setcounter{equation}{0}

\subsection{Identification of the active factors}

\begin{proposition}\label{prop:lofo}
Suppose Conditions~\ref{con:D1} and \ref{con:M1}--\ref{con:M2} hold. If $A\subseteq B\subseteq[p]$, then $R(B)\le R(A)$. There is a unique smallest set $A_0$ under inclusion satisfying $R(A_0)=R([p])$, and
\begin{align*}
R(A)=R([p])\quad\Longleftrightarrow\quad A_0\subseteq A.
\end{align*}
Consequently, $\Delta_\ell>0$ if and only if $\ell\in A_0$. If $A_0\ne\varnothing$, then $\min_{\ell\in A_0}\Delta_\ell>0$.
\end{proposition}

\begin{proof}
Conditional minimization and Condition~\ref{con:M1} give
\begin{align*}
0\le R(A)\le\E[d_\M^2(Y,m_{\oplus,\varnothing})]<\infty.
\end{align*}
If $A\subseteq B$, the prediction $m_{\oplus,A}(\mX_A)$ is measurable with respect to $\mX_B$. Thus
\begin{align*}
\E[d_\M^2(Y,m_{\oplus,B}(\mX_B))\mid\mX_B]\le\E[d_\M^2(Y,m_{\oplus,A}(\mX_A))\mid\mX_B].
\end{align*}
Taking expectations proves monotonicity. If $R(A)=R([p])$, the nonnegative random variable
\begin{align*}
M_\oplus(\mX,m_{\oplus,A}(\mX_A))-M_\oplus(\mX,m_\oplus(\mX))
\end{align*}
has expectation zero. Uniqueness of the full conditional mean therefore implies
\begin{align}\label{eq:risk_equality_targets}
m_{\oplus,A}(\mX_A)=m_\oplus(\mX)\quad\text{almost surely}.
\end{align}

Suppose $R(A)=R(B)=R([p])$. Write $C=A\cap B$, $U=A\setminus B$ and $V=B\setminus A$. Let $\nu_J$ denote normalized product surface measure on the factors in $J$, with unit mass for the empty product. Condition~\ref{con:D1} makes $P_{\mX}$ and $\nu_{[p]}$ equivalent. By \eqref{eq:risk_equality_targets}, and after integration over the factors outside $A\cup B$,
\begin{align*}
m_{\oplus,A}(\mx_C,\mx_U)=m_{\oplus,B}(\mx_C,\mx_V)
\end{align*}
for $\nu_C\otimes\nu_U\otimes\nu_V$-almost every triple, with each displayed argument placed in its original factor order. Fubini's theorem supplies a fixed $\mathbf v_0\in\Sph^{\md_V}$ for which the equality holds for $\nu_C\otimes\nu_U$-almost every pair. The measurable function $q_C(\mx_C):=m_{\oplus,B}(\mx_C,\mathbf v_0)$ therefore satisfies $q_C(\mX_C)=m_\oplus(\mX)$ almost surely. This reasoning also covers empty $C,U$ or $V$. By conditional minimization,
\begin{align*}
R([p])\le R(C)\le\E[d_\M^2(Y,q_C(\mX_C))]=R([p]).
\end{align*}
The finite family $\{A\subseteq[p]:R(A)=R([p])\}$ is consequently closed under intersections. Its intersection belongs to the family and is contained in every other member. This gives the unique smallest set $A_0$. Monotonicity gives the converse when $A_0\subseteq A$. Applying this equivalence to $A=[p]\setminus\{\ell\}$ proves the gap characterization. The positive minimum follows from the fixed finite number of active factors. Independence between predictor factors was not used.
\end{proof}

\subsection{Cross-fitted risk gaps}\label{supp:selection}

For this subsection write $m_0=m_{\oplus,\varnothing}$ and let $\ell_A(Y,\mX):=d_\M^2(Y,m_{\oplus,A}(\mX_A))$, so $P\ell_A=R(A)$. The candidate sets are $\mathfrak C_p$ from Section~\ref{sec:selection_theory}.

\begin{lemma}\label{lem:gap_approximation}
Suppose the conditions of \Cref{thm:wasserstein} or \Cref{thm:grid_uniform} hold for the foldwise LC estimators with bandwidths $\mb_A$. Then
\begin{align*}
\max_{\ell\in[p]}|\widehat\Delta_\ell-\Delta_\ell|=O_\P(r_{n,\mb}([p])+n^{-1/2}).
\end{align*}
\end{lemma}

\begin{proof}
Let $N_k=n-|I_k|$. The fixed folds have $N_k\asymp n$. When applying the estimation theorem to a training sample, its bandwidth is the prescribed $\mb(n)$, not $\mb(N_k)$. Replacing $n$ by $N_k$ in its probability bounds only changes the orders by constants because $\log N_k\asymp\log n$. Its response maximum is bounded by the full-sample maximum. The bandwidth and tail conditions therefore hold on every training sample. Since $K$ and $p$ are fixed, the model-specific LC bounds hold simultaneously. By \eqref{eq:subset_bandwidth_orders},
\begin{align*}
R_n:=\max_{A\in\mathfrak C_p}\max_{k\in[K]}\sup_{\mx_A}d_\M(\widehat m_{\mb_A,0}^{(-k)}(\mx_A),m_{\oplus,A}(\mx_A))=O_\P(r_{n,\mb}([p]))=o_\P(1).
\end{align*}
For general responses the target radii are bounded by Lemma~\ref{lem:target_continuity}, and for Wasserstein responses they are bounded by the diameter of the response space. Choose $C_m\ge1$ such that all target radii from $m_0=m_{\oplus,\varnothing}$ are at most $C_m$. By \eqref{eq:metric_increment}, the fitted-loss error is bounded pathwise by
\begin{align*}
\max_{A\in\mathfrak C_p}|\widehat R_{\mathrm{CF}}(A)-P_n\ell_A|\le2R_n\left(P_nd_\M(Y,m_0)+C_m\right)+R_n^2=O_\P(r_{n,\mb}([p])).
\end{align*}
The law of large numbers applies to $d_\M(Y,m_0)$ by the second-moment condition. No independence of the fitted functions and validation observations is needed for this comparison.

For a fixed $\ell$, put $u=m_{\oplus,[p]\setminus\{\ell\}}(\mX_{[p]\setminus\{\ell\}})$ and $w=m_\oplus(\mX)$. Then
\begin{align*}
|\ell_{[p]\setminus\{\ell\}}-\ell_{[p]}|&\le d_\M(u,w)(d_\M(Y,u)+d_\M(Y,w))\\
&\le4C_m(d_\M(Y,m_0)+C_m)\le4C_m(1+C_m)D_\M(Y).
\end{align*}
Thus each paired loss difference has finite variance under the existing second-moment bound. Chebyshev's inequality and the fixed number $p$ give
\begin{align*}
\max_{\ell\in[p]}|(P_n-P)(\ell_{[p]\setminus\{\ell\}}-\ell_{[p]})|=O_\P(n^{-1/2}).
\end{align*}
Finally,
\begin{align*}
\widehat\Delta_\ell-\Delta_\ell&=(P_n-P)(\ell_{[p]\setminus\{\ell\}}-\ell_{[p]})\\
&\quad+\widehat R_{\mathrm{CF}}([p]\setminus\{\ell\})-P_n\ell_{[p]\setminus\{\ell\}}-\widehat R_{\mathrm{CF}}([p])+P_n\ell_{[p]}.
\end{align*}
The two bounds prove the conclusion. The paired evaluation on the same observations is what permits the second-moment argument.
\end{proof}

If the additional moment $\E[d_\M^4(Y,m_0)]<\infty$ holds, take the same uniform target-radius bound $C_m$ as in the preceding proof. Then $\ell_A^2\le8d_\M^4(Y,m_0)+8C_m^4$ has finite expectation. Chebyshev's inequality gives $\max_{A\in\mathfrak C_p}|(P_n-P)\ell_A|=O_\P(n^{-1/2})$, and the fitted-loss calculation in Lemma~\ref{lem:gap_approximation} yields $\max_{A\in\mathfrak C_p}|\widehat R_{\mathrm{CF}}(A)-R(A)|=O_\P(e_n)$. This additional conclusion is not used to prove selection consistency.

\begin{proof}[Proof of \Cref{thm:selection}]
Proposition~\ref{prop:lofo} gives the unique active set and identifies its positive gaps. In the general metric-space branch its target conditions are included in the hypotheses. In the Wasserstein branch Lemmas~\ref{lem:quantile_projection} and \ref{lem:quantile_bias} give the required existence, uniqueness, continuity and second moments.

Write $H(y,z)=d_\M^2(y,z)/2$. The triangle inequality and the finite second moment give $\E[H(Y,Y')]<\infty$. Also,
\begin{align*}
\widehat V_Y=\binom n2^{-1}\sum_{1\le i<j\le n}H(Y^{(i)},Y^{(j)}).
\end{align*}
For completeness, set $H_T=\min(H,T)$ and let $U_{n,T}$ be the corresponding U-statistic. In the variance of $U_{n,T}$, terms from disjoint observation pairs have covariance zero. There are $O(n^3)$ overlapping pairs, each bounded by $T^2$, and $\binom n2^2$ normalization. Hence $\operatorname{Var}(U_{n,T})\le CT^2/n$. For fixed $T$, $U_{n,T}-\E[H_T(Y,Y')]\to0$ in probability, while
\begin{align*}
\E[|\widehat V_Y-U_{n,T}|]=\E[H(Y,Y')-H_T(Y,Y')]\longrightarrow0
\end{align*}
as $T\to\infty$, uniformly in $n$. Markov's inequality and dominated convergence prove $\widehat V_Y\to V_Y$ in probability without a fourth moment.

It follows that
\begin{align*}
\widehat\tau_n/a_n=c\widehat V_Y=cV_Y+o_\P(1),\quad \widehat\tau_n=o_\P(1).
\end{align*}
By Lemma~\ref{lem:gap_approximation}, $e_n=o(a_n)$ and $cV_Y>0$,
\begin{align*}
\frac{\max_{\ell\in[p]}|\widehat\Delta_\ell-\Delta_\ell|}{\widehat\tau_n}=o_\P(1)
\end{align*}
on the event $\widehat\tau_n>0$, whose probability tends to one. If $A_0=\varnothing$, every population gap is zero and no factor is selected whenever the maximum gap error is less than $\widehat\tau_n$.

If $A_0\ne\varnothing$, put $\Delta_{\min}=\min_{\ell\in A_0}\Delta_\ell>0$ within this proof. With probability tending to one,
\begin{align*}
\max_{\ell\in[p]}|\widehat\Delta_\ell-\Delta_\ell|<\widehat\tau_n,\quad 0<2\widehat\tau_n<\Delta_{\min}.
\end{align*}
Every inactive gap is then less than $\widehat\tau_n$, and every active estimated gap is at least $\Delta_{\min}-\widehat\tau_n>\widehat\tau_n$. Hence $\P(\widehat A_{\mb}=A_0)\to1$.
\end{proof}

\subsection{Selected estimators}

\begin{proof}[Proof of \Cref{thm:post_selection}]
By \eqref{eq:risk_equality_targets}, $m_{\oplus,A_0}(\mX_{A_0})=m_\oplus(\mX)$ almost surely. For general responses, Lemma~\ref{lem:target_continuity} makes both targets continuous. For Wasserstein responses the same continuity follows from Lemma~\ref{lem:quantile_bias}. Thus the function
\begin{align*}
\mx\longmapsto d_\M(m_{\oplus,A_0}(\mx_{A_0}),m_\oplus(\mx))
\end{align*}
is continuous. If it were positive at some point, it would be positive on a nonempty open set. Condition~\ref{con:D1} assigns positive probability to that set, contradicting almost-sure equality. Consequently the targets agree at every $\mx\in\Sph^{\md}$.

On $\widehat A_{\mb}=A_0$, the selected estimator is the estimator on $A_0$ with bandwidth $\mh_{A_0}$. Since the projection of the full predictor space onto $\Sph^{\md_{A_0}}$ is onto,
\begin{align*}
\sup_{\mx}d_\M(\widehat m_{\mb,\mh,s}(\mx),m_\oplus(\mx))=\sup_{\mx_{A_0}}d_\M(\widehat m_{\mh_{A_0},s}(\mx_{A_0}),m_{\oplus,A_0}(\mx_{A_0})).
\end{align*}
For every $C>0$, the probability that the left-hand side exceeds $Cr_{n,\mh}(A_0)$ is at most
\begin{align*}
\P(\widehat A_{\mb}\ne A_0)+\P\left(\sup_{\mx_{A_0}}d_\M(\widehat m_{\mh_{A_0},s}(\mx_{A_0}),m_{\oplus,A_0}(\mx_{A_0}))>Cr_{n,\mh}(A_0)\right).
\end{align*}
The first term tends to zero by \Cref{thm:selection}. The second gives boundedness in probability by the appropriate uniform theorem. This proves the asserted rate for either $s=0$ or $s=1$. The argument does not require independence between selection and the final fit. If $A_0=\varnothing$, use the empirical mean and the empty-model rate.

For the final consequence, set $h_n=(\log n/n)^{1/(|\md_{A_0}|+4)}$. The assumed active bandwidth order gives $b_{\mh}(A_0)\asymp h_n^2$ and $v_{\mh}(A_0)\asymp h_n^{|\md_{A_0}|}$. With $\beta_\oplus=2$,
\begin{align*}
r_{n,\mh}(A_0)\asymp h_n^2+\sqrt{\frac{\log n}{nh_n^{|\md_{A_0}|}}}\asymp(\log n/n)^{2/(|\md_{A_0}|+4)}.
\end{align*}
This substitution is made under all the hypotheses of the theorem. In particular, the active bandwidth order alone does not imply the full-product bandwidth and tail conditions for the other factors. The conclusion is a stochastic-order statement, not a bound on expected prediction loss.
\end{proof}

\subsection{Bandwidth selection from finite grids}\label{supp:finite_grids}

This subsection permits bandwidths to be selected separately for different models and training folds. For each relevant $A$, let $\cH_{A,n}^{\mathrm{sel}}$ and $\cH_{A,n}^{\mathrm{fin}}$ be deterministic nonempty finite sets of positive vectors indexed by $A$. Their cardinalities are bounded by a constant independent of $n$. For a candidate $\mathbf u=(u_\ell:\ell\in A)$, use
\begin{align*}
v_{\mathbf u}(A)=\prod_{\ell\in A}u_\ell^{d_\ell},\quad b_{\mathbf u}(A)=\sum_{\ell\in A}u_\ell^2,\quad r_{n,\mathbf u}(A)=b_{\mathbf u}(A)^{1/(\beta_\oplus-1)}+\left(\frac{\log n}{nv_{\mathbf u}(A)}\right)^{1/[2(\beta_\oplus-1)]}.
\end{align*}
These definitions require only the retained coordinates. For $A=\varnothing$, each grid contains the empty vector and the empirical mean convention applies. Put
\begin{align*}
r_n^{\mathrm{sel}}:=\max_{A\in\mathfrak C_p}\max_{\mathbf u\in\cH_{A,n}^{\mathrm{sel}}}r_{n,\mathbf u}(A),\quad r_n^{\mathrm{fin}}(A):=\max_{\mathbf u\in\cH_{A,n}^{\mathrm{fin}}}r_{n,\mathbf u}(A).
\end{align*}

\begin{proposition}\label{prop:finite_grid_extension}
Suppose the distribution, kernel and response-space assumptions of the relevant fixed-model uniform theorem hold. In (ii), impose the estimation hypotheses for the LC estimators used in selection and for the chosen final estimator. For each nonempty relevant $A$ and each grid $\cH_{A,n}$ used at either stage, suppose
\begin{align*}
\max_{\mathbf u\in\cH_{A,n}}\max_{\ell\in A}u_\ell\to0,\quad \frac n{\log n}\min_{\mathbf u\in\cH_{A,n}}v_{\mathbf u}(A)\to\infty.
\end{align*}
For the general metric-space branch, also suppose
\begin{align*}
\max_{1\le i\le n}D_\M(Y^{(i)})=O_\P\left(\min_{\mathbf u\in\cH_{A,n}}\sqrt{\frac{nv_{\mathbf u}(A)}{\log n}}\right).
\end{align*}
The following conclusions hold.
\begin{enumerate}[label=(\roman*),leftmargin=2.2em]
\item For fixed $s$, every measurable bandwidth choice from the grid satisfies
\begin{align*}
\sup_{\mx_A}d_\M(\widehat m_{\widehat{\mathbf u}_A,s}(\mx_A),m_{\oplus,A}(\mx_A))=O_\P\left(\max_{\mathbf u\in\cH_{A,n}}r_{n,\mathbf u}(A)\right).
\end{align*}
The assertion also holds on each of the fixed training folds.
\item Choose each selection bandwidth measurably from $\cH_{A,n}^{\mathrm{sel}}$ using only the corresponding training observations, and compute the same cross-fitted LOFO gaps and threshold. If $0<V_Y<\infty$, $a_n\to0$ and $r_n^{\mathrm{sel}}+n^{-1/2}=o(a_n)$, the resulting selected set equals $A_0$ with probability tending to one. If the final bandwidth is chosen measurably from $\cH_{\widehat A,n}^{\mathrm{fin}}$, the selected estimator with index $s$ has uniform error $O_\P(r_n^{\mathrm{fin}}(A_0))$.
\end{enumerate}
\end{proposition}

\begin{proof}
Enumerate each grid deterministically and pad the list, if needed, by repeating its first element so that it has a fixed length $J_A$ bounded uniformly in $n$. Every enumerated candidate is a deterministic sequence satisfying the fixed-model bandwidth and, when needed, tail conditions. The proofs of the fixed-model uniform theorems use precisely these local conditions, rather than values of omitted components. They therefore give an individual rate for each candidate sequence. Condition~\ref{con:M3}(E) supplies simultaneous empirical minimizers for each candidate and training fold almost surely. There are finitely many candidates, factor sets and folds, so a finite union makes these bounds simultaneous.

If $\widehat{\mathbf u}_A$ belongs to the grid, then pathwise
\begin{align*}
\sup_{\mx_A}d_\M(\widehat m_{\widehat{\mathbf u}_A,s}(\mx_A),m_{\oplus,A}(\mx_A))\le\max_{\mathbf u\in\cH_{A,n}}\sup_{\mx_A}d_\M(\widehat m_{\mathbf u,s}(\mx_A),m_{\oplus,A}(\mx_A)).
\end{align*}
The simultaneous candidate bounds prove (i). Sample sizes in the training folds are comparable to $n$, and their response maxima are bounded by the full-sample maximum, exactly as in the proof of Lemma~\ref{lem:gap_approximation}.

For (ii), the maximum training-fold LC error is $O_\P(r_n^{\mathrm{sel}})$. Substituting this bound in the pathwise fitted-loss calculation of Lemma~\ref{lem:gap_approximation} gives gap error $O_\P(r_n^{\mathrm{sel}}+n^{-1/2})$. The proof of \Cref{thm:selection} then applies without change to the threshold comparison. On the event of correct selection, (i) on the oracle-model final grid bounds the error by $O_\P(r_n^{\mathrm{fin}}(A_0))$. The probability decomposition in the proof of \Cref{thm:post_selection} completes the argument.
\end{proof}

The proposition allows different selected bandwidths across models and folds, and different final choices for LC and LL. It does not require these choices to be subvectors of one common random full-model vector. It also does not prove that cross-validation selects an optimal bandwidth. The active-dimension rate follows when all eligible final candidates on $A_0$ have that rate, or when a separate result establishes the rate of the selected candidate. The presence of one optimal candidate is not sufficient for this conclusion. Fixed numerical grids at the reported sample sizes are distinct from the shrinking grid sequences required by the proposition. No growing-cardinality extension is asserted here.

\clearpage
\section{Additional details for simulation analysis}\label{supp:c_sensitivity}

For the proposed and Euclidean estimators, we minimize the cross validation criterion by selecting one bandwidth at a time while holding the others fixed. For each bandwidth, we compare its $20$ candidate values and retain the one giving the smallest cross validation score. We repeat this procedure over all bandwidths until the selected values no longer change or the maximum number of iterations is reached. The proposed estimators use at most five iterations for factor selection and ten for final estimation, starting from the middle grid values. The Euclidean estimators use $23$ fixed starting vectors to reduce sensitivity to initialization, with at most ten iterations per start. The bandwidth vector with the smallest cross validation score is retained. The Riemannian estimators use grid search.

During cross validation, bandwidth candidates for local linear estimation are excluded if the reciprocal condition number of the corresponding coordinate design matrix is below $10^{-6}$. For spherical responses, candidates with failed numerical minimization are also excluded. At prediction points, unstable local linear estimates or failed local linear minimizations are replaced by the corresponding local constant estimates. The training Fr\'echet mean is used when a compactly supported kernel has empty local support or, for spherical responses, local constant minimization fails.

We assess the sensitivity of factor selection to eight equally spaced multipliers $c\in\{0.0025j:j=1,\ldots,8\}$ using the saved cross-fitted LOFO risk gaps.

\subsection{Case of scalar responses}\label{supp:scalar_sensitivity}

Table~\ref{tab:supp_scalar_c_sensitivity} reports exact factor-selection frequencies for scalar responses. The multiplier $0.005$ was selected using $200$ Monte Carlo replications for each of the six strong signal settings with $n\in\{500,1000,2000\}$ and $\sigma\in\{0.25,1\}$. It was the smallest candidate with mean false inclusion frequency at most $0.05$ and $X_2$ inclusion frequency at least $0.90$ in every calibration setting.

\begin{table}[!htbp]
\centering
\caption{Exact factor-selection frequencies for scalar responses over eight multiplier candidates. Each entry pools $200$ replications for each $n\in\{500,1000,2000\}$ and $\sigma\in\{0.25,1\}$. The calibrated multiplier is $c=0.005$.}
\label{tab:supp_scalar_c_sensitivity}
\begingroup
\setlength{\tabcolsep}{4pt}
\begin{tabular}{l*{8}{c}}
\toprule
Signal & 0.0025 & 0.005 & 0.0075 & 0.01 & 0.0125 & 0.015 & 0.0175 & 0.02\\
\midrule
Strong & 0.9542 & 0.9867 & 0.9950 & 0.9975 & 0.9992 & 1.0000 & 1.0000 & 1.0000\\
Weak & 0.9475 & 0.9833 & 0.9908 & 0.9967 & 0.9975 & 0.9983 & 0.9975 & 0.9950\\
\bottomrule
\end{tabular}
\endgroup
\end{table}

\subsection{Case of spherical responses}\label{supp:sphere_sensitivity}

Table~\ref{tab:supp_sphere_c_sensitivity} reports exact factor-selection frequencies for spherical responses. At $c=0.02$, the minimum $X_2$ inclusion frequency across the $12$ settings is $0.900$.

\begin{table}[!htbp]
\centering
\caption{Exact factor-selection frequencies for spherical responses over eight multiplier candidates. Each entry pools $200$ replications for each $n\in\{500,1000,2000\}$ and $\kappa\in\{400,100\}$. The calibrated multiplier is $c=0.005$.}
\label{tab:supp_sphere_c_sensitivity}
\begingroup
\setlength{\tabcolsep}{4pt}
\begin{tabular}{l*{8}{c}}
\toprule
Signal & 0.0025 & 0.005 & 0.0075 & 0.01 & 0.0125 & 0.015 & 0.0175 & 0.02\\
\midrule
Strong & 0.9758 & 0.9950 & 1.0000 & 1.0000 & 1.0000 & 1.0000 & 1.0000 & 1.0000\\
Weak & 0.9558 & 0.9883 & 0.9983 & 0.9992 & 0.9983 & 0.9933 & 0.9900 & 0.9833\\
\bottomrule
\end{tabular}
\endgroup
\end{table}

\end{document}